\documentclass[11pt]{amsart}

\usepackage{amssymb, amsmath, amscd, amsthm, graphicx, psfrag,rotating}
\usepackage[matrix,arrow,curve]{xy}
\xyoption{dvips}              
\usepackage{wasysym}

\newcommand{\R}{{\mathbb R}}

\newcommand{\calO}{{\mathcal{O}}}
\newcommand{\calI}{{\mathcal{I}}}

\newcommand{\calP}{{\mathcal{P}}}

\newcommand{\calY}{{\mathcal{Y}}}

\newcommand\rth{\refstepcounter{equation}}
\newcommand\numb{\rth{\rm \theequation}}
\numberwithin{equation}{section}

\DeclareMathOperator*{\holim}{holim}
\DeclareMathOperator{\hoTot}{\Tot}
\DeclareMathOperator{\Tot}{Tot}
\DeclareMathOperator{\sTot}{sTot}
\DeclareMathOperator{\hosTot}{\sTot}

\newcommand{\unit}{{\mathbf{1}}}

\newcommand{\Assoc}{{{\mathcal A}ssoc}}

\newcommand{\Top}{\mathrm{Top}}

\newcommand{\calC}{{\mathcal{C}}}
\newcommand{\calX}{{\mathcal{X}}}

\newcommand{\id}{\mathrm{id}}

\newcommand{\calK}{{\mathcal{K}}}
\newcommand{\calQ}{{\mathcal{Q}}}

\newcommand{\sfX}{{\mathsf{X}}}
\newcommand{\sfY}{{\mathsf{Y}}}

\newcommand{\Ebar}{\overline{Emb}}

\newcommand{\Operad}{\operatorname{Operads}}
\newcommand{\hoOperad}{\widetilde{\Operad}}
\newcommand{\Pentagon}{\text{\Large{{\pentagon}}}}
\newcommand{\hoSquare}{\widetilde{\square}}
\newcommand{\hoTriangle}{\widetilde{\triangle}}
\newcommand{\hoPentagon}{\widetilde{\Pentagon}}
\newcommand{\WBimod}{\operatorname{Wbimod}}
\newcommand{\Bimod}{\operatorname{Bimod}}
\newcommand{\hoBimod}{\widetilde{\Bimod}}

\newcommand{\pushout}{\text{\makebox[0pt][r]{\raisebox{13pt}[1pt]{$\lrcorner$\hspace{2pt}}}}}

\theoremstyle{plain}
\newtheorem{thm}{Theorem}[section]
\newtheorem{theorem}[thm]{Theorem}

\newtheorem{proposition}[thm]{Proposition}
\newtheorem{lemma}[thm]{Lemma}

\newtheorem{cor-def}[thm]{Corollary-Definition}

\newtheorem*{theorem*}{Theorem}

\theoremstyle{definition}
\newtheorem{definition}[thm]{Definition}
\newtheorem{remark}[thm]{Remark}
\newtheorem{notation}[thm]{Notation}

\theoremstyle{remark}

\newtheorem*{notation*}{Notation}

\begin{document}


\title{Delooping totalization of a multiplicative operad}



\author{Victor Turchin}
\address{Kansas State University, USA.}
\email{turchin@ksu.edu} \urladdr{http://www.math.ksu.edu/\~{}turchin/}
\thanks{The author was supported in part by NSF research grant DMS-0968046.}

%


\begin{abstract}
The paper shows that under some conditions the totalization of a cosimplicial space obtained from a multiplicative operad is a double loop space  of the space of derived morphisms from the associative operad to the operad itself.
\end{abstract}

\maketitle

\sloppy


\section{Introduction}\label{s:introduction}
Let $\calO$ be a topological non-$\Sigma$ operad endowed with a map
$
\Assoc\to\calO
$
from the  associative operad. In such case $\calO$ is called {\it multiplicative} operad. One can also say that $\calO$ is a non-$\Sigma$ operad in the category of
 based topological spaces. This structure endows $\{\calO(n),\, n\geq 0\}$ with a structure of a cosimplicial space. McClure and Smith~\cite[Theorem~3.3]{McSmith} showed
 that the homotopy totalization $\hoTot\calO(\bullet)$ of such cosimplicial space admits an action of an operad equivalent to the operad of little squares. In this
 paper we show that in the case $\calO(0)\simeq\calO(1)\simeq *$ the space $\hoTot\calO(\bullet)$ is homotopy equivalent to the double loop space of the space of derived
 morphisms of operads $\hoOperad(\Assoc,\calO)$. (Here and below for a category $\calC$ we denote by $\calC(a,b)$ the space of morphisms between two
 objects $a$ and $b$ in $\calC$; and by $\widetilde{\calC}(a,b)$ we denote the space of derived morphisms.) The same result was obtained a few month earlier by Dwyer and
 Hess~\cite{DwHe}. Our proof seems to be more geometrical and less relying on the techniques of the homotopy theory. We construct an explicit cofibrant
 model $\hoPentagon$ of the operad $\Assoc$, and a homotopy equivalence
$$
\Omega^2\Operad(\hoPentagon,\calO)\stackrel{\simeq}{\longrightarrow}\hoTot\calO(\bullet).
\eqno(\numb)\label{eq:11}
$$

The construction involves several steps. First $\hoTot\calO(\bullet)$ can also be expressed as the space of (derived) morphisms in the category $\underset{\Assoc}{\WBimod}$ of {\it weak bimodules} over the operad $\Assoc$:\footnote{The category  $\underset{\Assoc}{\WBimod}$ is equivalent to the category
of cosimplicial spaces, see Lemma~\ref{l:cosimplicial_semicos}.}
$$
\hoTot\calO(\bullet)\simeq \widetilde{\underset{\Assoc}{\WBimod}}(\Assoc,\calO)\simeq \underset{\Assoc}{\WBimod}(\hoTriangle,\calO),
\eqno(\numb)\label{eq:12}
$$
where $\hoTriangle$ is a cofibrant model of $\Assoc$ in the category $\underset{\Assoc}{\WBimod}$, see Section~\ref{ss:deg_wb}. Secondary, one shows that $\underset{\Assoc}{\WBimod}(\hoTriangle,\calO)$ is a loop space of the space of derived morphisms in the category $\underset{\Assoc}{\Bimod}$ of bimodules over $\Assoc$:
$$
\underset{\Assoc}{\WBimod}(\hoTriangle,\calO)\simeq\Omega\,\underset{\Assoc}{\hoBimod}(\Assoc,\calO),
\eqno(\numb)\label{eq:13}
$$
see Section~\ref{s:first_deloop_tilde}. This homotopy equivalence is obtained by constructing an explicit cofibrant model $\hoSquare$ of $\Assoc$ in the category $\underset{\Assoc}{\Bimod}$, and a map
$$
\Omega\,\underset{\Assoc}{\Bimod}(\hoSquare,\calO)\stackrel{\simeq}{\longrightarrow}\underset{\Assoc}{\WBimod}(\hoTriangle,\calO).
\eqno(\numb)\label{eq:14}
$$
For the second delooping of $\hoTot\,\calO(\bullet)$, we construct an explicit cofibrant model $\hoPentagon$ of $\Assoc$ in the category of non-$\Sigma$ operads, and a  homotopy equivalence
$$
\Omega\,\Operad(\hoPentagon,\calO)\stackrel{\simeq}{\longrightarrow}\underset{\Assoc}{\Bimod}(\hoSquare,\calO),
\eqno(\numb)\label{eq:15}
$$
see Section~\ref{s:second_deloop_tilde}.

The construction of the cofibrant models $\hoTriangle$, $\hoSquare$, $\hoPentagon$ of $\Assoc$ in the categories $\underset{\Assoc}{\WBimod}$, $\underset{\Assoc}{\Bimod}$, $\Operad$ respectively is quite complicated (we actually do it only in the second part of the paper). For example, $\hoTriangle$ is the same as a contractible in each degree cofibrant cosimplicial space (cofibrant in the sense the corresponding functor from the simplicial indexing category $\Delta$ to spaces is cofibrant or free in the sense of Farjoun~\cite{Farjoun}). Because the degeneracies should also act freely, all the components of $\hoTriangle$ are infinite dimensional $CW$-complexes. To make life easier in Part~\ref{part1} we make all the above constructions ignoring degeneracies.  It is well known that the homotopy totalization of a cosimplicial space is weakly equivalent to its semicosimplicial homotopy totalization~\cite[Lemma~3.8]{DrDw}:
$$
\hoTot\calO(\bullet)\simeq\hosTot\calO(\bullet).
\eqno(\numb)\label{eq:16}
$$
In our terms a semicosimplicial space is the same thing as a weak bimodule over the operad $\Assoc_{>0}$ of associative non-unitary monoids, see Lemma~\ref{l:cosimplicial_semicos}, and the homotopy equivalence~\eqref{eq:16} can be rewritten as
$$
\underset{\Assoc}{\widetilde{\WBimod}}(\Assoc,\calO)\simeq \underset{\Assoc_{>0}}{\widetilde{\WBimod}}(\Assoc,\calO).
\eqno(\numb)\label{eq:17}
$$
In Part~\ref{part1} we deloop the right-hand side. We construct a cofibrant model $\triangle$ of $\Assoc$ in the category $\underset{\Assoc_{>0}}{\WBimod}$ of weak bimodules over $\Assoc_{>0}$, a cofibrant model $\square$ of $\Assoc_{>0}$ in the category $\underset{\Assoc_{>0}}{\Bimod}$ of bimodules over $\Assoc_{>0}$, a cofibrant model $\Pentagon$  of $\Assoc_{>0}$ in the category of operads. We also construct homotopy equivalences
$$
\Omega\,\underset{\Assoc_{>0}}{\Bimod}(\square,\calO)\stackrel{\simeq}{\longrightarrow}\underset{\Assoc_{>0}}{\WBimod}(\triangle,\calO),
\eqno(\numb)\label{eq:18}
$$
and
$$
\Omega\,\Operad(\Pentagon,\calO)\stackrel{\simeq}{\longrightarrow}\underset{\Assoc_{>0}}{\Bimod}(\square,\calO).
\eqno(\numb)\label{eq:19}
$$
see Sections~\ref{s:first_deloop},~\ref{s:second_deloop}. These two homotopy equivalences show that the space of derived morphisms $\widetilde{\Operad}(\Assoc_{>0},\calO)$ is a double delooping of $\hosTot\calO(\bullet)\simeq\hoTot\calO(\bullet)$. We mention that $\triangle$, $\square$, and $\Pentagon$ are polytopes in each degree. For example, the $n$-th component $\triangle(n)$, $n\geq 0$, is an $n$-simplex. The component $\square(n)$, $n\geq 1$, is an $(n-1)$-cube. The component $\Pentagon(n)$, $n\geq 2$, is an $(n-2)$-dimensional Stasheff polytope. Notice that the dimension of $\triangle(n)$, $\square(n)$, and $\Pentagon(n)$  is $n$, $n-1$, and $n-2$ respectively. The fact that each time we loose one degree of freedom explains why each next space of derived morphisms is a delooping of the previous one.

The Part~\ref{part2} is more technical since $\hoTriangle$, $\hoSquare$, and $\hoPentagon$ have a more complicated description. Our main motivation to include this part was that the method produces deloopings of the partial homotopy totalizations $\hoTot_N\calO(\bullet)$, $N\geq 0$, in terms of the derived morphisms of truncated bimodules or truncated operads. It is well known that the partial semicosimplicial totalizations do not \lq\lq approximate" in the right sense the totalization. To be precise the (homotopy) limit of the tower
$$
\hosTot_0\calO(\bullet)\leftarrow\hosTot_1\calO(\bullet)\leftarrow\hosTot_2\calO(\bullet)\leftarrow\hosTot_3\calO(\bullet)\leftarrow\ldots
\eqno(\numb)\label{eq:110}
$$
is $\hosTot\calO(\bullet)$, but the maps between the stages of the tower become higher and higher connected only for very degenerate (semi)cosimplicial spaces. On the contrary if we consider the tower of partial totalizations
$$
\hoTot_0\calO(\bullet)\leftarrow\hoTot_1\calO(\bullet)\leftarrow\hoTot_2\calO(\bullet)\leftarrow\hoTot_3\calO(\bullet)\leftarrow\ldots,
\eqno(\numb)\label{eq:1_11}
$$
then under a natural convergency condition the maps in the tower do become higher and higher connected with $N$.

We show that under the condition $\calO(0)\simeq\calO(1)\simeq*$ the partial (semicosimplicial) totalizations of a multiplicative operad $\calO$ deloop as spaces of derived morphisms of truncated operads:
$$
\hosTot_N\calO(\bullet)\simeq\Omega^2\widetilde{\Operad}_N(\Assoc_{>0}|_N,\calO|_N),
\eqno(\numb)\label{eq:deloop_part_sTot}
$$
$$
\hoTot_N\calO(\bullet)\simeq\Omega^2\widetilde{\Operad}_N(\Assoc|_N,\calO|_N).
\eqno(\numb)\label{eq:deloop_part_Tot}
$$
To the best of our understanding the approach of Dwyer and Hess~\cite{DwHe} can also be  extended to prove~\eqref{eq:deloop_part_sTot}-\eqref{eq:deloop_part_Tot}, but still requires an additional amount of work.

We reiterate that this work is based on the following very simple observation. The partial semicosimplicial totalization
$\hosTot_N\calO(\bullet)$  is the space of tuples of maps $\bigl(f_k\colon\triangle(k)\to\calO(k)\bigr)_{k=0\ldots N}$ which are compatible in the sense the behavior of each $f_k$ on the boundary of the $k$-simplex $\triangle(k)$ is completely determined by the previous $f_{k-1}$ (since these maps should respect the truncated semicosimplicial structure). This implies that the preimage of any point under the fibration
$$
\hosTot_N\calO(\bullet)\to \hosTot_{N-1}\calO(\bullet)
$$
is either empty or homotopy equivalent to $\Omega^N\calO(N)$. On the other hand  we can also consider the Stasheff operad $\Pentagon$ and the space of maps of $N$-truncated operads $\Operad_N(\Pentagon|_N,\calO|_N)$. The latter space consists of tuples of maps $\bigl(h_k\colon\Pentagon(k)\to\calO(k)\bigr)_{k=2\ldots N}$ compatible in the sense $h_k|_{\partial\pentagon(k)}$ is determined by $h_2,\ldots,h_{k-1}$. As a consequence the preimage of the fibration
$$
\Operad_N\left(\Pentagon|_N,\calO|_N\right)\to\Operad_{N-1}\left(\Pentagon|_{N-1},\calO|_{N-1}\right)
$$
over any point is either empty or homotopy equivalent to $\Omega^{N-2}\calO(N)$. Thus we obtain that the fibers of the tower~\ref{eq:110}
are double loop spaces of the fibers of the tower
$$
\Operad_1\left(\Pentagon|_1,\calO|_1\right)\leftarrow\Operad_2\left(\Pentagon|_2,\calO|_2\right)\leftarrow\Operad_3\left(\Pentagon|_3,\calO|_3\right)\leftarrow\ldots
$$
This naturally leads to a question whether the same  is true for the stages of the tower, which we  prove in Part~\ref{part1}.

In  Part~\ref{part2} we adjust the construction taking into account degeneracies. Because of their presence the objects $\hoTriangle$ and $\hoPentagon$ are much more complicated than $\triangle$ and $\Pentagon$. In particular in each degree these objects are no more polytopes, but infinite dimensional $CW$-complexes. However the idea of Part~\ref{part1} can still be applied as it turns out $\hoTriangle$ and $\hoPentagon$ are freely generated by almost the same countable set of cells and the generating cells of $\hoPentagon$ has dimension 2 less than the corresponding cells of $\hoTriangle$. We said almost because in this correspondence one has to exclude  the generating cells of $\hoTriangle$ of dimension~0 and~1.

The algebraic structures such as an operad or a (weak) bimodule are governed by appropriate colored operads. The necessary language of cellular cofibrant algebras over colored operads is developed in Section~\ref{s:col_operads}.

We mention that one of the motivations and examples for this work is the space $\Ebar(\R^1,\R^d)$, $d\geq 4$, of {\it long knots modulo immersions}~\cite{Sinha-OKS}, which according to Sinha~\cite{Sinha-OKS} is the homotopy totalization of the Kontsevich operad $\calK_d$:
$$
\Ebar(\R^1,\R^d)\simeq\hoTot\calK_d(\bullet).
$$
To recall the operad $\calK_d$ is weakly equivalent to the operad $\calC_d$ of little $d$-cubes. For this example the tower~\eqref{eq:1_11} is the Goodwillie-Weiss embedding tower~\cite[Theorem~1.1]{Sinha-OKS}, for which our paper gives an explicit double delooping. The idea that the space of long knots should be related to the space of derived morphisms of operads $\widetilde{\Operad}(\calC_1,\calC_d)$ is due to Kontsevich.

\section{Cofibrant models of algebras over colored operads}\label{s:col_operads}
\subsection{Colored operads}\label{ss:col_oper}
All throughout the paper when we say {\it operad} (without the adjective \lq\lq colored") we mean a non-$\Sigma$ operad. However we will also need to use the \lq\lq metalanguage" of so called {\it colored operads}. The main difference with the usual operads is that the inputs and output of an operation might have not one but many different colors. When one composes such operations one has to take into account these colors allowing to insert an output to an input only if their colors match. We refer to~\cite{BergMoerd} as an introduction to this notion. The authors of the latter paper (C.~Berger and I.~Moerdijk)  produce  a (cofibrantly generated) model structure on the category of algebras over a colored operad~\cite[Theorem~2.1]{BergMoerd}. To avoid a heavy language of the homotopy theory we will not be using this model structure,  but a reader familiar with this notion can see that this language is appropriate for our situation. The book of M.~Hoovey~\cite{Hovey} gives a good introduction to this notion.  When we say that an algebra is a cofibrant model over a colored operad we mean that it was obtained by a sequence of attachments of cells, which is a very explicit construction described below. The main properties of  such cofibrant replacements are given by Lemmas~\ref{l:extension_space},~\ref{l:fiber_equiv},~\ref{l:cof_mod},~\ref{l:cof_maps},~\ref{l:truncation},~\ref{l:trunc_maps} below. Models cofibrant in our sense are actually cofibrant in the sense of the model structure defined in~\cite{BergMoerd}.

We will be considering only colored operads whose set of colors is finite or at most countable, and whose components (spaces of operations) in all multidegrees  are discrete sets. In all our examples these sets are countable. Such operads will be called {\it discrete} colored operads. We will be also dealing with {\it $\Sigma$-free} operads which means that the permutation of entries of the same color on each component is given by a free action of the corresponding group (which is a product of symmetric groups). The algebras over such operads will be taken in the category $\Top$ of compactly generated topological spaces. The (monoidal) product in the latter category is the weak product of spaces.  The  assumption that the operad is $\Sigma$-free is nice to have when one has to deal with the cellular structure of the components of cofibrant operads. 
Our examples are the colored operads, algebras over which form one of the following categories:
\begin{itemize}
\item $\underset{\Assoc_{>0}}{\WBimod}$, $\underset{\Assoc}{\WBimod}$ weak bimodules over $\Assoc_{>0}$, $\Assoc$ respectively;
\item $\underset{\Assoc_{>0}}{\Bimod}$, $\underset{\Assoc}{\Bimod}$ bimodules over $\Assoc_{>0}$, $\Assoc$ respectively;
\item $\Operad$ -- non-$\Sigma$ operads;
\item truncated analogues of these structures:  $\underset{\Assoc_{>0}}{\WBimod}{}_N$, $\underset{\Assoc}{\WBimod}{}_N$, $\underset{\Assoc_{>0}}{\Bimod}{}_N$, $\underset{\Assoc}{\Bimod}{}_N$, $\Operad_N$, $N\geq 0$.
\end{itemize}
In the latter case the set of colors is finite ($=N+1$). The algebraic structures in question are considered in Sections~\ref{ss:wb},~\ref{ss:bimodules},~\ref{ss:trunc_operads}.

\subsection{Cofibrant models}\label{ss:cof_models}
Let $\calX$ be a $\Sigma$-free discrete colored operad with $\calI$ as a set of its colors. Let $A=\{A(i),\, i\in\calI\}$ be a collection of topological spaces. We will call such family of spaces by an $\calI$-{\it collection}. We will denote by $\calX(A)$ the free $\calX$-algebra generated by $A$. In particular $\calX(\emptyset)$ will denote the free $\calX$-algebra generated by the collection of empty sets. Notice that $\calX(\emptyset)$ might have non-empty components, since the operad $\calX$ is allowed to have operations of arity~0.

By abuse of notation $\calX$ will also denote the category of $\calX$ algebras.  The space of $\calX$-maps between $\calX$-algebras $A$ and $B$ is thus denoted by $\calX (A,B)$. The algebra $\calX(\emptyset)$ is the initial object in this category.

Let $D=\{D(i),\, i\in\calI\}$ be a collection of topological spaces where each $D(i)$ is a disjoint union of discs possibly of different dimensions. Let also
$\partial D=\{\partial D(i),\, i\in\calI\}$ be a collection of (disjoint unions of) spheres each $\partial D(i)$ being the boundary of $D(i)$. One has the inclusion of $\calX$-algebras
$$
\xymatrix{\calX(\partial D)\ar@{^{(}->}[r]^{\calX(\iota)}&\calX(D)}
$$
induced by the inclusion $\partial D\stackrel{\iota}{\hookrightarrow}D$ of $\calI$-collections. We say that an $\calX$-algebra $C$ is obtained from an $\calX$-algebra $B$ by a {\it free attachment of cells} (or simply by an {\it attachment of cells}) if they fit into the following cocartesian diagram of $\calX$-algebras:
$$
\xymatrix{
\calX(\partial D)\ar[r]^{\calX(\iota)}\ar[d]_{\calX(f)}&\calX(D)\ar[d]\\
A\ar[r]&\text{\makebox[0pt][r]{\raisebox{13pt}[1pt]{$\lrcorner$\hspace{2pt}}}}B.
}
\eqno(\numb)\label{eq:21}
$$
In the above $\calX(f)$ denotes the map of $\calX$-algebras induced by the (attaching) map $f$ of $\calI$-collections:
$$
f_i\colon \partial D(i)\to A(i),\, i\in\calI.
$$

One has the composition of the maps of $\calI$-collections
$
D\hookrightarrow\calX(D)\to B,
$
which is inclusion on the interior of each $D(i)$. The image in $B$ of each connected component of the interior of $D(i)$ will be called a {\it generating cell} in $B$ of {\it color $i$}.

\begin{definition}\label{d:cofibrant_alg}
An algebra over a $\Sigma$-free discrete colored operad $\calX$ is said {\it cellular cofibrant} (or simply {\it cofibrant}) if it is obtained from $\calX(\emptyset)$ by a possibly infinite sequence of attachments of cells.
\end{definition}

The word \lq\lq sequence" in the above definition might mean any ordinal. In all our examples the sequence will still be finite or at most countable ($={\aleph}$), moreover each attachment will be an attachment of finitely many cells of the same dimension and to the same component (thus it can be split into a sequence of attachments of only one cell). From this definition one can show that the components of any cofibrant algebra over a discrete colored operad are generalized $CW$-complexes (\lq\lq generalized" means that the attachment maps do not have to respect the dimension of the cells). Though in all our examples of such algebras the components will be honest $CW$-complexes.

As an algebra in sets any cellular cofibrant $\calX$ algebra is a free $\calX$ algebra generated by its generating cells.

It might happen that a colored operad $\calX$ has only unary non-trivial operations. In that case an $\calX$-algebra is the same thing as a functor from some discrete category whose set of objects is $\calI$, and the set of morphisms between $i$ and $j$ is given by the unary operations with input colored by $i$ and output colored by $j$. We mention that the idea of a cofibrant cellular functor is due to Dr.~Farjoun~\cite{Farjoun} and the described construction is a straightforward generalization of his work. This construction implicitly (and actually in  higher generality) appears in~\cite{BergMoerd} where the authors show that under some conditions the category of algebras over a colored operad has a structure of a cofibrantly generated model category.

Let $A\subset B$ be an inclusion of $\calX$-algebras, and let $g\colon A\to C$ be a morphism of $\calX$-algebras. We will denote by $\calX_g((B,A),C)$ the space of $\calX$-maps $B\to C$ that coincide with $g$ restricted on $A$.

\begin{lemma}\label{l:extension_space} Let $\calX$ be a discrete $\Sigma$-free colored operad with the set of colors $\calI$ and
 $B$ be an $\calX$ algebra obtained from another $\calX$ algebra  $A$ by an attachment of cells as in~\eqref{eq:21}, then one has the following homeomorphism of spaces:
$$
\calX_g((B,A),C)\simeq \prod_{i\in \calI} \Top_{g_i\circ f_i}((D(i),\partial D(i)), C(i)),
\eqno(\numb)\label{eq:2fiber}
$$
where each factor of the right-hand side is the space of continuous maps $D(i)\to C(i)$ that restrict as $g_i\circ f_i$ on each $\partial D(i)$.
\end{lemma}

\begin{proof}
From the universal property of the pushout, $\calX_g((B,A),C)=\calX_{g\circ\calX(f)}\left((\calX(D),\calX(\partial D)),C\right)$. The latter space is obviously the right-hand side of~\eqref{eq:2fiber} since $\calX(-)$ is the left adjoint functor to the forgetful functor from the category of $\calX$-algebras to the category of $\calI$-collections of topological spaces.
\end{proof}

\begin{lemma}\label{l:fiber_equiv}
In the previous settings assume that for each $i\in\calI$, $D_0(i)\subset D(i)$ is a closed subset in $D(i)$ such that for each connected component of $D(i)$ its intersection with $D_0(i)$ is a closed ball of the same dimension. We assume that the $\calX$-algebras $B$ and $A_+$ are obtained from $A$ using the following pushout diagrams in the category of $\calX$-algebras
$$
\xymatrix{
\calX(\partial D)\ar[r]^{\calX(\iota)}\ar[d]_{\calX(f)}&\calX(D)\ar[d]\\
A\ar[r]&\pushout B,
}
\qquad
\xymatrix{
\calX(\partial D)\ar[r]^-{\calX(\iota)}\ar[d]_{\calX(f)}&\calX(D\setminus Int(D_0))\ar[d]\\
A\ar[r]&\pushout A_+
}
\eqno(\numb)\label{eq:fiber_equiv}
$$
(In the above $D\setminus Int(D_0)$ denotes the family of spaces $D(i)\setminus Int(D_0(i))$, $i\in\calI$, where $Int$ stays for the interior of the corresponding union of balls.)
Thus one has inclusion of $\calX$-algebras $A\subset A_+\subset B$. Let $g\colon A_+\to C$ be a morphism of operads,
then the following natural inclusion of mapping spaces is a homotopy equivalence:
$$
 \calX_g((B,A_+),C)\stackrel{\simeq}{\hookrightarrow}\calX_g((B,A),C).
 \eqno(\numb)\label{eq:restr_map}
$$
\end{lemma}

One can show that in the settings of Lemma~\ref{l:fiber_equiv} the algebra $A_+$ is weakly equivalent to $A$ which means the inclusion $A(i)\hookrightarrow A_+(i)$ is a weak equivalence for every color $i\in\calI$.

\begin{definition}\label{d:punctured}
In the settings of Lemma~\ref{l:fiber_equiv} we will be saying  that $A_+$ is obtained from $A$ by a {\it free attachment of punctured discs}.
\end{definition}

\begin{proof}[Proof of Lemma~\ref{l:fiber_equiv}]
Notice that $B$ is obtained from $A_+$ by an attachment of the same collection of discs, as $B$ is obtained from $A$. Thus the left-hand side of~\eqref{eq:restr_map} is similar to~\eqref{eq:2fiber}. Explicitly the first space is
$$
\prod_{i\in \calI} \Top_{g_i\circ f_i}\bigl(\left(D(i),D(i)\setminus Int(D_0(i))\right), C(i)\bigr),
$$
while the second one is
$$
\prod_{i\in \calI} \Top_{g_i\circ f_i}\bigl((D(i),\partial D(i)), C(i)\bigr),
$$
which are homotopy equivalent since each factor of the first product is naturally homotopy equivalent to the corresponding factor of the second product.
The homotopy is obtained by using the homotopy between $\partial D_0(i)$ and $\partial D(i)$ inside $D(i)\setminus Int(D_0(i))$.
\end{proof}

\begin{definition}\label{d:cof_model}
Let $A$ be an algebra over a discrete colored  operad $\calX$. We say that a cofibrant $\calX$-algebra $\widetilde A$ is a {\it cofibrant model} of $A$ if one has an $\calX$-map
$$
\widetilde{A}\to A
$$
which is  a weak equivalence for each component $A(i)$, $i\in\calI$, where $\calI$ is the set of colors of $\calX$.
\end{definition}

The following is an easy consequence of Definitions~\ref{d:cofibrant_alg},~\ref{d:cof_model}.

\begin{lemma}\label{l:cof_mod}
If $\widetilde A$ and $\widetilde{A}'$ are cofibrant models of an $\calX$-algebra $A$ with the maps $\widetilde{A}\stackrel{f}{\to} A$ and $\widetilde{A}'\stackrel{g}{\to} A$ as in Definition~\ref{d:cof_model}, then there exists an $\calX$-map $\tilde A\stackrel{h}{\to}\tilde A'$ making the following diagram of spaces
$$
\xymatrix{
\widetilde{A}(i)\ar[d]_{h_i}\ar[rd]^{f_i}&\\
\widetilde{A}'(i)\ar[r]_{g_i}&A(i)
}
$$
commute up to homotopy for each $i\in\calI$. Moreover each $h_i$ is a weak equivalence.
\end{lemma}

\begin{lemma}\label{l:cof_maps}
If $\widetilde{A}$ and $\widetilde{A}'$ are cofibrant models of the same $\calX$-algebra, then the spaces of $\calX$-maps $\calX(\widetilde{A},B)$ and ${\calX}(\widetilde{A}',B)$ are homotopy equivalent for any $\calX$-algebra $B$.
\end{lemma}

\begin{proof}
Let $h\colon \widetilde{A}\to\widetilde{A}'$ be an $\calX$-map from Lemma~\ref{l:cof_mod}, and let $h'\colon \widetilde{A}'\to\widetilde{A}$ be its homotopy inverse (which exists by the same lemma). The composition with $h$ and $h'$ induces the maps:
$$
\xymatrix{
{\calX}(\widetilde{A},B)\ar@/^1pc/[r]^{h'\circ{-}}&
{\calX}(\widetilde{A}',B)\ar@/^1pc/[l]^{h\circ{-}},
}
$$
which are homotopy inverses to each other.
\end{proof}

\begin{definition}\label{d:sp_der_morph}
If $A$ and $B$ are algebras over a $\Sigma$ free discrete operad $\calX$ then by the {\it space of derived morphisms between $A$ and $B$} we understand the space
$$
\widetilde{\calX}(A,B):=\calX(\widetilde{A},B)
$$
of $\calX$-maps, where $\widetilde{A}$ is a cofibrant model of $A$. By Lemma~\ref{l:cof_maps} this space up to homotopy is uniquely determined by $A$ and $B$.
\end{definition}

\subsection{Restriction functor}
Let $\calX$ be a  colored operad with $\calI$ as a set of colors. Assuming $\calI_0\subset\calI$ one can define an operad $\calX|_{\calI_0}$ whose set of colors is $\calI_0$ and the operations are exactly those operations of $\calX$ whose all inputs and output are in $\calI_0$. One has an obvious restriction functor
$$
(-)|_{\calI_0}\colon \calX\to\calX|_{\calI_0}
$$
from the category of $\calX$-algebras to the category of $\calX|_{\calI_0}$-algebras.\footnote{We keep abusing notation using the same notation for an operad and the category of algebras over it. Notice that the morphism of operads goes in opposite direction.}

The following obvious results will be frequently used in the paper.

\begin{lemma}\label{l:truncation}
Let $A$ be a cellular cofibrant algebra over a discrete $\Sigma$-free operad $\calX$ such that all its generating cells have color from a subset $\calI_0$, then
$A|_{\calI_0}$ is also a cellular cofibrant $\calX|_{\calI_0}$-algebra with the same generating cells as $A$.
\end{lemma}


\begin{lemma}\label{l:trunc_maps}
Let $A$, $\calX$ and $\calI_0$ be as in the settings of Lemma~\ref{l:truncation}, let also $B$ be any $\calX$-algebra, then one has the following homeomorphism of mapping spaces:
$$
\calX(A,B)\cong \calX|_{\calI_0}(A|_{\calI_0},B|_{\calI_0}).
$$
\end{lemma}

\subsection{Refinement}\label{ss:refin}
Notice that any cofibrant $\calX$-algebra $A$ according to Definition~\ref{d:cofibrant_alg} is naturally endowed with a filtration
$$
\calX(\emptyset)=A_0\subset A_1\subset A_2\subset\ldots\subset A,
$$
where each $A_i$ is obtained from $A_{i-1}$ by an attachment of cells.

\begin{definition}\label{d:refinement}
Given two cofibrant $\calX$-algebras $A$ and $B$, we say that $A$ is a {\it refinement} of $B$ (and $B$ is a {\it coarsening} of $A$) if  for any filtration term $B_i$ in $B$ there exists a filtration term $A_{\alpha(i)}$ in $A$ and a homeomorphism of $\calX$-algebras $f_i\colon B_i\to A_{\alpha(i)}$, such that
for any $i<j$, one has $\alpha(i)<\alpha(j)$, and the following diagram commutes
$$
\xymatrix{
B_i\ar@{^{(}->}[r]\ar[d]_{f_i}&B_{j}\ar[d]^{f_{j}}\\
A_{\alpha(i)}\ar@{^{(}->}[r]&A_{\alpha(j)},
}
$$
and moreover in the limit the inclusions $f_{i}$ induce a homeomorphism $f\colon B\to A$.
\end{definition}

Notice that the last condition simply means that the elements $\alpha(i)$ are not bounded in the indexing ordinal for filtration in~$B$.

%
%

\part{Delooping ignoring degeneracies}\label{part1}

\section{Weak bimodules and (semi)cosimplicial spaces}\label{s:wb}

\subsection{Weak bimodules}\label{ss:wb}

Let $\calO$ be a non-$\Sigma$ topological operad. A weak bimodule $M$ over $\calO$ is a sequence of topological spaces $M=\{M(n),\, n\geq 0\}$ together with the composition maps:
$$
\overline{\circ}_i\colon \calO(n)\times M(k)\to M(n+k-1),\,\, i=1\ldots n, \text{ (weak left action)};
\eqno(\numb)\label{eq_w_left_action}
$$
$$
\underline{\circ}_i\colon M(k)\times \calO(n)\to M(k+n-1),\,\, i=1\ldots k, \text{ (weak right action)}.
\eqno(\numb)\label{eq_w_right_action}
$$
We will be using the same sign $\circ_i$ for $\overline{\circ}_i$, $\underline{\circ}_i$. It will be clear from the context which composition is considered. The left and right weak actions can be depicted as follows:

\vspace{.3cm}

\begin{figure}[h]
\psfrag{OM}[0][0][1][0]{$o\circ_3 m$}
\psfrag{MO}[0][0][1][0]{$m\circ_2 o$}
\psfrag{o}[0][0][1][0]{$o$}
\psfrag{m}[0][0][1][0]{$m$}
\includegraphics[width=16cm]{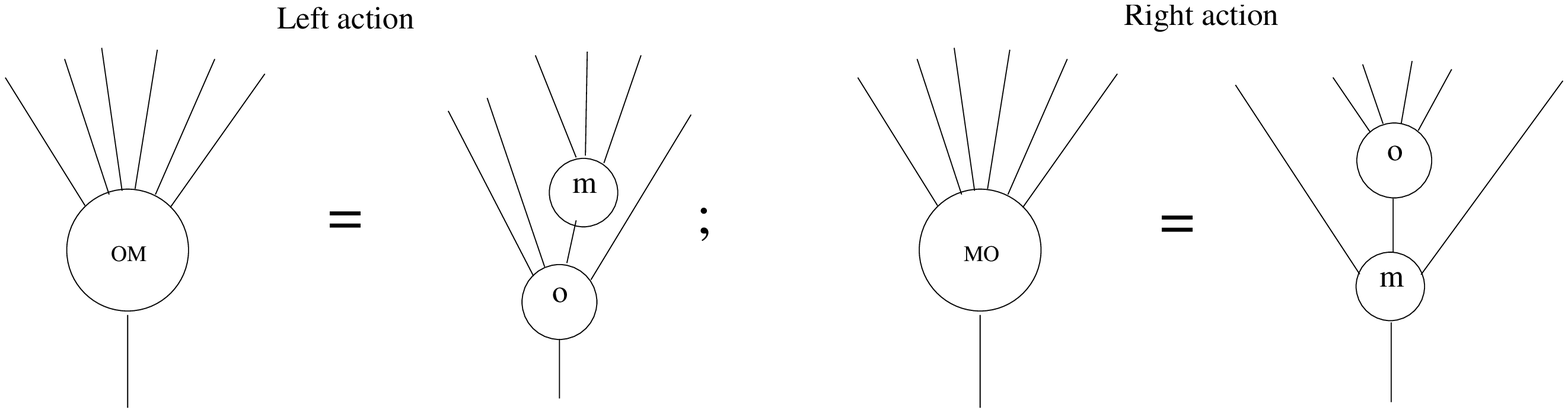}
\caption{}\label{fig1}
\end{figure}

The right and left actions have to satisfy natural unity and associativity conditions. Let $o_1\in O(i)$, $o_2\in O(j)$, and $m\in M(k)$. We will be denoting by $id\in O(1)$ the unit element of the operad $O$. One has the following axioms of the weak bimodule structure:

\vspace{.3cm}

(1) Unity condition with respect to the left and right actions:
$$
id\circ_1 m=m=m\circ_p id,\quad 1\leq p\leq k.
$$

(2) Associativity of the left action:
$$
(o_1\circ_po_2)\circ_{p+q-1}m=o_1\circ_p(o_2\circ_q m), \quad 1\leq p\leq i, \quad 1\leq q\leq j.
$$

(3) Associativity of the right action:
$$
(m\circ_p o_1)\circ_{p+q-1} o_2=m\circ_p(o_1\circ_q o_2), \quad 1\leq p\leq k, \quad 1\leq q\leq i.
$$

(4) Commutativity of the right action on different inputs:
$$
(m\circ_p o_1)\circ_{q+i-1}o_2=(m\circ_q o_2)\circ_p o_1, \quad 1\leq p<q\leq k.
$$

(5) Associativity between the left and right actions:
$$
(o_1\circ_pm)\circ_{p+q-1}o_2=o_1\circ_p(m\circ_q o_2), \quad 1\leq p\leq i, \quad 1\leq q\leq k.
$$

(6) Compatibility between the bi-action and the operad composition:
$$
(o_1\circ_p m)\circ_q o_2=(o_1\circ_q o_2)\circ_{p+q-1} m,\quad 1\leq q<p\leq i.
$$
$$
(o_1\circ_pm)\circ_{q+k-1}o_2=(o_1\circ_q o_2)\circ_p m, \quad 1\leq p<q\leq i.
$$

The weak bimodules were introduced and used by Greg Arone and the second author in~\cite{ArTu,Tu-HDHLK}. In the differential graded context these objects were called {\it infinitesimal bimodules over an operad} by Merkulov and Vallette~\cite{MerkVal}.

\begin{lemma}[\cite{Tu-HDHLK}]\label{l:cosimplicial_semicos}
\textup{(i)} The structure of a weak bimodule over $\Assoc$ is equivalent to the structure of a cosimplicial space.

\textup{(ii)} The structure of a weak bimodule over $\Assoc_{>0}$ is equivalent to the structure of a semi-cosimplical space.
\end{lemma}

For a non-$\Sigma$ operad $\calO$, by $\underset{\calO}{\WBimod}$ we will denote the category of weak bimodules over $\calO$. Given an integer $N\geq 0$, we will also consider the category $\underset{\calO}{\WBimod}{}_N$ of {\it $N$-truncated}  weak bimodules over $\calO$. The objects of  $\underset{\calO}{\WBimod}{}_N$ are finite sequences of spaces $\{M(n),\, 0\leq n\leq N\}$ together with the compositions~\eqref{eq_w_left_action}, and~\eqref{eq_w_right_action} (with the restriction $k\leq N$, $n+k-1\leq N$). One has a natural {\it restriction functor}:
 $$
 (-)|_N\colon \underset{\calO}{\WBimod}\to\underset{\calO}{\WBimod}{}_N.
 $$
 It turns out that a (truncated) $\calO$-weak bimodule is the same thing as a functor from a certain category enriched in topological spaces.  Explicitly this category is described in~\cite{ArTu}. In case $\calO$ is one of the operads $\Assoc$ or $\Assoc_{>0}$, the corresponding category is the simplicial indexing category $\Delta$, or its semisimplicial analogue $s\Delta$, or one of their truncations $\Delta[N]$, $s\Delta[n]$. Since the category of functors from a discrete category is the same thing as the category
of algebras over some discrete colored operad (whose all operations are unary), the results and constructions from Section~\ref{s:col_operads} can be applied to the categories $\underset{\Assoc}{\WBimod}$, $\underset{\Assoc_{>0}}{\WBimod}$, and their truncated analogues  $\underset{\Assoc}{\WBimod}{}_N$, $\underset{\Assoc_{>0}}{\WBimod}{}_N$.

\subsection{A cofibrant model of $\Assoc$ as a weak bimodule over $\Assoc_{>0}$}\label{ss:Triangle}

It is well known that a sequence of simplices $\triangle=\{\triangle(n),\, n\geq 0\}$ defines a cofibrant contractible in each degree semicosimplicial space (cofibrant in the sense that the corresponding functor from the semisimplicial category is cofibrant in the sense of Farjoun~\cite{Farjoun}). Proposition~\ref{p:cof_triangle} below states this fact in our language of cellular cofibrant algebras.

To encode the faces of $\triangle(n)$ we will use planar rooted trees that have one distinguished vertex called {\it bead}:

\vspace{.3cm}

\begin{figure}[h]
\psfrag{D0}[0][0][1][0]{$\triangle(0)$}
\psfrag{D1}[0][0][1][0]{$\triangle(1)$}
\psfrag{D2}[0][0][1][0]{$\triangle(2)$}
\includegraphics[width=6cm]{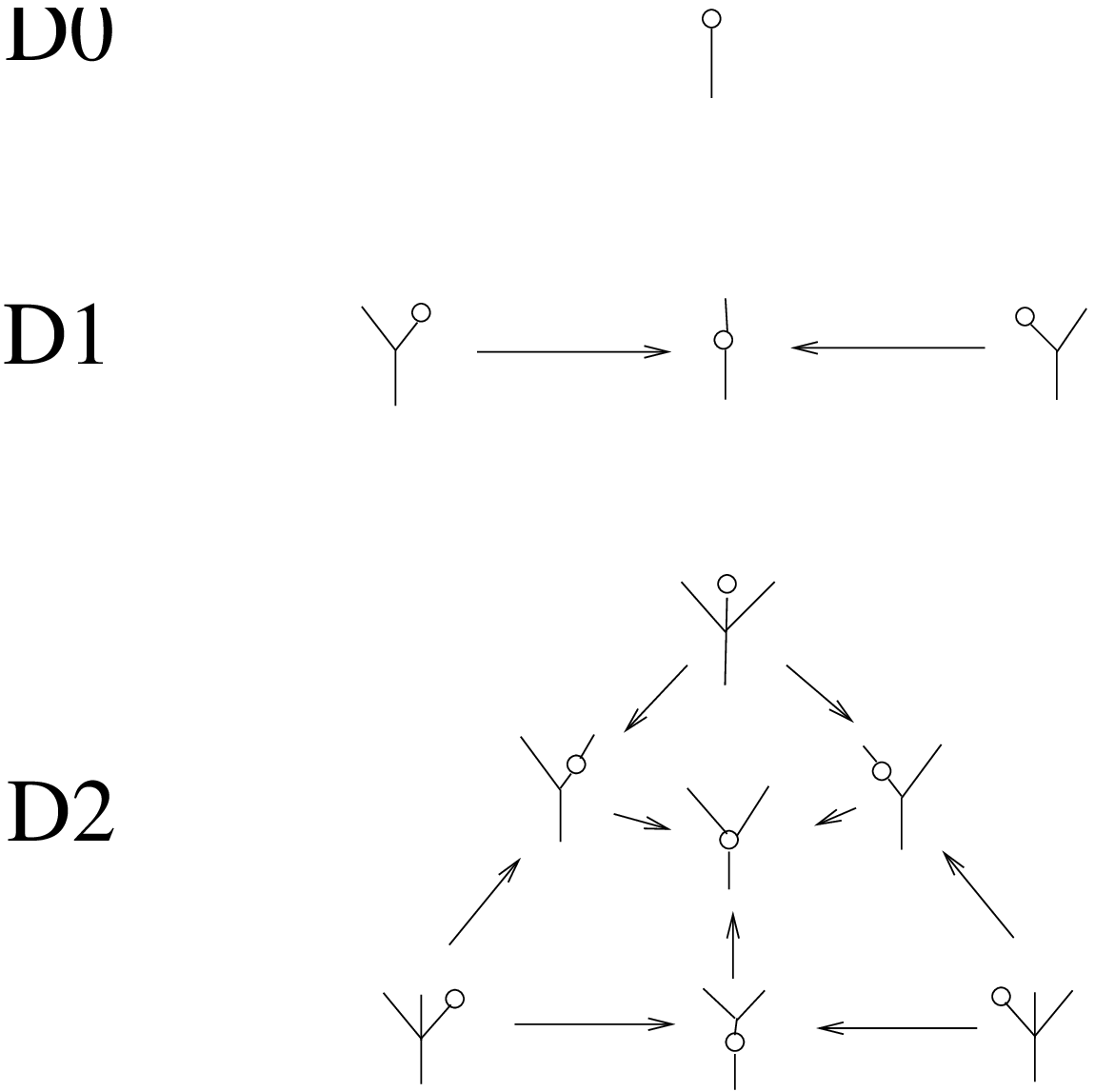}
\caption{Categories of faces of $\triangle(0)$, $\triangle(1)$, $\triangle(2)$.}\label{fig2}
\end{figure}

The trees that are used to encode faces of $\triangle(n)$ have vertices of 4 types: 1 {\it root} of valence 1, $n$ {\it leaves} counted in a clockwise order, one {\it bead} of any positive valence depicted by a little circle (even if its valence is 1, it is not counted for a leaf; similarly the root does not count for a leaf neither), and some number of {\it inner vertices} of any valence~$\geq 3$ (the latter vertices correspond to the action of $\Assoc_{> 0}$. Not all such trees are allowed. One has a restriction that inner vertices can not be joined by an edge (such an edge is automatically contracted by associativity). The dimension of the face corresponding to such tree is the number of edges outgoing from the bead. A typical tree looks like this:
$$
\psfrag{t1}[0][0][1][0]{$t_1$}
\psfrag{t2}[0][0][1][0]{$t_2$}
\psfrag{t3}[0][0][1][0]{$t_3$}
\psfrag{t4}[0][0][1][0]{$t_4$}
\psfrag{t5}[0][0][1][0]{$t_5$}
\psfrag{t6}[0][0][1][0]{$t_6$}
\psfrag{t7}[0][0][1][0]{$t_7$}
\psfrag{t8}[0][0][1][0]{$t_8$}
\psfrag{t9}[0][0][1][0]{$t_9$}
\psfrag{t10}[0][0][1][0]{$t_{10}$}
\psfrag{t11}[0][0][1][0]{$t_{11}$}
\psfrag{t12}[0][0][1][0]{$t_{12}$}
\includegraphics[width=6cm]{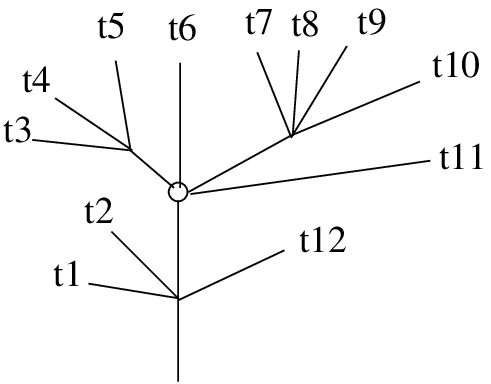}
$$

In coordinates the simplex $\triangle(n)$ is a configuration space of $n$ points on a unit interval:
$$
\triangle(n)=\{(t_1,\ldots,t_n)\, |\, 0\leq t_1\leq t_2\leq \ldots\leq t_n\leq 1\}.
$$
The face encoded by the tree from the above figure is described by the (in)equalities:
$$
0=t_1=t_2<t_3=t_4=t_5<t_6<t_7=t_8=t_9=t_{10}<t_{11}<t_{12} =1.
$$
To obtain these (in)equalities we put the variables $t_1,\ldots,t_{12}$ on the leaves of the tree as in the figure above.
And then let these variables fall one edge down. We get the equality $t_i=0$ (respectively $t_i=1$) if $t_i$ lands below the bead and the $i$-th leaf is to the left (respectively right) of the bead. We get the equality $t_i=t_{i+1}$ if both $t_i$ and $t_{i+1}$ land to the same vertex which is not the bead and in addition the bead is not  between the $i$-th and $(i+1)$-st leaves, otherwise one gets  inequality $t_i<t_{i+1}$.

Notice that as we mentioned earlier, the dimension of a cell equals the number of edges outgoing from the bead.

Denote by $a_k$ the only element in $\Assoc(k)$, $k\geq 0$. We will represent this operation by a tree with a univalent root, and only one inner vertex with $k$ outgoing edges (such tree will be called a {\it $k$-corolla}). For a tree $T$ encoding a face (that we also denote by $T$) of $\triangle(n)$, the tree encoding faces $a_k\circ_i T$ and $T\circ_i a_k$ of $\triangle(n+k-1)$ are obtained by grafting $T$ and the $k$-corolla, and then, if necessary, contracting an edge connecting inner vertices, see example below.
$$
\psfrag{a3}[0][0][1][0]{$a_3$}
\psfrag{c2}[0][0][1][0]{$\circ_2$}
\psfrag{c1}[0][0][1][0]{$\circ_1$}
\includegraphics[width=11cm]{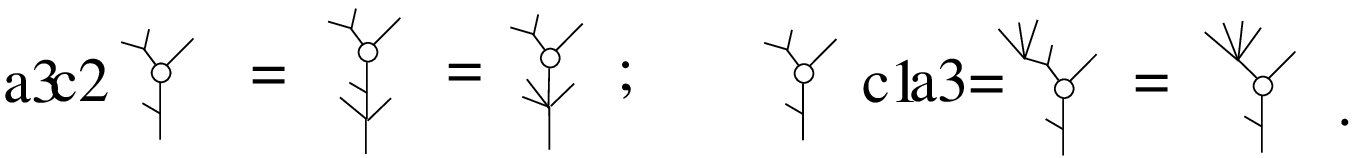}
$$

In coordinates
$$
a_k\circ_i(t_1,\ldots,t_n)=(\underbrace{0\ldots 0}_{i-1},t_1,\ldots,t_n,\underbrace{1\ldots1}_{k-i});
$$
$$
(t_1,\ldots,t_n)\circ_ia_k=(t_1,\ldots,t_{i-1},\underbrace{t_i\ldots t_i}_{k},t_{i+1},\ldots,t_n).
$$

The space $\triangle(n)$ is the configuration space of $n$ points $0\leq t_1\leq\ldots\leq t_n\leq 1$. The left action of $\Assoc_{>0}$ adds endpoints to a configuration, the right action of $\Assoc_{>0}$ multiply points in a configuration.

\begin{proposition}\label{p:cof_triangle}
The weak $\Assoc_{> 0}$ bimodule $\triangle$ and its truncation $\triangle|_N$, $N\geq 0$, are cofibrant models of $\Assoc$, $\Assoc|_N$, respectively, in the categories $\underset{\Assoc_{>0}}{\WBimod}$,  $\underset{\Assoc_{>0}}{\WBimod}{}_N$, respectively.
\end{proposition}

\begin{proof}
Basically one has to show that the sequence of simplices $\triangle(n)$, $n\geq 0$, forms a cofibrant semi-cosimplicial space. This fact is well known. One has that $\triangle$ is obtained from the empty weak $\Assoc_{>0}$ bimodule by a sequence of cell attachments. Let $\triangle_k(n)$ denote the $k$-skeleton of $\triangle(n)$. One has a filtration of $\Assoc_{>0}$ weak bimodules:
$$
\emptyset=\triangle_{-1}\subset\triangle_0\subset\triangle_1\subset\triangle_2\subset\triangle_3\subset\ldots
\eqno(\numb)\label{eq:triangle_filtration}
$$
From the description of cell decomposition in terms of trees it is easy to see that $\triangle_k$ is obtained from $\triangle_{k-1}$ by attaching a $k$-cell to the $k$-th component. Thus all the terms in filtration~\eqref{eq:triangle_filtration} are cofibrant weak $\Assoc_{>0}$ bimodules, and so is  $\triangle$.

The truncated version follows from Lemma~\ref{l:truncation} and from the fact that $\triangle_N|_N=\triangle|_N$.
\end{proof}

\subsection{Tower associated to $\underset{\Assoc_{>0}}{\WBimod}(\triangle,\calO)$}\label{ss:tower_wb}
Let $\calO$ be a weak bimodule over $\Assoc_{>0}$. The space of maps $\underset{\Assoc_{>0}}{\WBimod}(\triangle,\calO)$ is a subspace of the product
$$
\prod_{n\geq 0}\Top(\triangle(n),\calO(n))
$$
of tuples $(f_n)_{n\geq 0}$ of continuous maps $f_n\colon \triangle(n)\to \calO(n)$ such that the restriction of each $f_n$ on the boundary $\partial\triangle(n)$ is determined by the previous map $f_{n-1}$. More precisely for any $x\in\triangle(n-1)$ and $\{a_2\}=\Assoc(2)$, one has
$$
f_n(a_2\circ_ix)=a_2\circ_if_{n-1}(x),\quad i=1,2;
\eqno(\numb)\label{eq:wb_boundary1}
$$
$$
f_n(x\circ_j a_2)=f_{n-1}(x)\circ_ja_2,\quad j=1\ldots n-1.
\eqno(\numb)\label{eq:wb_boundary2}
$$
For any integer $N\geq 0$ we define the space $T_N^\triangle(\calO)$ as a subspace of
$$
\prod_{n=0}^N\Top(\triangle(n),\calO(n))
$$
of tuples $(f_n)_{0\leq n\leq N}$ satisfying~\eqref{eq:wb_boundary1}-\eqref{eq:wb_boundary2}. One obtains a tower of fibrations:
$$
T_0^\triangle(\calO)\leftarrow T_1^\triangle(\calO)\leftarrow T_2^\triangle(\calO)\leftarrow T_3^\triangle(\calO)\leftarrow\ldots
\eqno(\numb)\label{eq:tower_triangle}
$$
Notice that the stages of the tower can also be described as
\begin{multline}
T_N^\triangle(\calO)=\underset{\Assoc_{>0}}{\WBimod}(\triangle_N,\calO)= \underset{\Assoc_{>0}}{\WBimod}{}_N(\triangle|_N,\calO|_N)\\ \simeq\widetilde{\underset{\Assoc_{>0}}{\WBimod}}_N(\Assoc|_N,\calO|_N)\simeq\hosTot_N\calO(\bullet).
\label{eq:stages_triangle}
\end{multline}
The second equation follows from Lemma~\ref{l:truncation}. From the next to the last expression and from Lemma~\ref{l:cof_maps} the stages of the tower~\eqref{eq:tower_triangle} and the tower itself is independent of the cofibrant model.

For $\gamma\in T_{N-1}^\triangle(\calO)$, the fiber of the map $T_N^\triangle(\calO)\to T_{N-1}^\triangle(\calO)$ over $\gamma$ is the space
$\underset{\Assoc_{>0}}{\WBimod}{}_\gamma\left((\triangle_N,\triangle_{N-1}),\calO\right)$. As it follows from Lemma~\ref{l:extension_space} (and from the fact that $\triangle_N$ is obtained from $\triangle_{N-1}$ by a free attachment of one $N$-cell in degree $N$) this fiber is either empty or homotopy equivalent to~$\Omega^N\calO(N)$. Since all the maps in the tower are fibrations, its limit coincides with the homotopy limit and is
$$
T_\infty^\triangle(\calO)=\lim_{\stackrel{\longleftarrow}{N}}T_N^\triangle(\calO)\simeq\holim_{\stackrel{\longleftarrow}{N}}T_N^\triangle(\calO)\simeq \underset{\Assoc_{>0}}{\WBimod}(\triangle,\calO).
\eqno(\numb)\label{eq_limT_Tr}
$$

\section{Bimodules}\label{s:bimodules}

\subsection{(Truncated) bimodules}\label{ss:bimodules}

Let $\calO$ be a non-$\Sigma$ operad. A bimodule $M$ over it is a sequence of spaces $M=\{M(n),\, n\geq 0\}$ together with the structure composition maps:
$$
\gamma_{n;k_1,\ldots,k_n}\colon M(n)\times\calO(k_1)\times\ldots\times\calO(k_n)\to M(k_1+\ldots +k_n) \qquad \text{(right action)};
\eqno(\numb)\label{eq:r_act}
$$
$$
\tilde\gamma_{k;n_1,\ldots,n_k}\colon \calO(k)\times M(n_1)\times \ldots\times M(n_k)\to M(n_1+\ldots +n_k)\qquad \text{(left action)},
\eqno(\numb)\label{eq:l_act}
$$
satisfying natural unity and associativity conditions~\cite{BergMoerd,IwMi,MerkVal}.

Since any operad $\calO$ has the identity element $id\in\calO(1)$, the right action~\eqref{eq:r_act} is equivalent to the weak right action~\eqref{eq_w_right_action}. On the contrary the left action is different from the weak left action. Moreover neither one can be obtained one from another.
Though, in one important case the structure of a weak $\calO$-bimodule is induced from a structure of an $\calO$-bimodule, namely in case when an $\calO$-bimodule $M$ is endowed with a map of $\calO$-bimodules
$
\alpha\colon \calO\to M
$
(where $\calO$ is viewed as a bimodule over itself). Indeed, in this case one has the element $\alpha(id)\in M(1)$ that can be used to mimic empty insertions.

By an $N$-{\it truncated} bimodule we will understand a finite sequence of spaces $\{M(n),\, 0\leq n\leq N\}$ endowed with composition maps~\eqref{eq:r_act},~\eqref{eq:l_act} (in the range where they can be defined) that satisfy the natural unity and associativity properties.
By $\underset{\calO}{\Bimod}$, respectively $\underset{\calO}{\Bimod}{}_N$ we will denote the category of $\calO$-bimodules, respectively $N$-truncated
$\calO$-bimodules. One has an obvious restriction functor
$$
(-)|_N\colon \underset{\calO}{\Bimod}\to \underset{\calO}{\Bimod}{}_N.
$$
Notice that the structure of a bimodule over $\calO$ (as well as that of a truncated bimodule) is governed by some colored operad. In case $\calO$ is $\Assoc$ or $\Assoc_{>0}$ the corresponding colored operads are discrete and $\Sigma$-free. Thus the constructions from Section~\ref{s:col_operads} can be applied to the categories $\underset{\Assoc}{\Bimod}$, $\underset{\Assoc_{>0}}{\Bimod}$, $\underset{\Assoc}{\Bimod}{}_N$, $\underset{\Assoc_{>0}}{\Bimod}{}_N$.

\subsection{Cofibrant model of $\Assoc_{>0}$ as a bimodule over $\Assoc_{>0}$}\label{ss:cof_bim1}
Define $\square=\{\square(n),\, n\geq 0\}$ as
$$
\square(n)=
\begin{cases}
\emptyset,& n=0;\\
[0,1]^{n-1},& n\geq 1.
\end{cases}
\eqno(\numb)\label{eq:square1}
$$

\begin{figure}[h]
\psfrag{S1}[0][0][1][0]{$\square(1)$}
\psfrag{S2}[0][0][1][0]{$\square(2)$}
\psfrag{S3}[0][0][1][0]{$\square(3)$}
\includegraphics[width=6cm]{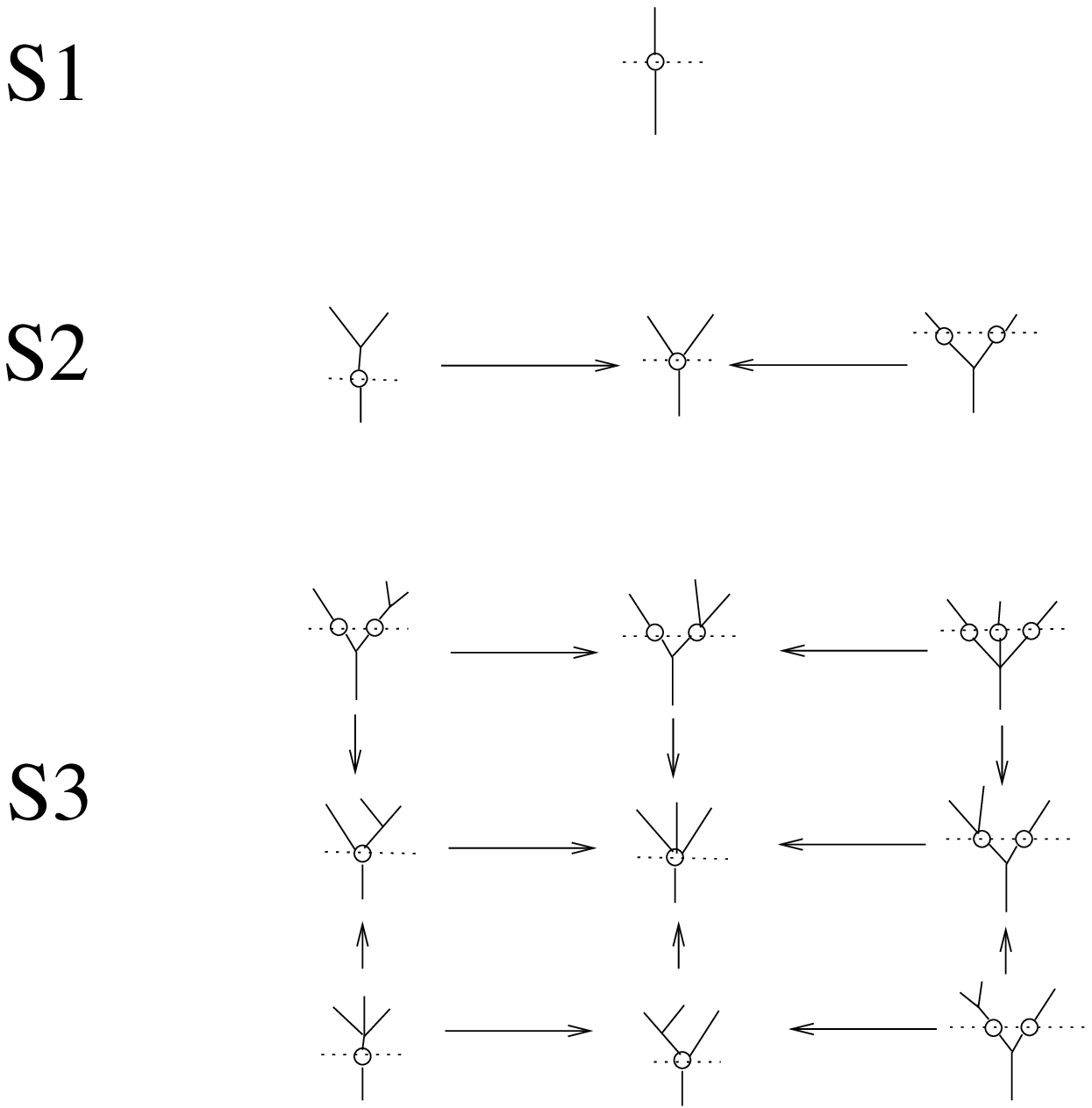}
\caption{Categories of faces of $\square(1)$, $\square(2)$, $\square(3)$.}\label{fig3}
\end{figure}

In other words $\square(n)$, $n\geq 1$, is an $(n-1)$-dimensional cube. The faces of $\square(n)$, $n\geq 0$, will be encoded by planar rooted trees, see Figure~\ref{fig3}.

\vspace{.3cm}

The corresponding trees have vertices of 4 different types: the {\it root} of valence 1, $n$ {\it leaves} counted in a clockwise order, a number of {\it beads} each of valence $\geq 2$ (all lying on the same horizontal line), and some number of {\it inner vertices} all of valence $\geq 3$.

There are two restrictions: an edge can not connect two inner vertices (such edge is automatically contracted by the associativity property); the path between any leaf and the root should pass through one and only one bead. The last condition means that one can draw any such tree $T$ in a way that all the beads lie on the same horizontal line.

Such trees form a poset of faces of a cube $\square(n)$. A face lies in the closure of the other if the tree corresponding to the second one can be obtained from the first tree by a sequence of contractions of edges  and splitting-contraction operations:
$$
\includegraphics[width=12cm]{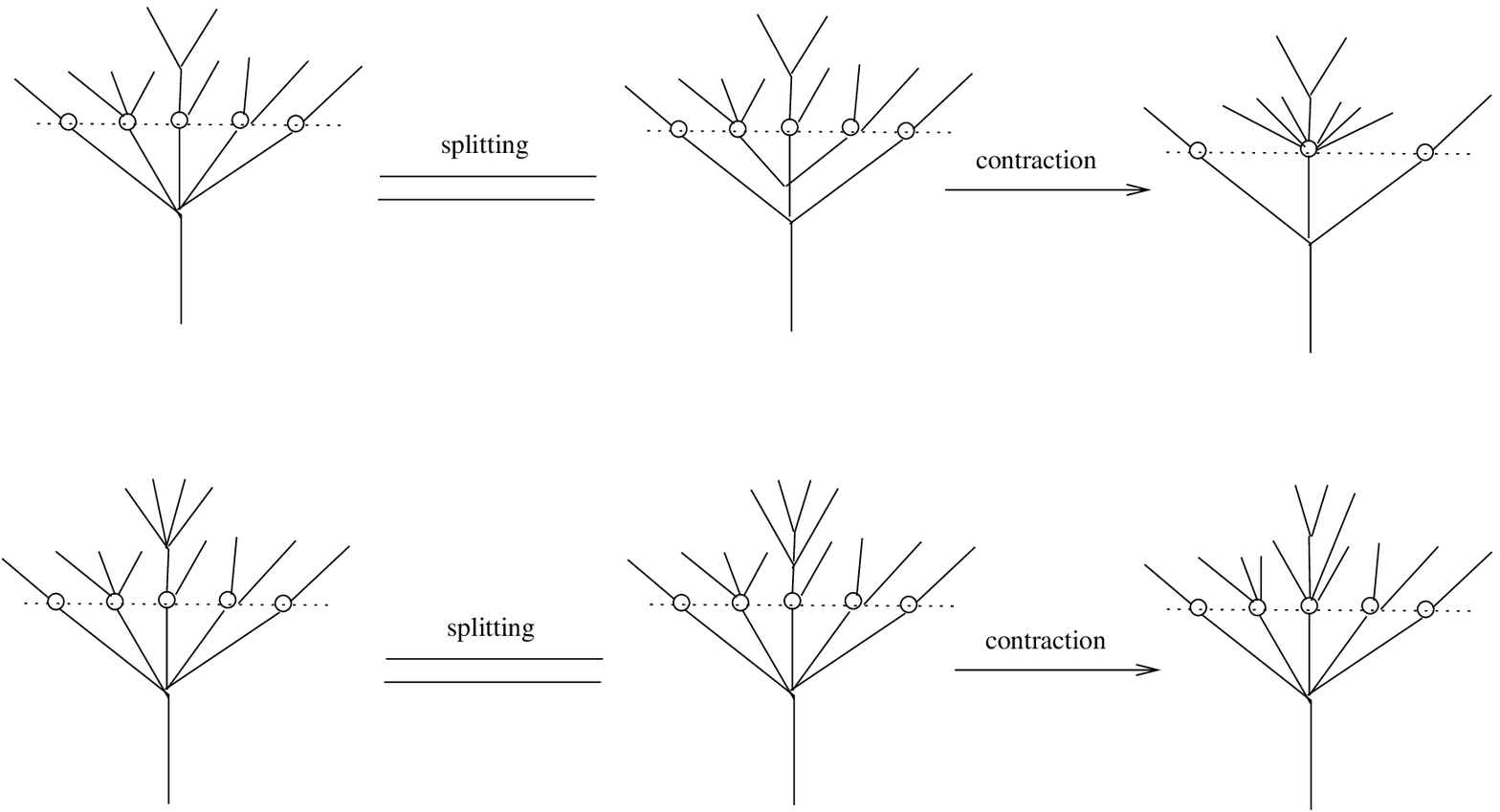}
$$

The interior of $\square(n)$ is encoded by the $n$-corolla \raisebox{-5pt}{\includegraphics[width=.85cm]{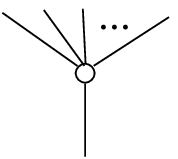}}.
The codimension one faces can have either one or two beads:
\begin{align*}
a_2(\square(\ell),\square(n-\ell))=&\raisebox{-5pt}{\includegraphics[width=1cm]{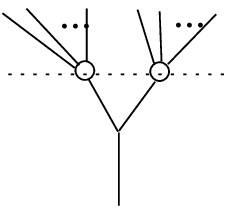}},& \ell=1,\ldots,n-1,\\
\square(\ell)\circ_ia_2=&\raisebox{-7pt}{\includegraphics[width=1cm]{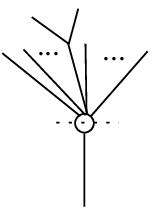}},& i=1\ldots\ell,
\end{align*}

As it was mentioned $\square(n)$ is an $(n-1)$-cube
$$
\square(n)=\{(d_1,\ldots,d_{n-1})\, |\, 0\leq d_i\leq 1,\, i=1\ldots n-1\}
\eqno(\numb)\label{eq:square2}
$$
To see which face of~\eqref{eq:square2} corresponds to a given tree, we put the variables $d_1,\ldots,d_{n-1}$ between the leaves of the tree and let them fall down to the first vertex:
$$
\psfrag{d1}[0][0][1][0]{$d_1$}
\psfrag{d2}[0][0][1][0]{$d_2$}
\psfrag{d3}[0][0][1][0]{$d_3$}
\psfrag{d4}[0][0][1][0]{$d_4$}
\psfrag{d5}[0][0][1][0]{$d_5$}
\includegraphics[width=4cm]{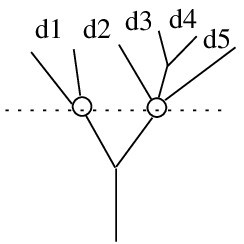}
$$
If $d_i$ falls below the horizontal line passing through the beads, one gets the equation $d_i=1$. If it falls above, one gets $d_i=0$. If it falls on one of the beads, there is no additional restrictions on the variable $d_i$, or in other words one has $0\leq d_i\leq 1$.
For example, in the case of the tree above, one gets the face given by the constraints:
$$
0\leq d_1\leq 1,\quad d_2=1, \quad 0\leq d_3\leq 1, \quad d_4=0,\quad 0\leq d_5\leq 1.
$$

The bimodule action of $\Assoc_{>0}$ on $\square$ is defined by grafting construction, similar to the weak bimodule action on $\triangle$, see Section~\ref{s:wb}. In coordinates
$$
a_2(\square(k),\square(n-k))\subset \square(n),
$$
$$
a_2\left((d_1\ldots d_k),(d_1'\ldots d_{n-k}')\right)=(d_1\ldots d_{k-1},1,d_1'\ldots d_{n-k}');
$$
$$
\square(n-1)\circ_ia_2\subset \square(n),
$$
$$
(d_1\ldots d_{n-2})\circ_ia_2=(d_1\ldots d_{i-1},0,d_i\ldots d_{n-2}).
$$

The space $\square(n)$ can also be viewed as the space of configurations of $n+1$ points on $\R^1$ modulo translations of $\R^1$ with the restriction that the distance between two neighbor points should be non-strictly between~0 and~1. The coordinate $d_i$ is exactly the distance between the $i$-th and $(i+1)$-st points in the configuration.

\begin{proposition}\label{p:cof_bimod}
The $\Assoc_{>0}$ bimodule $\square$ is a cofibrant model of $\Assoc_{>0}$ in the category $\underset{\Assoc_{>0}}\Bimod$. Similarly, $\square|_N$ is a cofibrant model of $\Assoc|_N$ in the category  $\underset{\Assoc_{>0}}{\Bimod}{}_N$.
\end{proposition}

\begin{proof}
The $\Assoc_{>0}$-bimodule $\square$ has a natural filtration
$$
\emptyset=\square_0\subset \square_1\subset \square_2\subset\square_3\subset\ldots\subset \square,
\eqno(\numb)\label{eq:square_filtration}
$$
where in degree $n$, $\square_N(n)$ is a subcomplex of $\square(n)$ consisting only of the faces encoded by the trees whose beads have no more than $N$ outgoing edges. It is clear that $\square_N$ is obtained from $\square_{N-1}$ by a free attachment of an $(N-1)$-cell in degree $N$ (in particular $\square_N(n)=\square(n)$ if $N\geq n$). Thus all the terms of the filtration~\eqref{eq:square_filtration} as well as $\square$ itself are cofibrant. The natural projection $\square\to\Assoc_{>0}$ is obviously a weak equivalence of $\Assoc_{>0}$ bimodules.

The truncated case follows from Lemma~\ref{l:truncation} and from the fact that $\square|_N=\square_N|_N$.
\end{proof}

\subsection{Cofibrant model of $\Assoc$ as an $\Assoc{-}\Assoc_{>0}$-bimodule}\label{ss:cof_bim2}

It might appear strange that in the previous section we define a cofibrant model of $\Assoc_{>0}$ and not of $\Assoc$. Why the degree zero component $\Assoc(0)$ has to disappear? There is a way to keep this component, but in that case we need to consider a slightly different category where a cofibrant model of $\Assoc$ will be defined. We mention that this slightly different approach might seem to have a more natural formulation for the first delooping result Theorem~\ref{t:first_deloop_stages}.

The definition of a bimodule over an operad has a modification in which the operad $\calO_1$ acting from the left is different from the operad $\calO_2$ acting from the right. Such bimodule will be called $\calO_1{-}\calO_2$-bimodule. The category of $\calO_1{-}\calO_2$-bimodules will be denoted by $\underset{\calO_1{-}\calO_2}{\Bimod}$.

We redefine the sequence $\square=\{\square(n),\, n\geq 0\}$ by changing only its zero component to be a point $\square(0):=*$. (Initially this component is defined as the empty set.) Abusing notation the obtained sequence will be also denoted by $\square$. This version of $\square$ will be used in Part~\ref{part2} (but nowhere else in Part~\ref{part1}).


\begin{proposition}\label{p:cof_bimod2}
The sequence of spaces $\square$ has a natural structure of an $\Assoc-\Assoc_{>0}$-bimodule, moreover $\square$ is a cofibrant model of $\Assoc$ in the category $\underset{\Assoc{-}\Assoc_{>0}}{\Bimod}$ of such bimodules. One also has that $\square|_N$ is a cofibrant model of $\Assoc_N$ in the category of $N$-truncated $\Assoc{-}\Assoc_{>0}$-bimodules.
\end{proposition}

The proof is essentially the same as that of Proposition~\ref{p:cof_bimod}. One only has to understand what is the difference between $\Assoc_{>0}$-bimodules and $\Assoc{-}\Assoc_{>0}$-bimodules. The difference is of course that in addition to $\Assoc_{>0}$ left and right actions we also have a left action of $a_0\in\Assoc(0)$ that produces an element $\unit=a_0()$ in the component~0 of an $\Assoc{-}\Assoc_{>0}$-bimodule $M$. This element has to satisfy the unity condition:
$$
a_2(\unit,x)=x=a_2(x,\unit)
\eqno(\numb)\label{eq:unit_bimod}
$$
for any element $x\in M(n)$, $n\geq 0$, which follows from the associativity of the left action:
$$
a_2(\unit,x)=a_2(a_0(),x)=(a_2\circ_1 a_0)\circ_1 x= a_1\circ_1 x=id\circ_1 x=x.
$$

The structure of a (truncated) $\Assoc{-}\Assoc_{>0}$-bimodule is governed by some colored operad. Notice that because of the presence of the operation of arity~0, the free algebra generated by empty sets in each degree always has a point $\unit$ in degree~0 (actually it is its only point). The role of the element $\unit$ in $\square$ is similar to the role of the identity in the Stasheff operad $\Pentagon$, see Subsection~\ref{ss:trunc_operads}.

\subsection{Tower associated to $\underset{\Assoc_{>0}}{\Bimod}(\square,\calO)$}\label{ss:tower_bimod}
Let $\calO$ be a bimodule over $\Assoc_{>0}$. The space of morphisms $\underset{\Assoc_{>0}}{\Bimod}(\square,\calO)$ is a subspace of the product
$$
\underset{\Assoc_{>0}}{\Bimod}(\square,\calO)\subset\prod_{n\geq 0}\Top(\square(n),\calO(n)).
\eqno(\numb)\label{eq:prod_bim_maps}
$$
Since $\square(0)=\emptyset$, the first factor $\Top(\square(0),\calO(0))$ is a point and can be omitted. For the other maps $g_n\colon\square(n)\to\calO(n)$ in a tuple $(g_n)_{n\geq 0}\in \underset{\Assoc_{>0}}{\Bimod}(\square,\calO)$, one has that their restriction on the boundary $g_n|_{\partial\square(n)}$, $n\geq 2$, is determined by the previous $g_1,\ldots,g_{n-1}$:
$$
g_n(a_2(x,y))=a_2(g_\ell(x),g_{n-\ell}(y)),\quad x\in\square(\ell),\, y\in\square(n-\ell),\, 1\leq \ell\leq n-1;
\eqno(\numb)\label{eq:bound_bimod1}
$$
$$
g_n(x\circ_ia_2)=g_{n-1}(x)\circ_ia_2,\quad x\in\square(n-1),\, i=1\ldots n-1.
\eqno(\numb)\label{eq:bound_bimod2}
$$
One can define the space $T_N^\square(\calO)\subset \prod_{n=0}^N \Top(\square(n),\calO(n))$ as the space of  tuples $(g_n)_{0\leq n\leq N}$ satisfying the above conditions. One gets a tower of fibrations
$$
T_0^\square(\calO)\leftarrow T_1^\square(\calO)\leftarrow T_2^\square(\calO)\leftarrow T_3^\square(\calO)\leftarrow \ldots,
\eqno(\numb)\label{eq:tower_Square}
$$
whose both the limit and the homotopy limit (since all maps are fibrations) is $\underset{\Assoc_{>0}}{\Bimod}(\square,\calO)$. The stages of the tower can also be described as
$$
T_N^\square(\calO)=\underset{\Assoc_{>0}}{\Bimod}(\square_N,\calO)=\underset{\Assoc_{>0}}{\Bimod}{}_N(\square|_N,\calO|_N)\simeq \widetilde{\underset{\Assoc_{>0}}{\Bimod}}{}_N(\Assoc_{>0}|_N,\calO|_N),
\eqno(\numb)\label{eq:bimod_stages}
$$
where $\square_N$ is the $N$-th term of the filtration~\eqref{eq:square_filtration}. The second equality follows from Lemma~\ref{l:trunc_maps}. From the last expression we see that the stages of the tower and the tower itself are independent of the cofibrant model.

The maps in the tower are fibrations and as it follows from Lemma~\ref{l:extension_space}
(and from the fact that $\square_N$ is obtained from $\square_{N-1}$ by a free attachment of one $(N-1)$-cell in degree $N$) the fiber over any point under the map
$$
T_N^\square(\calO)\to T_{N-1}^\square(\calO)
$$
is either empty or homotopy equivalent to $\Omega^{N-1}\calO(N)$.


\section{Operads}\label{s:operads}

\subsection{(Truncated) operads}\label{ss:trunc_operads}
We assume that the reader is familiar with the notion of an operad. The category of topological non-$\Sigma$ operads will be denoted $\Operad$. By an $N$-truncated operad we will understand a sequence of spaces $\{\calO(n),\, n=0\ldots N\}$ together with the operadic compositions (in the range where they can be defined) and an identity element $id\in\calO(1)$ (if $N\geq 1$), satisfying the usual associativity and unity conditions. The category of $N$-truncated non-$\Sigma$ operads in topological spaces will be denoted by $\Operad_N$. One has a restriction functor
$$
(-)|_N\colon\Operad\longrightarrow\Operad_N,
$$
that forgets all the components of an operad except the first $(N+1)$ ones. Both the structure of an operad and of its truncated analogue are governed by some $\Sigma$-free discrete colored operads. We conclude that the constructions of Section~\ref{s:col_operads} can be applied to the categories $\Operad$, $\Operad_N$.

\subsection{Cofibrant model for the operad $\Assoc_{>0}$}\label{ss:cof_operad}
\begin{figure}[h]
\psfrag{p2}[0][0][1][0]{$\Pentagon(2)$}
\psfrag{p3}[0][0][1][0]{$\Pentagon(3)$}
\psfrag{p4}[0][0][1][0]{$\Pentagon(4)$}
\includegraphics[width=7cm]{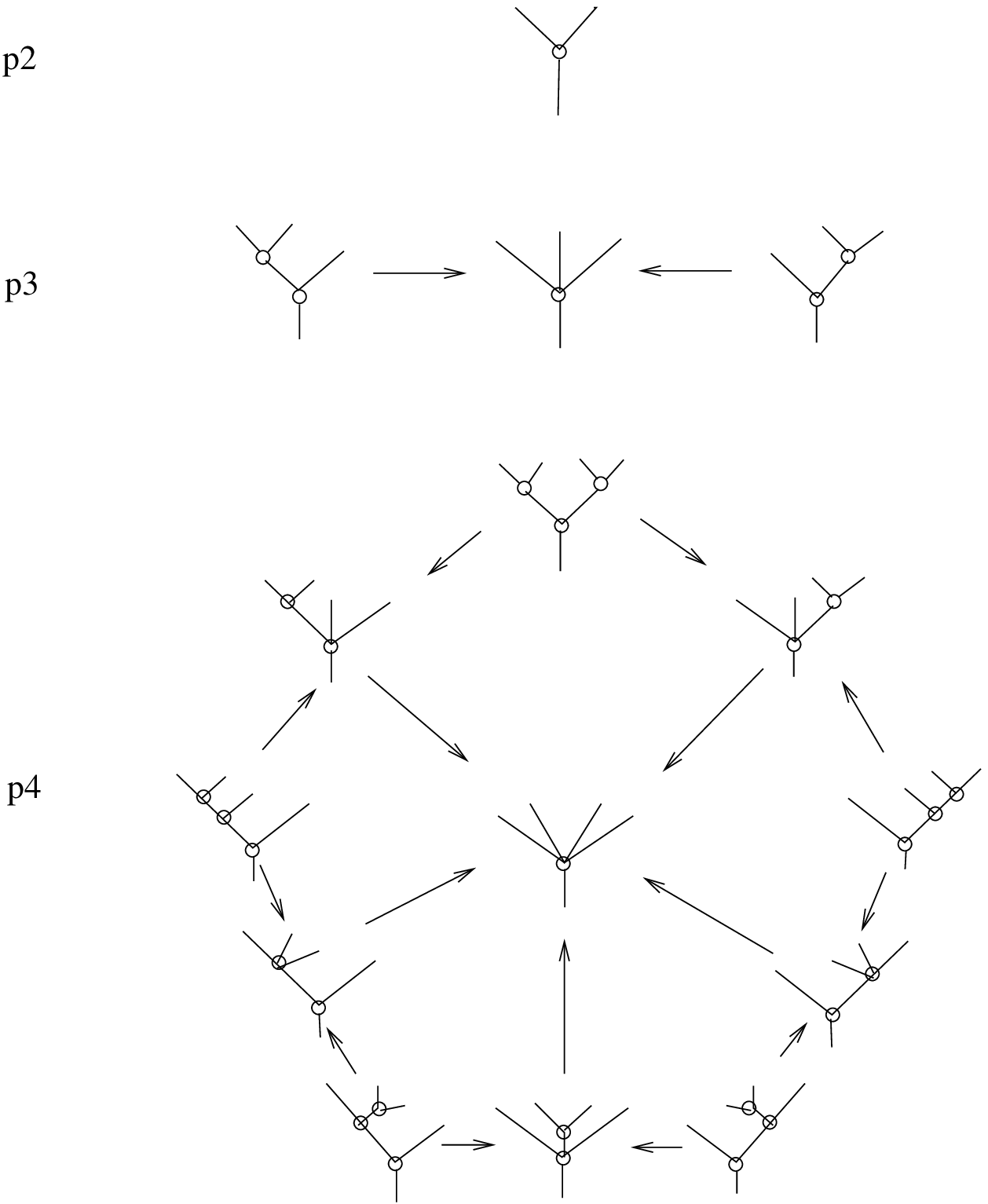}
\caption{Categories of faces of $\Pentagon(2)$, $\Pentagon(3)$, $\Pentagon(4)$.}\label{fig4}
\end{figure}

It is well known that the Stasheff operad, that we will denote by $\Pentagon$, is a cofibrant model for the operad $\Assoc_{>0}$~\cite{Stasheff-HAH,Stasheff-HHPV}. One has
$$
\Pentagon(n)=\begin{cases}
\emptyset,& n=0,\\
*,& n=1,\\
K_{n-2},& n\geq 2,
\end{cases}
\eqno(\numb)\label{eq:pentagon_components}
$$
where $K_{n-2}$ is the $(n-2)$-dimensional Stasheff polytope~\cite{Stasheff-HAH,Stasheff-HHPV}. The faces of $\Pentagon(n)$, $n\geq 2$, are encoded by planar rooted trees with $n$ {\it leaves} (labeled by $1\ldots n$ in a clockwise order), a root of valence 1, and a bunch of other vertices of valence $\geq 3$. For a convenience of exposition these other vertices will be referred as {\it beads} and will be depicted by little circles. The category of faces of $\Pentagon(n)$ will be denoted by $\Psi_n$. Figure~\ref{fig4}  shows the categories $\Psi_2$, $\Psi_3$, and $\Psi_4$.

\vspace{.3cm}

The dimension of the face corresponding to a tree $T$ is $\sum_{b\in B(T)}(|b|-2)$, where $B(T)$ is the set of beads of $T$, and $|b|$ denotes the number of edges outgoing from $b$ (which is valence minus one).

\begin{proposition}\label{p:operad_cof}
The Stasheff operad $\Pentagon$ and its restriction $\Pentagon|_N$, $N\geq 0$, provide cofibrant models for $\Assoc_{>0}$ and $\Assoc_{>0}|_N$ in the categories $\Operad$, $\Operad_N$ respectively.
\end{proposition}

\begin{proof}
Consider a natural filtration
$$
\Pentagon_0\subset\Pentagon_1\subset\Pentagon_2\subset\Pentagon_3\subset\ldots\subset\Pentagon,
\eqno(\numb)\label{eq:pentag_filtr}
$$
where $\Pentagon_0=\Pentagon_1$ are the operads having only identity element (they are free objects in $\Operad$ generated by empty sets in each degree). The  component $\Pentagon_N(n)$ consists of those faces of $\Pentagon(n)$ whose corresponding trees have only beads of valence $\leq N+1$. Assuming $N\geq 2$, the operad $\Pentagon_N$ is obtained from $\Pentagon_{N-1}$ by attaching an $(N-2)$-cell to the $N$-th component. Thus all the operads in the filtration~\eqref{eq:pentag_filtr} and $\Pentagon$ itself are cellular cofibrant. On the other hand one also has that the projection $\Pentagon\to\Assoc_{>0}$ is a homotopy equivalence in each degree.

The truncated case follows from Lemma~\ref{l:truncation} and from the fact that $\Pentagon|_N=\Pentagon_N|_N$.
\end{proof}

\subsection{Tower associated to {\sc$\Operad(\Pentagon,\calO)$}}\label{s:tower_operad}
The space of morphisms of operads $\Operad(\Pentagon,\calO)$ is a subspace in the product
$$
\Operad(\Pentagon,\calO)\subset\prod_{n\geq 0}\Top(\Pentagon(n),\calO(n)).
\eqno(\numb)\label{eq:oper_subprod}
$$
Notice that the first 2 factors in the above product can be ignored since $\Pentagon(0)=\emptyset$, and $\Pentagon(1)$ consists of only one element $id$, which should be sent to the identity element in $\calO(1)$. Thus a morphism of operads $h\colon\Pentagon\to\calO$ is given by a tuple $(h_n)_{n\geq 2}$ with $h_n\colon\Pentagon(n)\to\calO(n)$, whose restriction $h_n|_{\partial\pentagon}$ is determined by the previous maps $h_2,\ldots,\, h_{n-1}$:
$$
h_n(x\circ_i y)=h_k(x)\circ_i h_{n-k+1}(y),\quad x\in\Pentagon(k),\, y\in\Pentagon(n-k+1),\, i=1\ldots k,\, k=2\ldots n-1.
\eqno(\numb)\label{eq:pentag_bound}
$$
Define the tower
$$
T_0^{\pentagon}(\calO)\gets T_1^{\pentagon}(\calO)\gets T_2^{\pentagon}(\calO)\gets T_3^{\pentagon}(\calO)\gets  \ldots,
\eqno(\numb)\label{eq:pent_tower}
$$
where $T_0^{\pentagon}(\calO)=T_1^{\pentagon}(\calO)=*$, and $T_N^{\pentagon}(\calO)$, $N\geq 2$, is the space of tuples $(h_n)_{2\leq n\leq N}$ satisfying the above properties~\eqref{eq:pentag_bound}. The space $T_N^{\pentagon}(\calO)$ can also be viewed as the space of morphisms of (truncated) operads:
$$
T_N^{\pentagon}(\calO)=\Operad(\Pentagon_N,\calO)=\Operad_N(\Pentagon|_N,\calO|_N)\simeq \widetilde{\Operad}_N(\Assoc_{>0}|_N,\calO|_N),
$$
where $\Pentagon_N$ is the $N$-th term of the filtration~\eqref{eq:pentag_filtr}. The second equation follows from Lemma~\ref{l:trunc_maps}. The last expression describes the stages of the tower as spaces of derived morphisms which shows that the stages of the tower and the tower itself are independent of the cofibrant model of the operad $\Assoc_{>0}$ (see Lemma~\ref{l:cof_maps}).

All the maps
$
T_N^{\pentagon}(\calO)\to T_{N-1}^{\pentagon}(\calO)
$
of the tower~\eqref{eq:pent_tower} are fibrations and according to Lemma~\ref{l:extension_space} (and to the fact that for $N\geq 2$, the operad $\Pentagon_N$ is obtained from $\Pentagon_{N-1}$ by a free attachment of one $(N-2)$-dimensional cell in degree $N$) the preimage of any point under such map is either empty or homotopy equivalent to $\Omega^{N-2}\calO(N)$.


\section{First delooping}\label{s:first_deloop}
Let $\eta\colon\Assoc\longrightarrow\calO$
be a morphism of $\Assoc_{>0}$-bimodules. We will be assuming that $\calO(0)\simeq *$. We stress the fact that we need $\eta$ to be a map from $\Assoc$ and not from $\Assoc_{>0}$. Technically it means that one has to fix an additional element $\eta(a_0)\in\calO(0)$ satisfying
$$
a_2(\eta(a_0),\eta(a_1))=\eta(a_1)=a_2(\eta(a_1),\eta(a_0)).
\eqno(\numb)\label{eq:wunit_condition}
$$
If we denote $\eta(a_0)$ by $\unit$, the property~\eqref{eq:wunit_condition} is a partial case of~\eqref{eq:unit_bimod} automatically satisfied in case $\calO$ is an $\Assoc{-}\Assoc_{>0}$-bimodule.

One has a map of $\Assoc_{>0}$-bimodules $p\colon\square\to\Assoc$,
which is uniquely defined since $\Assoc$ is the terminal object in $\underset{\Assoc_{>0}}{\Bimod}$. By $\Omega\,\underset{\Assoc_{>0}}{\Bimod}(\square,\calO)$ we will understand the loop space with the base point $\eta\circ p$. In the same way we define $\Omega\, T_N^\square(\calO)$ as a loop space with the base point $(\eta\circ p)|_N$. One gets a tower of fibrations
$$
\Omega\, T_0^\square(\calO)\gets \Omega\, T_1^\square(\calO)\gets\Omega\, T_2^\square(\calO)\gets\Omega\, T_3^\square(\calO)\gets\ldots
\eqno(\numb)\label{eq:tower_loop_square}
$$
On the other hand, due to the morphism $\eta\colon\Assoc\to\calO$,
the sequence of spaces $\calO$ inherits a structure of a weak $\Assoc_{>0}$-bimodule, and the tower~\eqref{eq:tower_triangle}
can be defined.

In this section we will construct a map of towers
$$
\xymatrix{
\Omega\, T_0^\square(\calO)\ar[d]_{\xi_0}&\Omega\, T_1^\square(\calO)\ar[d]_{\xi_1}\ar[l]&\Omega\, T_2^\square(\calO)\ar[d]_{\xi_2}\ar[l]&\Omega\, T_3^\square(\calO)\ar[d]_{\xi_3}\ar[l]&\ldots\ar[l]\\
T_0^\triangle(\calO)&T_1^\triangle(\calO)\ar[l]&T_2^\triangle(\calO)\ar[l]&T_3^\triangle(\calO)\ar[l]&\ldots\ar[l],
}
\eqno(\numb)\label{eq:map1_of_towers}
$$
that induces a map of their limits
$$
\Omega\,\underset{\Assoc_{>0}}{\Bimod}(\square,\calO)\stackrel{\xi_\infty}{\longrightarrow}\underset{\Assoc_{>0}}{\WBimod}(\triangle,\calO).
\eqno(\numb)\label{eq:first_deloop_xi}
$$

\begin{theorem}\label{t:first_deloop_stages}
Let $\calO$ be an $\Assoc_{>0}$-bimodule with $\calO(0)\simeq*$ and endowed with a map of $\Assoc_{>0}$-bimodules $\eta\colon\Assoc\to\calO$. Then the map
$\xi_N\colon\Omega\, T_N^\square(\calO)\to T_N^\triangle(\calO)$ constructed below is a homotopy equivalence for each $N\geq 0$.
\end{theorem}

An immediate corollary is the following.

\begin{theorem}\label{t:first_deloop}
In the settings of Theorem~\ref{t:first_deloop_stages}, the induced map of limits~\eqref{eq:first_deloop_xi} is a homotopy equivalence.
\end{theorem}

\begin{proof}[Proof of Theorem~\ref{t:first_deloop}]
In both towers~\eqref{eq:tower_triangle} and~\eqref{eq:tower_loop_square} all maps are fibrations. Therefore for each of them the limit coincides with the homotopy limits. On the other hand by Theorem~\ref{t:first_deloop_stages} the morphism of towers $\xi_\bullet\colon \Omega T_\bullet^\square(\calO)\to T_\bullet^\triangle(\calO)$ is a homotopy equivalence for each stage, therefore the induced map of homotopy limits is a weak equivalence.
\end{proof}

Before constructing the maps $\xi_\bullet$, let us compare the stages $\Omega\, T_1^\square(\calO)$ and $T_1^\triangle(\calO)$. Notice by the way that $T_0^\triangle(\calO)=\Top(\triangle(0),\calO(0))=\calO(0)\simeq *$, and so are $T_0^\square(\calO)=\Omega\, T_0^\square(\calO)=*$. One has $T_1^\square(\calO)=
\Top(\square(1),\calO(1))=\Top(*,\calO(1))=\calO(1)$. So, $\Omega\, T_1^\square(\calO)=\Omega\,\calO(1)$. The space $T_1^\triangle(\calO)=\underset{\Assoc_{>0}}{\WBimod_1}(\triangle|_1,\calO|_1)$ is homotopy equivalent to the fiber of the map $T_1^\triangle(\calO)\to T_0^\triangle(\calO)$, since the target space is contractible and the above map is a fibration. The fiber of this map over $\eta(a_0)\in\calO(0)$  consists of maps $f_1\colon \triangle(1)\to\calO(1)$ having a fixed behavior at the boundary $\partial\triangle(1)$.
But by~\eqref{eq:wunit_condition} $f_1$ at both ends of the interval $\triangle(1)$ is equal to $\eta(a_1)$ which is viewed as a base-point of $\calO(1)$. As a consequence we see that this fiber is homeomorphic to $\Omega\,\calO(1)$. It turns out that $\xi_1$ sends $\Omega\, T_1^\square(\calO)$ homeomorphically to the fiber above. As a consequence $\xi_1$ is a homotopy equivalence.



To define $\xi_N$ we will decompose each $\triangle(n)$, $1\le n\leq N$, into a union of polytopes, such that each of the polytopes will have a natural map to $\calO(n)$ that arises from $g\in\Omega\,\underset{\Assoc_{>0}}{\Bimod_n}(\square|_n,\calO|_n)$. Explicitly such $g$ is a tuple $(g_i)_{1\leq i\leq n}$, with $g_i\colon\square(i)\times [0,1]\to\calO(i)$, satisfying
$$
g_i(x,0)=g(x,1)=\eta(a_i),
$$
for all $x\in\square(i)$. Sometimes we will view each $g_i$ as a map from the  suspension
$$
g_n\colon\Sigma\square(n)\to\calO(n).
$$
For each $t\in[0,1]$, the tuple $(g_n(-,t))_{1\leq n\leq N}$, where $g_n(-,t)\colon\square(n)\to\calO(n)$, satisfies the boundary conditions~\eqref{eq:bound_bimod1}-\eqref{eq:bound_bimod2}.

Let $\Xi_n$ be the subcategory of faces of $\square(n)$ corresponding only to those trees that do not have inner vertices above beads.
One can easily see that $\Xi_n$ is an $(n-1)$-cube or in other words the $(n-1)$-th power of the category $\{0\to 1\}$. Define two functors: a covariant
$$
\lambda_\square\colon\Xi_n\to\Top,
$$
and a contravariant
$$
\chi_\blacktriangle\colon\Xi_n\to\Top.
$$
The functor $\lambda_\square$ assigns to a tree $T\in Obj(\Xi_n)$ the (closed) face of $\square(n)$ corresponding to $T$, and to a morphism between trees the corresponding inclusion of faces. Explicitly
$$
\lambda_\square(T)=\prod_{b\in B(T)}\square(|b|),
\eqno(\numb)\label{eq:lambda_square}
$$
where $B(T)$ is the set of beads of $T$, and $|b|$ is the number of edges outgoing from $b$. This functor describes the left action of $\Assoc_{>0}$ on $\square$.
For example, for $n=3$ the corresponding $\Xi_3$ diagram is  the 2-cube:
$$
\raisebox{.9cm}[1.4cm]{\xymatrix{
\square(1)\times\square(2)\ar[d]_{a_2}&\square(1)\times\square(1)\times\square(1)\ar[l]_-{id\times a_2}\ar[d]^{a_2\times id}\\
\square(3)&\square(2)\times\square(1)\ar[l]_-{a_2}.
}}
$$
If the number of beads $|B(T)|$ in a tree $T$ is $k$, then the dimension of the cube $\lambda_\square(T)$ is $n-k$. If $b_1,\ldots,b_k$ are the beads of $T$ appearing in a clockwise order, a point in $\lambda_\square(T)$ will be denoted by $(x_1,\ldots,x_k)$, where each $x_i\in\square(|b_i|)$.

The functor $\chi_\blacktriangle$ assigns the simplex $\Delta^k$ to a tree $T$ with $k$ beads:
$$
\chi_\blacktriangle(T)=\{(t_1,t_2,\ldots,t_k)\, |\, 0\leq t_1\leq t_2\leq \ldots \leq t_k\leq 1\}.
$$
In other words points in $\chi_\blacktriangle(T)$ are non-decreasing functions
from the set of beads $B(T)$ to $[0,1]$. Given a morphism of trees in $\Xi_n$:
$$
\alpha\colon T_1\to T_2,
$$
the corresponding map of simplices
$$
\chi_\blacktriangle(\alpha)\colon\chi_\blacktriangle(T_2)\to\chi_\blacktriangle(T_1)
\eqno(\numb)\label{eq:chi_face_incl}
$$
is a face inclusion. More precisely $\alpha$ induces a surjective map of the sets of beads
$$
\alpha_*\colon B(T_1)\to B(T_2).
$$
The map~\eqref{eq:chi_face_incl} is defined as a composition with $\alpha_*$:
$$
\chi_\blacktriangle(\alpha)(f)=f\circ\alpha_*.
$$
For example for the morphism
$$
\includegraphics[width=5cm]{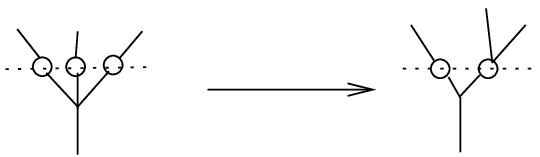}
$$
the inclusion is $\Delta^2\hookrightarrow\Delta^3$ given in coordinates by $(t_1,t_2)\mapsto (t_1,t_2,t_2)$.

One has a natural evaluation map
$$
\bar\xi_T^g\colon\lambda_\square(T)\times\chi_\blacktriangle(T)\to\calO(n),
\eqno(\numb)\label{eq:ev_map1}
$$
given by
$$
\bar\xi_T^g\bigl((x_1,\ldots,x_k),(t_1,\ldots,t_k)\bigr)=a_k\bigl(g_{|b_1|}(x_1,t_1),\ldots,g_{|b_k|}(x_k,t_k)\bigr).
$$

Define $\overline{Wb\square}(T)$ as $\lambda(T)\times\chi_\blacktriangle(T)$, which is a product of a $k$-simplex with an $(n-k)$-cube, thus it is an $n$-dimensional prism. Define also
$$
\overline{Wb\square}(n):=\lambda_\square\otimes_{\Xi_n}\chi_\blacktriangle.
$$
This space has a natural cell decomposition coming from the cell decomposition of its components  $\overline{Wb\square}(T)$.

\begin{proposition}\label{p:cell_dec_W_Square}
The space $\overline{Wb\square}(n)=\lambda_\square\otimes_{\Xi_n}\chi_\blacktriangle$ as a cell complex  is homeomorphic to an $n$-simplex divided into $2^{n-1}$ prisms (each prism corresponding to some  $\overline{Wb\square}(T)$) by $(n-1)$ hyperplanes obtained as follows: One fixes an edge $e$ of the $n$-simplex and a parallel to it line $\ell$ having non-trivial intersection with the interior of the simplex. The $(n-1)$ planes are chosen to contain $\ell$ and to be parallel to one of the $(n-1)$ facets of the simplex containing~$e$.
\end{proposition}

\vspace{.3cm}

\begin{figure}[h]
\psfrag{S1}[0][0][1][0]{$\overline{Wb\square}(1)$}
\psfrag{S2}[0][0][1][0]{$\overline{Wb\square}(2)$}
\psfrag{S3}[0][0][1][0]{$\overline{Wb\square}(3)$}
\includegraphics[width=14cm]{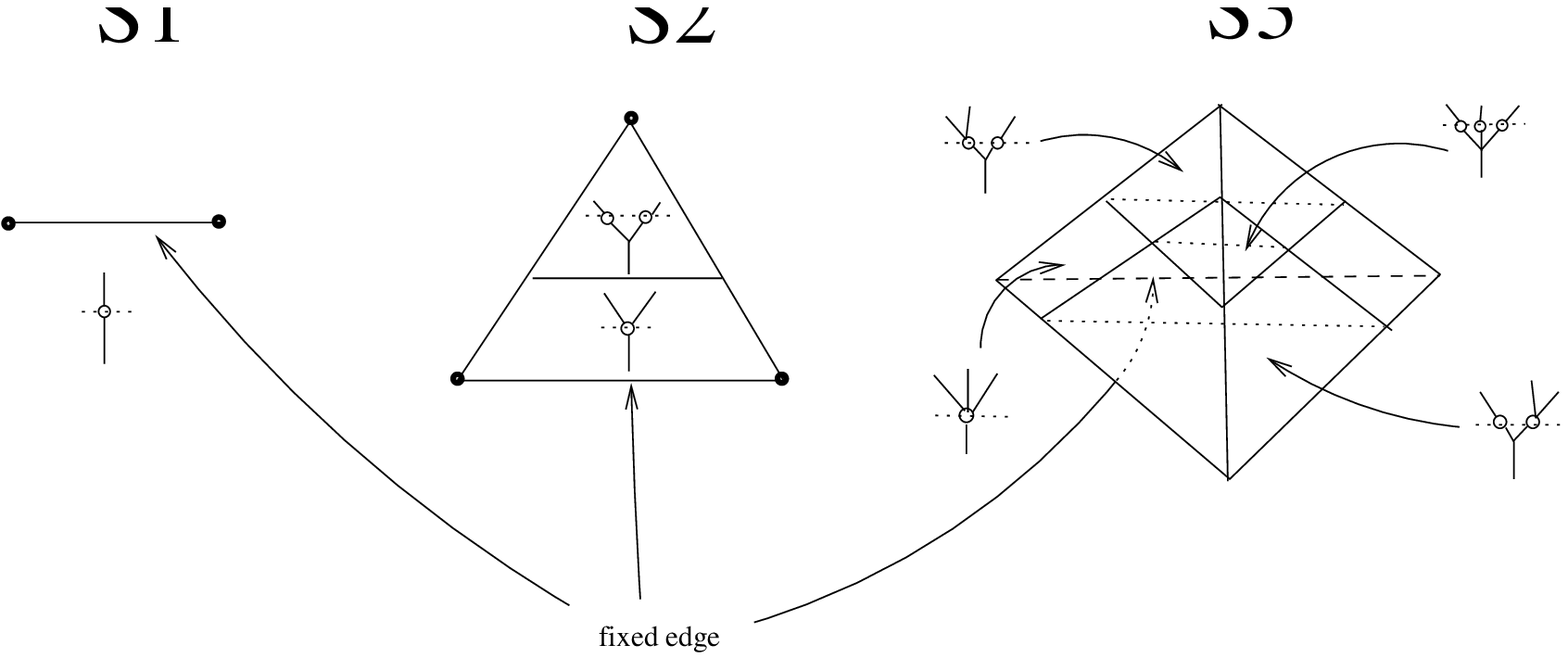}
\caption{Subdivision of $\overline{Wb\square}(1)$, $\overline{Wb\square}(2)$, $\overline{Wb\square}(3)$ into prisms. Each tree $T$ in the above picture denote the corresponding prism $\overline{Wb\square}(T)$.}\label{fig5}
\end{figure}

\begin{proof}
The proof goes by induction over $n$. For small $n=1$, 2, 3, see figures above. Notice that the prism containing the fixed edge $e$ is $Wb\square(C_n)=\square(n)\times [0,1]=[0,1]^n$ where $C_n$ is the $n$-corolla - the terminal element in $\Xi_n$.
\end{proof}

\begin{remark}\label{r:ev_w_sq}
The evaluation maps~\eqref{eq:ev_map1} coincide on the faces of $\overline{Wb\square}(T)$'s that are glued together in $\overline{Wb\square}(n)$ and all together define a map
$$
\bar\xi_n^g\colon\overline{Wb\square}(n)\to\calO(n).
\eqno(\numb)\label{eq:ev_map2}
$$
This is so because for any $t\in[0,1]$ the map $g(t,-)\colon\square\to\calO$ is a morphism of left $\Assoc_{>0}$ modules.
\end{remark}

The sequence of spaces $\overline{Wb\square}=\{\overline{Wb\square}(n),\, n\geq 1\}$ is a sequence of simplices, but this sequence still does not have a structure of a weak $\Assoc_{>0}$ bimodule. The problem is that such structure could not be compatible with the cell decomposition of the spaces in the sequence. Indeed, if  $\overline{Wb\square}$ had such structure in a way that it would be homeomorphic to $\triangle$ as a weak $\Assoc_{>0}$ bimodule, then all the codimension 1 faces of each  simplex $\overline{Wb\square}(n)$ would have the same cell decomposition as that of  $\overline{Wb\square}(n-1)$. But we can see that it is not true even for  $\overline{Wb\square}(2)$. To remedy this we will define  a  quotient of  $\overline{Wb\square}(n)$ that will be denoted by $Wb\square(n)$. First for $T\in Obj(\Xi_n)$, define $Wb\square(T)= \overline{Wb\square}(T)/\sim$, where
$$
\bigl((x_1,\ldots,x_k),(t_1,\ldots,t_k)\bigr)\sim\bigl((x'_1,\ldots,x'_k),(t'_1,\ldots,t'_k)\bigr),
\eqno(\numb)\label{eq:W_relation}
$$
if for all $i=1\ldots k$ one has either $(x_i,t_i)=(x'_i,t'_i)$, or $t_i=t'_i=\varepsilon$ with $\varepsilon\in\{0,1\}$. Notice that one still can define evaluation map
$\xi^g_T\colon Wb\square(T)\to\calO(n)$, such that the following diagram commutes
$$
\xymatrix{
\overline{Wb\square}(T)\ar[rr]^{\bar\xi_T^g}\ar[rd]&&\calO(n)\\
&Wb\square(T)\ar[ru]_{\xi_T^g}
}.
\eqno(\numb)\label{eq:eval_map3}
$$
This is so, because $g_i(x,\varepsilon)=\eta(a_i)$ for any $x\in\square(i)$, $i\geq 1$, $\varepsilon\in\{0,1\}$.

The above quotient map collapses only those codimension 1 faces of  $\overline{Wb\square}(T)$ that appear  in the boundary of the simplex $\overline{Wb\square}(n)$. The space $Wb\square(n)$ is defined as a quotient space of  $\overline{Wb\square}(n)$ by the above relations. As a consequence one has an evaluation map
$$
\xi_n^g\colon Wb\square(n)\to\calO(n).
$$
To complete the sequence $Wb\square$, define $Wb\square(0)$ as a point $*$. The map $\xi_0^g$ is defined to have image $\eta(a_0)$.

The figure below describes the quotient map $\overline{Wb\square}(n)\to Wb\square(n)$, $n=1$, 2, 3.

$$
\psfrag{BS1}[0][0][1][0]{$\overline{Wb\square}(1)$}
\psfrag{BS2}[0][0][1][0]{$\overline{Wb\square}(2)$}
\psfrag{BS3}[0][0][1][0]{$\overline{Wb\square}(3)$}
\psfrag{S1}[0][0][1][0]{$Wb\square(1)$}
\psfrag{S2}[0][0][1][0]{$Wb\square(2)$}
\psfrag{S3}[0][0][1][0]{$Wb\square(3)$}
\includegraphics[width=14cm]{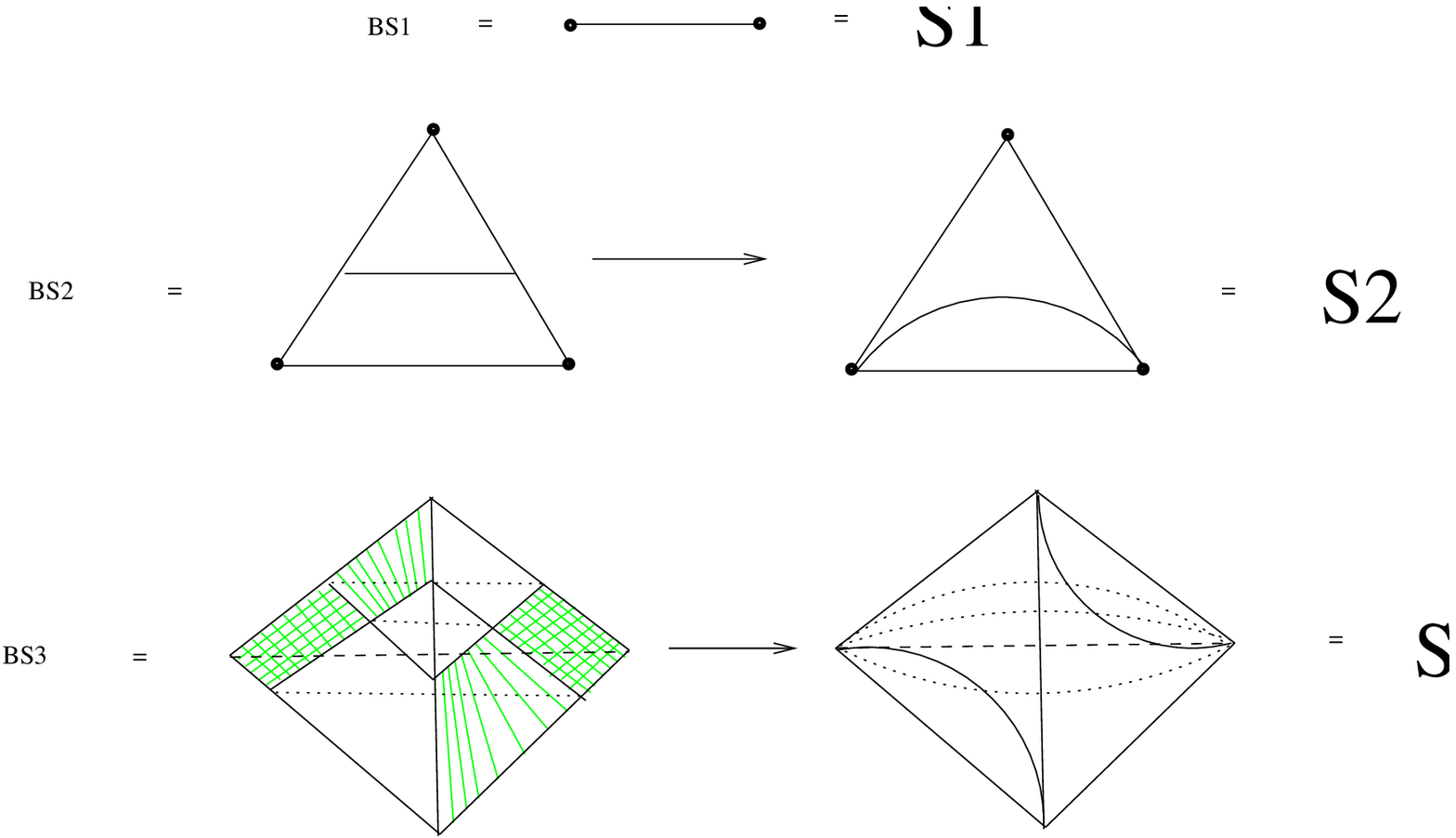}
$$
In $\overline{Wb\square}(3)$ the shaded 2-cells are those  that are  collapsed via the projection to $Wb\square(3)$.

\begin{proposition}\label{p:w_square}
\textup{(i)} The sequence of spaces $Wb\square(n)$, $n\geq 0$, has a structure of a weak bimodule over $\Assoc_{>0}$. Moreover with this structure for any $g\in\underset{\Assoc_{>0}}{\WBimod}(\square,\calO)$, the evaluation map $\xi^g\colon Wb\square\to\calO$ is a morphism of weak bimodules.

\textup{(ii)} The weak bimodule $Wb\square$ is cofibrant and contractible in each degree.

\textup{(iii)} As a weak bimodule $Wb\square$ is homeomorphic to $\triangle$, moreover $Wb\square$ is a refinement of $\triangle$, see Definition~\ref{d:refinement}.
\end{proposition}

\begin{proof}
First we will describe a natural way to encode the cells of $Wb\square(n)$, $n\geq 0$, in terms of trees. A grafting construction similar to the one given in Section~\ref{ss:Triangle}  will be used to define a weak $\Assoc_{>0}$-bimodule structure on $Wb\square$. Let $T\in Obj(\Xi_n)$. We start by encoding faces of $\overline{Wb\square}(T)=\lambda_\square(T)\times\chi_\blacktriangle(T)$. Since by construction $\lambda_\square(T)$ is a face of $\square(n)$, any face of $\lambda_\square(T)$ is also a face of $\square(n)$ and thus can be encoded by the corresponding tree $T'$ (which is an expansion of $T$).
$$
\includegraphics[width=10cm]{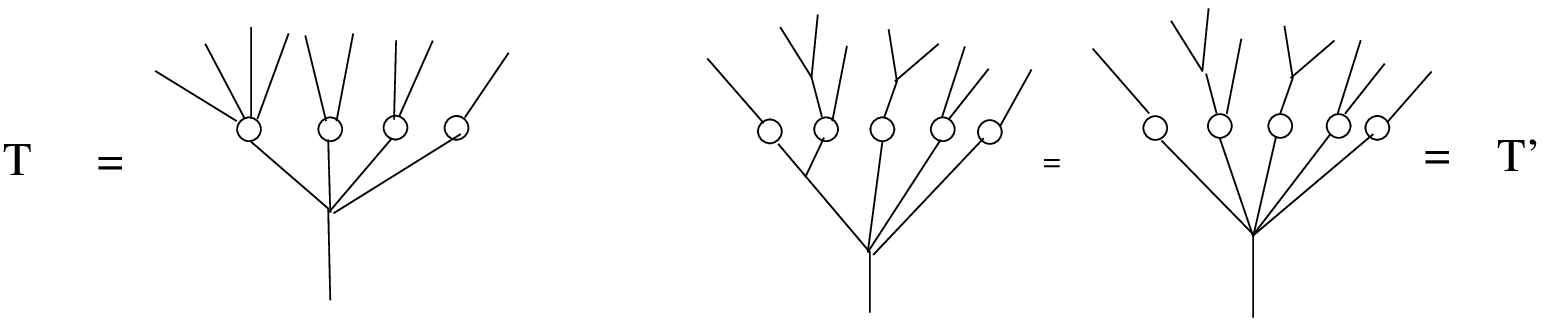}
$$

To recall a point of $\chi_\blacktriangle(T)$ is a configuration $0\leq t_1\leq t_2\leq \ldots\leq t_k\leq 1$, where  each $t_i$ corresponds to a bead of $T$. Thus a face of $\chi_\blacktriangle(T)$ is described by a collision of some points (including collision with 0 or 1). A face of  $\chi_\blacktriangle(T)$ will be encoded by the same tree $T$ in which we encircle by dotted lines the beads whose corresponding parameters collide.  In addition to that we put 0 or 1 on the leftmost, respectively rightmost encircled group of beads if the corresponding parameters collided with~0, respectively~1. Finally to encode a face of
$\lambda_\square(T)\times\chi_\blacktriangle(T)$ we pullback the encircled groups of beads from $T$ to $T'$. (One has a natural projection map of the sets of beads $B(T')\to B(T)$. Each encircled group in $B(T')$ is defined as a preimage of an encircled group in $B(T)$.) We also put 0 or 1 on the corresponding group if it is a preimage of a group with this label. The figure below gives an example of such tree for the case $T$ and $T'$ are as in the figure above, and the face of $\chi_\blacktriangle(T)$ is given by the equations $t_3=t_4=1$.
$$
\includegraphics[width=2.3cm]{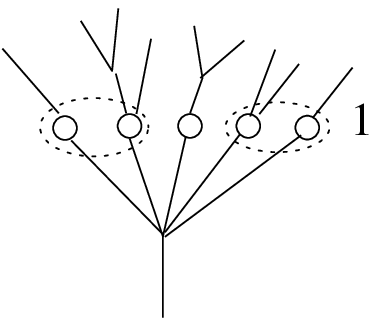}
$$

Such trees will be called trees with encircled beads. It is easy to see that the natural maps $\overline{Wb\square}(T)\to \overline{Wb\square}(n)$ are all inclusions of cell complexes, and moreover for two different trees $T_1,T_2\in Obj(\Xi_n)$ a face of $\overline{Wb\square}(T_1)$ is identified with a face of $\overline{Wb\square}(T_2)$ in $Wb\square(n)$ if and only if they are encoded by the same tree with encircled beads.

Now we will describe how to encode cells of $Wb\square(n)=\overline{Wb\square}(n)/\sim$. Given a tree with encircled beads (that encodes a cell of $\overline{Wb\square}(n)$) it can be of one of the 4 types:

\begin{itemize}
\item there is at least one encircled group not labeled neither by 0, nor by 1;
\item there are exactly two encircled groups: the left one labeled by 0, and the right one labeled by 1;
\item there is only one encircled group which is labeled by 0;
\item there is only one encircled group which is labeled by 1.
\end{itemize}

In the first case the tree corresponding to the image (quotient) cell  in $Wb\square(n)$ is obtained from the tree with encircled beads by turning all the beads in the encircled groups labeled by 0 and 1 into internal vertices and then by contracting all the edges that were adjacent to those vertices:
$$
\includegraphics[width=6cm]{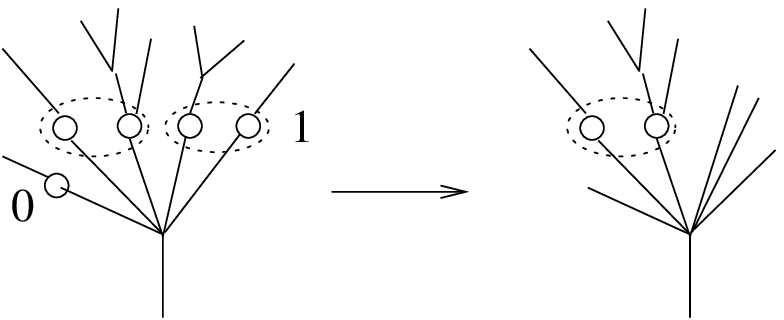}
$$
The leaves that used to come from the encircled group labeled by~0 (respectively~1) will be called \textit{left-most} (respectively \textit{right-most}).

In the  second case the corresponding tree is obtained using the same contraction of edges (that are adjacent to the beads from the group labeled by 0 and 1), but in addition we add a valence one bead to separate the leaves coming from two different groups:
$$
\includegraphics[width=5.5cm]{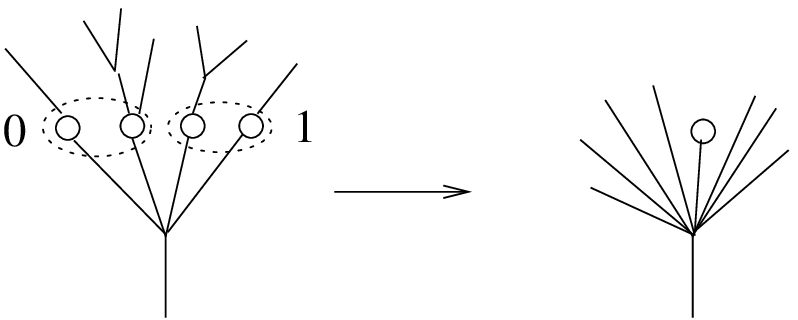}
$$
The corresponding cell is 0-dimensional.

In the third and forth cases the cell is also collapsed to a single point in $Wb\square(n)$. In the first case the corresponding 0-cell will be encoded by the tree as follows:
$$
\includegraphics[width=4.8cm]{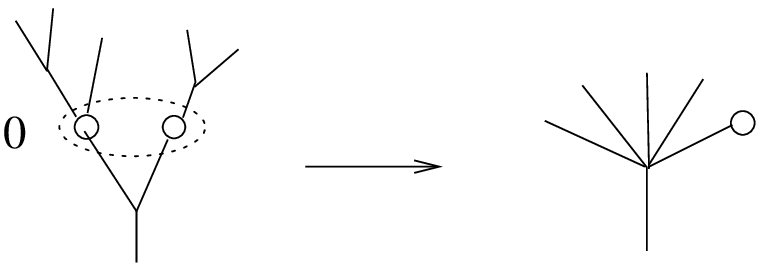}
$$

In the forth case, the corresponding tree will be as follows:
$$
\includegraphics[width=4.8cm]{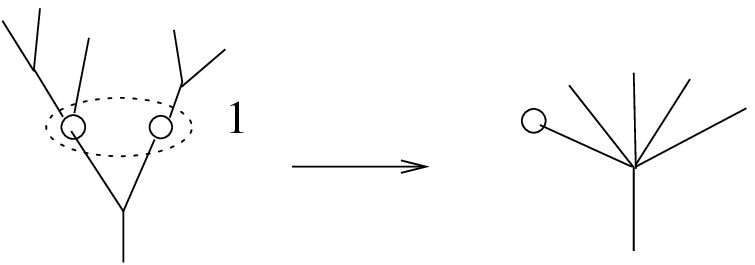}
$$

To recall $Wb\square(0)$ is a point, its only cell will be encoded by \raisebox{-3.5pt}{\includegraphics[width=.15cm]{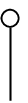}}. The trees encoding cells of $Wb\square(n)$, $n\geq 0$, will be also called trees with encircled beads.

To define a weak $\Assoc_{>0}$-bimodule structure on $Wb\square$ we use the same grafting construction as in Subsection~\ref{ss:Triangle}. It is easy to see that $Wb\square$ is cofibrant. The generating cells are encoded by the trees with encircled beads that do not have leaves coming from the inner vertices.

The map $\xi^g$ respects the right action of $\Assoc_{>0}$ since $g(-,t)$ respects the right $\Assoc_{>0}$ action for all $t\in [0,1]$. (Actually, $\overline{Wb\square}$ has a structure of a right $\Assoc_{>0}$ module and by construction the evaluation map from $\overline{Wb\square}$ also respects this structure.)  We recall that the (weak) left action on $\calO$ arises from a (strong) left action and the morphism of bimodules
$\eta|_{\Assoc_{>0}}\colon\Assoc_{>0}\to\calO$. One can show that the left weak action is preserved by $\xi^g$ because $g(x,0)=g(x,1)=\eta(a_n)$ for any $x\in\square(n)$.

To see that $Wb\square(n)$ is contractible in each degree we notice that $\overline{Wb\square}(n)$ is contractible being an $n$-simplex, and the preimage of any point under the projection $\overline{Wb\square}(n)\to Wb\square(n)$ is contractible. If such point lies in the interior of a cell encoded by a tree with $\ell$ left-most  and $r$ right-most leaves, its preimage is homeomorphic to $\square(\ell)\times\square(r)$, where $\square(0)=*$ as in Subsection~\ref{ss:cof_bim2}.

We will be sketchy in the proof of (iii) since it is enough for our purposes that $Wb\square$ is a cofibrant model of $\Assoc$ in the category $\underset{\Assoc_{>0}}{\WBimod}$, which follows from (i) and (ii). To see (iii) we can argue by induction. Define a filtration
$$
Wb\square_0\subset Wb\square_1\subset Wb\square_2\subset\ldots
\eqno(\numb)\label{eq:filt_w_sq}
$$
in $Wb\square$ similar to~\eqref{eq:triangle_filtration} where $Wb\square_i$ is a weak sub-bimodule of $Wb\square$ generated by $Wb\square(0)$, $Wb\square(1)$, $\ldots$, $Wb\square(i)$. We will show that this filtration of weak bimodules is homeomorphic to~\eqref{eq:triangle_filtration}. By induction hypothesis we can assume that $Wb\square_{N-1}$ is homeomorphic to $\triangle_{N-1}$. This in particular means that  $Wb\square_{N-1}(N)$ is homeomorphic to $\triangle_{N-1}(N)=\partial\triangle(N)\simeq S^{N-1}$. On the other hand $Wb\square(N)$ is obtained from $\overline{Wb\square}(N)$ by a quotient map. Moreover  $\overline{Wb\square}(N)$ is an $N$-simplex and the quotient map in question identifies only the points on the boundary of  the simplex $\overline{Wb\square}(N)$. Thus the preimage of $Wb\square_{N-1}(N)$ via the quotient map $\overline{Wb\square}(N)\to Wb\square(N)$ is exactly $\partial\overline{Wb\square}(N)$. From the above we obtain that $Wb\square(N)$ is obtained from $Wb\square_{N-1}(N)\simeq S^{N-1}$ by attaching an $N$-disc along a map
$$
\partial\overline{Wb\square}(N)=S^{N-1}\to S^{N-1}=Wb\square_{N-1}(N)
\eqno(\numb)\label{eq:W_sq_att_map}
$$
which is of degree~1. (The map~\eqref{eq:W_sq_att_map} is the quotient map $\overline{Wb\square}(N)\to Wb\square(N)$ restricted to the boundary.) Induction step will be proved if we show that $Wb\square(N)$ is homeomorphic to an $N$-disc. Indeed, this would imply that $Wb\square_N|_N$ is homeomorphic to $\triangle_N|_N$ as $N$-truncated weak bimodules. On the other hand the generating cells of $Wb\square_N$ and $\triangle_N$ all lie in degrees $\leq N$, as a consequence of Lemma~\ref{l:trunc_maps} $Wb\square_N$ should be homeomorphic to $\triangle_N$.

\begin{lemma}\label{l:stand_collapse}
Let a product of two closed discs $D^i\times D^j$ be embedded inside the boundary of a closed disc $D^N$ by a standard embedding $i\colon D^i\times D^j\hookrightarrow D^N$, $N>i+j$. Consider the quotient $D^N/\sim$ that identifies $i(x,y)\sim i(x,y')$ for all $x\in D^i$, $y,y'\in D^j$. Then $D^N/\sim$ is homeomorphic to $D^N$.
\end{lemma}

By \lq\lq standard embedding" above we mean an embedding which sends $D^i\times D^j$ to a product of coordinate discs in some coordinate chart of the boundary $\partial D^N=S^{N-1}$.

Even though technical the proof of the above lemma is straightforward by which reason it is omitted here.

As a consequence of this lemma we will have.

\begin{lemma}\label{l:w_sq(T)-disc}
For each $T\in Obj(\Xi_n)$, $Wb\square(T)$ is an $n$-disc.
\end{lemma}

\begin{proof}
$\overline{Wb\square}(T)$ is an $n$-disc since it is a product $\lambda_\square(T)\times\chi_\blacktriangle(T)$ of a cube and of a simplex. On the other hand it is easy to see that $Wb\square(T)$ is obtained from $\overline{Wb\square}(T)$ by a sequence of collapses from Lemma~\eqref{l:stand_collapse} (hint to the reader: the number of collapses needed is twice the number of the beads in $T$).
\end{proof}

From the way $Wb\square(T)$ are glued together, see Proposition~\ref{p:cell_dec_W_Square}, one can see that $Wb\square(n)$ is an $n$-disc, which finishes the proof of Proposition~\ref{p:w_square}~(3).

\end{proof}

\begin{proof}[Proof of Theorem~\ref{t:first_deloop_stages}]
We will induct over $N$. For $N=1$, $\xi_1$ sends $\Omega\, T_1^\square(\calO)=\Omega\,\calO(1)$ homeomorphically to the fiber of the map $T_1^\triangle(\calO)\to T_0^\triangle(\calO)$ over $\eta(a_0)\in\calO(0)=T_0^\triangle(\calO)$. Since the map $T_1^\triangle(\calO)\to T_0^\triangle(\calO)$ is a fibration with contractible base space, one gets that $\xi_1$ is a homotopy equivalence. For $N\geq 2$, we assume that
$$
\xi_{N-1}\colon\Omega\, T_{N-1}^\square(\calO)\to T_{N-1}^\triangle(\calO)
$$
is a homotopy equivalence. Since both maps $T_N^\triangle(\calO)\to T_{N-1}^\triangle(\calO)$ and $\Omega\, T_N^\square\to \Omega\, T_{N-1}^\square(\calO)$ are fibrations, it is enough to check that for any $g\in\Omega\, T_{N-1}^\square(\calO)$ the induced map of fibers is a homotopy equivalence. The fiber for the first map is $\underset{\Assoc_{>0}}\WBimod{}_{\xi_{N-1}^g}\left((\triangle_N,\triangle_{N-1}),\calO\right)=\underset{\Assoc_{>0}}\WBimod{}_{\xi_{N-1}^g}\left((Wb\square_N,Wb\square_{N-1}),\calO\right)$. It turns out that the fiber over the second map can also be expressed as a mapping space of weak bimodules: $\underset{\Assoc_{>0}}\WBimod{}_{\xi_{N-1}^g}\left((Wb\square_N,Wb\square_{N-1/2}),\calO\right)$, where the weak bimodule $Wb\square_{N-1/2}$ is described below.
The induction step will then follow from Lemma~\ref{l:fiber_equiv}.

 To recall $Wb\square_N$ is obtained from $Wb\square_{N-1}$ by a free attachment of an $N$-cell along the map~\eqref{eq:W_sq_att_map}, so that one has the following pushout diagram in the category of weak bimodules over $\Assoc_{>0}$:
$$
\xymatrix{
\WBimod(\partial\overline{Wb\square}(N))\ar[r]\ar[d]&\WBimod(\overline{Wb\square}(N))\ar[d]\\
Wb\square_{N-1}\ar[r]&\text{\makebox[0pt][r]{\raisebox{13pt}[1pt]{$\lrcorner$\hspace{2pt}}}}Wb\square_N,
}
$$
where $\WBimod(\partial\overline{Wb\square}(N))$ (respectively $\WBimod(\overline{Wb\square}(N))$) denotes a free weak $\Assoc_{>0}$ bimodule generated by the sequence that has the sphere $\partial\overline{Wb\square}(N)$ (respectively the ball $\overline{Wb\square}(N)$) in degree $N$ and the empty set in all the other degrees. Let $C_N$ denote the $N$-corolla (the terminal object of the category $\Xi_N$). Consider the interior $Int(\overline{Wb\square}(C_N))$ of $\overline{Wb\square}(C_N)$ viewed as an open subset of $\overline{Wb\square}(N)$. Define $Wb\square_{N-1/2}$ as a weak bimodule that fits into the following pushout diagram:
$$
\xymatrix{
\WBimod(\partial\overline{Wb\square}(N))\ar[r]\ar[d]&\WBimod(\overline{Wb\square}(N)\setminus Int(\overline{Wb\square}(C_N)))\ar[d]\\
Wb\square_{N-1}\ar[r]&\text{\makebox[0pt][r]{\raisebox{13pt}[1pt]{$\lrcorner$\hspace{2pt}}}}Wb\square_{N-1/2}.
}
$$
Again $\WBimod(\overline{Wb\square}(N)\setminus Int(\overline{Wb\square}(C_N)))$ denotes the free weak bimodule over $\Assoc_{>0}$ generated by the sequence with emptysets in all degrees except the degree $N$ where it has the space $\overline{Wb\square}(N)\setminus Int(\overline{Wb\square(C_N)})$. This diagram tells us that $Wb\square_{N-1/2}$ is obtained from $Wb\square_{N-1}$ by a free attachment of a punctured disc. Thus Lemma~\ref{l:fiber_equiv} can be applied and one obtains that the inclusion
$$
\underset{\Assoc_{>0}}\WBimod{}_{\xi_{N-1}^g}\left((Wb\square_N,Wb\square_{N-1/2}),\calO\right)\hookrightarrow
\underset{\Assoc_{>0}}\WBimod{}_{\xi_{N-1}^g}\left((Wb\square_N,Wb\square_{N-1}),\calO\right)
$$
is a homotopy equivalence. On the other hand one can show using Lemma~\ref{l:extension_space} that the left-hand side is naturally homeomorphic to the fiber of the map $\Omega\, T_N^\square\to \Omega\, T_{N-1}^\square(\calO)$ and one also has that the right-hand side is the fiber of the map $T_N^\triangle(\calO)\to T_{N-1}^\triangle(\calO)$. Thus we see that the fibers are homotopy equivalent which finishes the proof of~\ref{t:first_deloop_stages}.
\end{proof}

\subsection{$Wb$-construction and more about $Wb\square_{N-1/2}$}\label{ss:w_construction}
The notation $Wb\square_{N-1/2}$ in the proof above comes from the fact that this weak bimodule is something between $Wb\square_{N-1}$ and $Wb\square_N$: it contains $Wb\square_{N-1}$ and is contained in $Wb\square_N$. It is cofibrant and is weakly equivalent to $Wb\square_{N-1}$ being obtained from the latter one by a free attachment of a punctured cell. In this subsection one gives a more natural interpretation for this intermediate object. This more conceptual interpretation is used to shorten the proof of similar results: Theorems~\ref{t:hom_equiv_towers2},~\ref{t:deloop1_stages_tilde},~\ref{t:deloop2_stages_tilde}.

First, given an $\Assoc_{>0}$-bimodule $\calQ$ one can construct a weak $\Assoc_{>0}$ bimodule $Wb(\calQ)$ similarly as one produces $Wb\square$ from $\square$. One defines first a functor
$$
\lambda_\calQ\colon\Xi_n\to\Top,
$$
by $\lambda_\calQ(T)=\prod_{b\in B(T)}\calQ(|b|)$. On the morphisms the functor is defined using the left action of $\Assoc_{>0}$ on $\calQ$. One defines
$$
Wb(\calQ)(n):=\left(\lambda_\calQ\otimes_{\Xi_n}\chi_\blacktriangle\right)/\sim,
$$
where the relations are~\eqref{eq:W_relation}.
In particular one can apply this construction to any filtration term $\square_N$ of $\square$. The obtained weak bimodule $Wb(\square_N)$ is exactly the part of $Wb\square$ on which the evaluation map $\xi^g$ is defined only using the information about $g_1,\ldots,g_{N}$. To recall $Wb\square_N$ consists of those cells labeled by the trees with encircled beads such that the sum of outgoing edges from all beads is $\leq N$. While the cells of $Wb(\square_N)$ are labeled by the trees with encircled beads such that the number of outgoing edges from each bead is $\leq N$. Thus one has a natural inclusion of cofibrant weak bimodules $Wb\square_N\subset Wb(\square_N)$.

One can notice that $Wb\square_{N-1/2}$ is exactly the part of $Wb\square_N$, where $\xi^g$ is defined using only the information about $g_1,\ldots,g_{N-1}$ and not about $g_N$. In other words $Wb\square_{N-1/2}$ can be obtained as intersection 
$
Wb\square_{N-1/2}=Wb\square_N\cap Wb(\square_{N-1}).
$


\section{Second delooping}\label{s:second_deloop}
Let $\calO$ be an operad such that $\calO(1)\simeq *$. We will also assume that $\calO$ is endowed with a morphism of operads $\eta\colon\Assoc_{>0}\to\calO$.
Because of this map $\calO$ can be regarded as a bimodule over $\Assoc_{>0}$. Denote by $p$ the projection
$
p\colon\Pentagon\to\Assoc_{>0}.
$
By $\Omega\,\Operad(\Pentagon,\calO)$ we will understand the loop space with the base point $\eta\circ p$. Similarly $\Omega\, T_N^{\pentagon}(\calO)=\Omega\, \Operad_N(\Pentagon|_N,\calO|_N)$ will denote the loop space with the base point $\eta\circ p|_N$. A point in
$\Omega\, T_N^{\pentagon}(\calO)$ is a tuple $(h_2,\ldots,h_N)$, where each $h_i$ is a map
$$
h_i\colon [0,1]\times\Pentagon(i)\to\calO(i),
$$
satisfying the boundary condition~\eqref{eq:pentag_bound} for all $t\in [0,1]$. We will also have
$$
h_i(0,-)=h_i(1,-)=\eta\circ p|_{\pentagon(i)}.
\eqno(\numb)\label{eq:eta_p_ends}
$$
This means that each map $h_i$ can be viewed as a map of suspension $h_i\colon\Sigma\Pentagon(i)\to\calO(i)$.

The aim of this section is to construct a homotopy equivalence of towers:
$$
\xymatrix{
\Omega\, T_0^{\pentagon}(\calO)\ar[d]_{\zeta_0}&\Omega\, T_1^{\pentagon}(\calO)\ar[d]_{\zeta_1}\ar[l]&\Omega\, T_2^{\pentagon}(\calO)\ar[d]_{\zeta_2}\ar[l]&\Omega\, T_3^{\pentagon}(\calO)\ar[d]_{\zeta_3}\ar[l]&\ldots\ar[l]\\
T_0^\square(\calO)&T_1^\square(\calO)\ar[l]&T_2^\square(\calO)\ar[l]&T_3^\square(\calO)\ar[l]&\ldots\ar[l],
}
\eqno(\numb)\label{eq:map2_of_towers}
$$
that induces an equivalence of their homotopy limits
$$
\Omega\Operad(\Pentagon,\calO)\stackrel{\zeta_\infty}{\longrightarrow}\underset{\Assoc_{>0}}{\Bimod}(\square,\calO).
\eqno(\numb)\label{eq:zeta_infty}
$$

\begin{theorem}\label{t:hom_equiv_towers2}
Let $\calO$ be an operad endowed with a map of operads $\eta\colon\Assoc_{>0}\to \calO$ and with $\calO(1)\simeq *$.  For each $N\geq 0$ the map $\zeta_N\colon\Omega\, T_N^{\pentagon}(\calO)\to T_N^\square(\calO)$ constructed below is a homotopy equivalence.
\end{theorem}

As a consequence we immediately obtain the following.

\begin{theorem}\label{t:hom_equiv2}
Under the settings of Theorem~\ref{t:hom_equiv_towers2}
the map~\eqref{eq:zeta_infty} is a homotopy equivalence.
\end{theorem}

The proof of the latter theorem  is analogous to the proof of Theorem~\ref{t:first_deloop}.

The construction of the map $\zeta_\bullet$ is similar to that of $\xi_\bullet$. The idea is to decompose each $\square(N)$, $N\geq 1$, into a union of polytopes (of the same dimension $N-1$), one of which is $\Sigma\Pentagon(N)$ and the other ones all together define a homotopy between $\partial\Sigma\Pentagon(N)$ and $\partial\square(N)$. On each subpolytop a map to $\calO(N)$ will be naturally defined using the maps $h_2,\ldots, h_N$.

\subsection{Bimodule {\sc$B\Pentagon$} over {\sc$\Pentagon$}}\label{ss:PB}
Let $\Psi_n$ denote the category/poset of faces of $\Pentagon(n)$. Its objects are planar trees and morphisms are contractions of edges, see Section~\ref{s:operads}. We will define a covariant functor
$$
\lambda_{\pentagon}\colon\Psi_n\to\Top,
$$
and a contravariant one:
$$
\chi_\blacktriangle\colon\Psi_n\to \Top.
$$
The functor $\lambda_{\pentagon}$ assigns to a tree $T$ the face of $\Pentagon(n)$ encoded by $T$:
$$
\lambda_{\pentagon}(T)=\prod_{b\in B(T)}\Pentagon(|b|).
$$
The contravariant functor $\chi_\blacktriangle$ assigns to a tree $T$ a polytope in $[0,1]^{B(T)}$ whose points are the order preserving maps $B(T)\to [0,1]$, where $B(T)$ has a natural partial order induced by the structure of the tree. For example
$$
\chi_\blacktriangle\Bigl(\raisebox{-15pt}{\includegraphics[width=1.5cm]{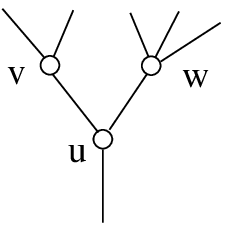}}\Bigr)=\{(t_u,t_v,t_w)\, |\, 0\leq t_u\leq t_v\leq 1,\, 0\leq t_u\leq t_w\leq 1\}.
\eqno(\numb)\label{eq:uvw}
$$
As we see from this example $\chi_\blacktriangle(T)$ might be different from a simplex. Given a morphism $\alpha\colon T_1\to T_2$ in $\Psi_n$ one has the induced map  of the sets of beads $\alpha_*\colon B(T_1)\to B(T_2)$. The map
$$
\chi_\blacktriangle(\alpha)\colon\chi_\blacktriangle(T_2)\to\chi_\blacktriangle(T_1)
$$
is a face inclusion defined by $f\mapsto f\circ\alpha_*$, for $f\in \chi_\blacktriangle (T_2)$ viewed as a map $f\colon B(T_2)\to [0,1]$.

Denote by $\overline{B\Pentagon}(T)=\lambda_{\pentagon}(T)\times\chi_\blacktriangle(T)$. A point in $\overline{B\Pentagon}(T)$ is a tuple $(x_b,t_b)_{b\in B(T)}$, where $x_b\in\Pentagon(|b|)$, $t_b\in [0,1]$. One can define a natural evaluation map
$$
\bar\zeta^h_T\colon \overline{B\Pentagon}(T)\to \calO(n)
\eqno(\numb)\label{eq:bar_zeta_T}
$$
by sending $(x_b,t_b)_{b\in B(T)}$ to the composition of the elements $h_{|b|}(t_b,x_b)\in\calO(|b|)$ according to the way they appear in the tree $T$. For example for the tree from~\eqref{eq:uvw} one has
$$\bar\zeta^h_T\bigl( (t_u,x_u), (t_v,x_v), (t_w,x_w)\bigr) = h_2(t_u,x_u)\bigl( h_2(t_v,x_v),h_3(t_w,x_w)\bigr).$$

Now define
$$
\overline{B\Pentagon}(n):=\lambda_{\pentagon}\otimes_{\Psi_n}\chi_\blacktriangle.
\eqno(\numb)\label{eq:overline_B_Pentagon}
$$
Notice that the natural map $\overline{B\Pentagon}(T)\to\overline{B\Pentagon}(n)$ is inclusion of cell complexes for each $T\in Obj(\Psi_n)$. Moreover
$\overline{B\Pentagon}(n)$ can be viewed as a union of polytopes $\overline{B\Pentagon}(T)$ (glued together by their faces) of the same dimension $n-1$.

The fact that $h(t,-)\colon\Pentagon\to\calO$ is a morphism of operads for all $t\in [0,1]$ implies that the evaluation maps~\eqref{eq:bar_zeta_T} can be glued together to define a map
$$
\bar\zeta^h_n\colon\overline{B\Pentagon}(n)\to\calO(n).
$$
We will also define $\overline{B\Pentagon}(1)$ to be a point. The evaluation map $\bar\zeta^h_1$ is defined to have image $\eta(a_1)$.

\begin{proposition}\label{p:over_B_Pent_disc}
The space $\overline{B\Pentagon}(n)$, $n\geq 1$, is homeomorphic to an $(n-1)$-disc.
\end{proposition}

\begin{proof}
$\overline{B\Pentagon}(n)$ is a union of $(n-1)$ dimensional polytopes $\overline{B\Pentagon}(T)$, $T\in Obj(\Psi_n)$. Indeed,
$$
\dim\lambda_{\pentagon}(T)=\sum_{b\in B(T)}(|b|-2)=\sum_{b\in B(T)}(|b|-1)-\# B(T)=n-1-\#B(T);
$$
$$
\dim \chi_\blacktriangle(T)=\# B(T).
$$
So, $\dim \overline{B\Pentagon}(n)=n-1$. Define a filtration in $\overline{B\Pentagon}(n)$:
$$
\emptyset=\overline{B\Pentagon}^0(n)\subset \overline{B\Pentagon}^1(n)\subset \overline{B\Pentagon}^2(n)\subset\ldots\subset \overline{B\Pentagon}^{n-1}(n)=
\overline{B\Pentagon}(n),
\eqno(\numb)\label{eq:filtr_over_B_Pent}
$$
where $\overline{B\Pentagon}^i(n)$ is the union of those $\overline{B\Pentagon}(T)$ for which $T$ has $\leq i$ beads.\footnote{We used upper indices to distinguish this filtration from the filtration~\eqref{eq:filt_b_pent}.} One has $\overline{B\Pentagon}^1(n)=
\overline{B\Pentagon}(C_n)=\Pentagon(n)\times [0,1]$ is an $(n-1)$-disc.\footnote{As before $C_n$ is the $n$-corolla, which is the terminal object in $\Psi_n$.} We
 will induct over $i$ to show that each $\overline{B\Pentagon}^i(n)$, $1\leq i\leq n-2$ is an $(n-1)$-disc. The poset $\Psi_n$ has a structure of a lattice with
 $T_1\vee T_2$ defined as a tree corresponding to the smallest face of $\Pentagon(n)$ that contains both $T_1$ and $T_2$ faces, and $T_1\wedge T_2$ defined as the tree corresponding to the
intersection of the $T_1$ and $T_2$ faces. It is easy to see that
 $\overline{B\Pentagon}(T_1)\cap \overline{B\Pentagon}(T_2)=\lambda_{\pentagon{}}(T_1\wedge T_2)\times\chi_\blacktriangle(T_1\vee T_2)\subset \overline{B\Pentagon}(T_1\vee T_2)$. This implies that if $T_1$ and $T_2$ both have $i$ beads then the intersection $\overline{B\Pentagon}(T_1)\cap \overline{B\Pentagon}(T_2)$ lies in the smaller filtration term $\overline{B\Pentagon}^{i-1}(n)$.
One can also notice that $\overline{B\Pentagon}(T_1)$ shares a facet (codimension one faces) with
 $\overline{B\Pentagon}(T_2)$ if and only if one tree ($T_1$ or $T_2$) is obtained from the other by a contraction of exactly one edge.
  Thus if $T_1$ and $T_2$ both have $i$ beads then $\overline{B\Pentagon}(T_1)\cap \overline{B\Pentagon}(T_2)$ lies in the lower term of filtration and has dimension $\leq n-3$. To finish the
 proof it is enough to show that for any tree $T$ with $i$ beads $\overline{B\Pentagon}(T)\cap\overline{B\Pentagon}_{i-1}(n)$ is an $(n-2)$-disc. But one can see
 that this intersection is the union of those facets of $\overline{B\Pentagon}(T)$ that contain the face
$\overline{B\Pentagon}(T)\cap \overline{B\Pentagon}(C_n)$. (The latter face of intersection is $\lambda_{\pentagon}(T)\times\chi_\blacktriangle(C_n)=\lambda_{\pentagon}(T)\times [0,1]$ --- it is a codimension $(i-1)$ face of $\overline{B\Pentagon}(C_n)$.) As a consequence this intersection is an $(n-2)$-disc.

%
%
%
%
%
%
%
%
%
\end{proof}

It turns out that the spaces $\overline{B\Pentagon}(n)$ , $n\geq 1$, form a bimodule over the operad $\Pentagon$. Below we make this structure explicit. The left action
$$
\Pentagon(n)\times \overline{B\Pentagon}(m_1)\times \ldots \times\overline{B\Pentagon}(m_n)\to \overline{B\Pentagon}(m_1+\ldots +m_n)
$$
is defined by gluing together the maps
$$
\Pentagon(n)\times \overline{B\Pentagon}(T_1)\times\ldots \times\overline{B\Pentagon}(T_n)\to \overline{B\Pentagon}(gr(T_1,\ldots,T_n)),
\eqno(\numb)\label{eq:over_B_pent_comp}
$$
where $T_i\in Obj(\Psi_{m_i})$, $1\leq i\leq n$, and $gr(T_1,\ldots,T_n)$ is a tree obtained from $T_1,\ldots,T_n$ by grafting grafteng them to the $n$-corolla $C_n$:
$$
gr(T_1,\ldots,T_n)=
\raisebox{-30pt}{\psfrag{T1}[0][0][1][0]{$T_1$}
\psfrag{T2}[0][0][1][0]{$T_2$}
\psfrag{Tn}[0][0][1][0]{$T_n$}
\includegraphics[width=3.4cm]{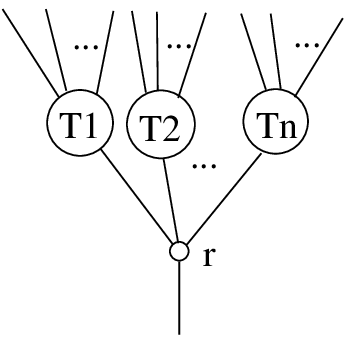}}
$$
Denote the lowest bead of $gr(T_1,\ldots,T_n)$ by $r$. Let $(\tilde t_b, \tilde x_b)_{b\in B(gr(T_1,\ldots,T_n))}$ be coordinates in
$\overline{B\Pentagon}(gr(T_1,\ldots,T_n))$ (notice that $B(gr(T_1,\ldots,T_n))=\sqcup_{i=1}^nB(T_i)\sqcup\{r\}$). Let $y$ be a point of $\Pentagon(n)$, and $(t_b,x_b)_{b\in \amalg_{i=1}^nB(T_i)}$ be a point  of $\prod_{i=1}^n \overline{B\Pentagon}(T_i)$. In coordinates the map~\eqref{eq:over_B_pent_comp} is given by
$$
\tilde t_b=
\begin{cases}
t_b,&b\neq r;\\
0,&b=r;
\end{cases}
$$
$$
\tilde x_b=
\begin{cases}
x_b,&b\neq r;\\
y,& b=r.
\end{cases}
$$

The right action
$$
\overline{B\Pentagon}(n)\times\Pentagon(m_1)\times\ldots\times\Pentagon(m_n)\to\overline{B\Pentagon}(m_1+\ldots +m_n)
$$
is defined similarly by gluing together the maps
$$
\overline{B\Pentagon}(T)\times\Pentagon(m_1)\times\ldots\times\Pentagon(m_n)\to\overline{B\Pentagon}(T(m_1,\ldots,m_n)),
\eqno(\numb)\label{eq:over_B_pent_comp2}
$$
where $T\in Obj(\Psi_n)$ and the tree $T(m_1,\ldots,m_n)$ is obtained by grafting the corollas $C_{m_1},\ldots C_{m_n}$ to $T$:
$$
T(m_1,\ldots,m_n)=
\raisebox{-30pt}{\psfrag{T}[0][0][1][0]{$T$}
\includegraphics[width=3cm]{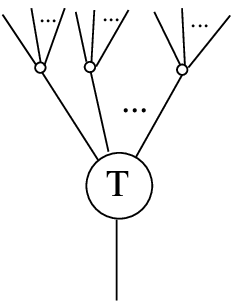}}
$$
Let us denote by $r_1,\ldots,r_n$ the new beads  of $T(m_1,\ldots,m_n)$ (i.e. those that do not come from $T$).
Let $(\tilde t_b,\tilde x_b)_{b\in B(T(m_1,\ldots,m_n))}$ be coordinates on  $\overline{B\Pentagon}(T(m_1,\ldots,m_n))$ (notice that $B(T(m_1,\ldots,m_n))=B(T)\sqcup\{r_1,\ldots,r_n\}$). Let also $y_i\in\Pentagon(m_i)$, $1\leq i\leq n$, and $(t_b,x_b)_{b\in B(T)}\in \overline{B\Pentagon}(T)$. The map~\eqref{eq:over_B_pent_comp2} in coordinates is given by
$$
\tilde t_b=
\begin{cases}
t_b,& b\in B(T);\\
1,& b=r_i,\, i=1\ldots n;
\end{cases}
$$
$$
\tilde x_b=
\begin{cases}
x_b,& b\in B(T);\\
y_i,& b=r_i,\, i=1\ldots n.
\end{cases}
$$

It is straightforward to see that the maps as above define a structure of a bimodule over $\Pentagon$ on $\overline{B\Pentagon}$.

\begin{remark}\label{r:multiplihedra}
One can show that $\overline{B\Pentagon}$ is a cofibrant bimodule over $\Pentagon$.  There is a well known sequence of polytopes called multiplihedra which defines a cofibrant bimodule over the Stasheff operad $\Pentagon$~\cite{Stasheff-HHPV,IwMi,Forcey}. One can show that $\overline{B\Pentagon}$ is a refinement (see Definition~\ref{d:refinement}) of the sequence of multiplihedra, in particular the two bimodules are homeomorphic.
\end{remark}

We finish this section by describing the cellular structure of $\overline{B\Pentagon}(n)$, $n\geq 0$, from which it will be clear that the bimodule action of $\Pentagon$ on $\overline{B\Pentagon}$ is free. We start by encoding the faces of $\overline{B\Pentagon}(T)=\lambda_{\pentagon}(T)\times\chi_\blacktriangle(T)$, where $T\in Obj(\Psi_n)$. Since $\lambda_{\pentagon}(T)$ is a face of $\Pentagon(n)$, any face of $\lambda_{\pentagon}(T)$ will be encoded by a tree $T'$ which is an expansion of $T$. A face of $\chi_\blacktriangle(T)$ can be encoded by a tree $T$ in which one encircles the beads whose parameters coincide (an encircled group forms always a subtree of~$T$). In addition if in one of the groups all the parameters are zeros or ones, one puts label 0, respectively 1 on this group:
$$
\includegraphics[width=4cm]{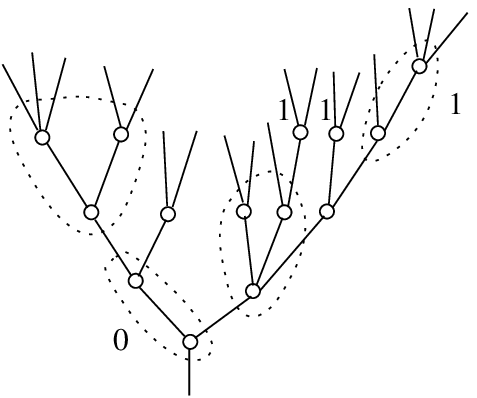}
$$
Finally to encode a face of $\lambda_{\pentagon}(T)\times\chi_\blacktriangle(T)$ which is a product of a face $T'$ and a face encoded by the tree $T$ with encircled beads, we incorporate this information in the tree $T'$ by encircling the groups of vertices of $T'$ which are pre-images of the corresponding encircled groups in $T$. We also put labels 0 and 1 appropriately. It is easy to see that the faces of $\overline{B\Pentagon}(T)$, $T\in Obj(\Psi_n)$, that are identified in $\overline{B\Pentagon}(n)$ are encoded by the same trees with encircled beads.

\subsection{Bimodule {\sc $B\Pentagon$} over $\Assoc_{>0}$}\label{ss:bimod_B_Pentagon}
In this section we define a quotient of  $\overline{B\Pentagon}$ that will naturally be  a cofibrant model of $\Assoc_{>0}$ as a bimodule over $\Assoc_{>0}$.
Let $T\in Obj(\Psi_n)$, we define $B\Pentagon(T)$ as the quotient of $\overline{B\Pentagon}(T)$ by the following relation:
$$
(t_b,x_b)_{b\in B(T)}\sim (t'_b,x'_b)_{b\in B(T)}
\eqno(\numb)\label{eq:relation_BPent}
$$
if and only if for all $b\in B(T)$ one has either $(t_b,x_b)=(t'_b,x'_b)$ or $t_b=t'_b=\epsilon$ with $\epsilon\in\{0,1\}$.

The cell which is the image of a cell encoded by a tree $T_0$ with encircled beads will be encoded by a similar tree with encircled beads obtained from $T_0$ by the following procedure:

1. contract all the edges inside encircled groups labeled by 0 and 1, and replace each bead from this group by a usual internal vertex;

2. put a bead on each edge connecting two internal vertices, one of which is labeled by 0 and another by 1; similarly one puts a bead on an edge connecting an internal vertex labeled by 0 and a leaf.

Here are two examples how this procedure goes:
$$
\includegraphics[width=14cm]{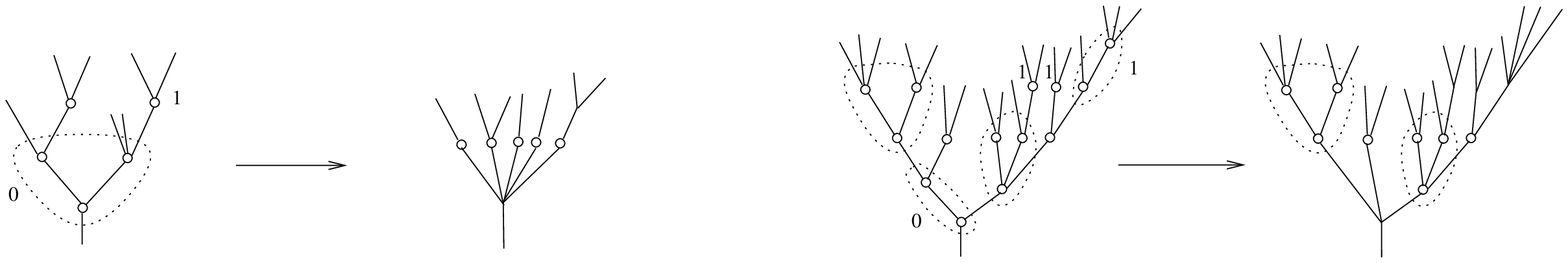}
$$

The space $B\Pentagon(n)$ is obtained from $\overline{B\Pentagon}(n)$ as a quotient by the above equivalence relations ($\overline{B\Pentagon}(n)$ being viewed as a union of $\overline{B\Pentagon}(T)$, $T\in Obj(\Psi_n)$). It is easy to see that two cells of $B\Pentagon(T_1)$ and of $B\Pentagon(T_2)$ are identified in $B\Pentagon(n)$ if and only if they are encoded by the same tree with encircled beads.

It is easy to see that for any $T\in Obj(\Psi_n)$ the evaluation map $\bar\zeta_T^h\colon \overline{B\Pentagon}(T)\to \calO$ factors through $B\Pentagon(T)$. As a consequence the evaluation map $\bar\zeta_n^h\colon \overline{B\Pentagon}(n)\to \calO$ factors through $B\Pentagon(n)$:

$$
\xymatrix{
\overline{B\Pentagon}(n)\ar[rr]^{\bar\zeta^h_n}\ar[rd]&&\calO(n)\\
&B\Pentagon(n)\ar[ru]_{\zeta^h_n}&
}
$$

\begin{proposition}\label{p:B_Pentagon_propert}
\textup{(i)} $B\Pentagon$ has a structure of a cofibrant $\Assoc_{>0}$-bimodule. With this structure the evaluation map $\zeta^h\colon B\Pentagon\to\calO$ defined above is a morphism of bimodules.

\textup{(ii)} Each component of $B\Pentagon$ is contractible.

\textup{(iii)} $B\Pentagon$ is homeomorphic to $\square$ as an $\Assoc_{>0}$-bimodule.
\end{proposition}

\begin{proof}
(i). The $\Assoc_{>0}$ bimodule structure on $B\Pentagon$ is defined by a grafting construction similar to the one described in Subsection~\ref{ss:cof_bim1}. With this definition the projection $\overline{B\Pentagon}\to B\Pentagon$ is a map of $\Pentagon$-bimodules (where $B\Pentagon$ is acted on by $\Pentagon$ through $\Assoc_{>0}$). It is easy to see that the bimodule action of $\Assoc_{>0}$  on $B\Pentagon$ is free. The generating cells are encoded by the trees with encircled beads that do not have internal vertices (that were labeled by 0 and 1). The fact that the evaluation map respects the bimodule structure follows from the construction.

(ii). We have seen that $\overline{B\Pentagon}(n)$ is an $(n-1)$-disc. On the other hand  the preimage of any point under the projection $\overline{B\Pentagon}(n)\to B\Pentagon(n)$ is a product of associahedra, thus is contractible. As a consequence $B\Pentagon(n)$ is homotopy equivalent to
$\overline{B\Pentagon}(n)$ and is also contractible.

(iii). Using Lemma~\ref{l:stand_collapse} one can prove that each $B\Pentagon(T)$, $T\in\Psi_n$, is an $(n-1)$-disc (hint: the number of collapses needed is twice the number of the beads of $T$). From the way how the polytopes  $B\Pentagon(T)$, $T\in\Psi_n$, are glued together (see the proof of Proposition~\ref{p:over_B_Pent_disc}) one gets that $B\Pentagon(n)$ is also an $(n-1)$-disc. Consider a filtration in $B\Pentagon$
$$
\emptyset= B\Pentagon_0\subset B\Pentagon_1\subset B\Pentagon_2\subset\ldots,
\eqno(\numb)\label{eq:filt_b_pent}
$$
where the term $B\Pentagon_n$ is a subbimodule of $B\Pentagon$ generated by $B\Pentagon(1),\ldots,B\Pentagon(n)$. Using the induction and the fact that each $B\Pentagon(n)$ is an $(n-1)$-disc one can show existence of a homeomorphism $i_n\colon\square_n\cong B\Pentagon_n$ that extends $i_{n-1}$. The argument is similar to the proof of Proposition~\ref{p:w_square}.
\end{proof}

Concerning  the proof above we mention that $\overline{B\Pentagon}(0)$ is empty, $\overline{B\Pentagon}(1)$ is a point, $\overline{B\Pentagon}(2)$ is segment and the maps
$\overline{B\Pentagon}(i)\to B\Pentagon(i)$, $i=0$, 1, 2,   are homeomorphisms. The quotient map $\overline{B\Pentagon}(3)\to B\Pentagon(3)$ is described by the following picture:

$$
\psfrag{BBP}[0][0][1][0]{$\overline{B\Pentagon}(3)$}
\psfrag{BP}[0][0][1][0]{$B\Pentagon(3)$}
\includegraphics[width=13cm]{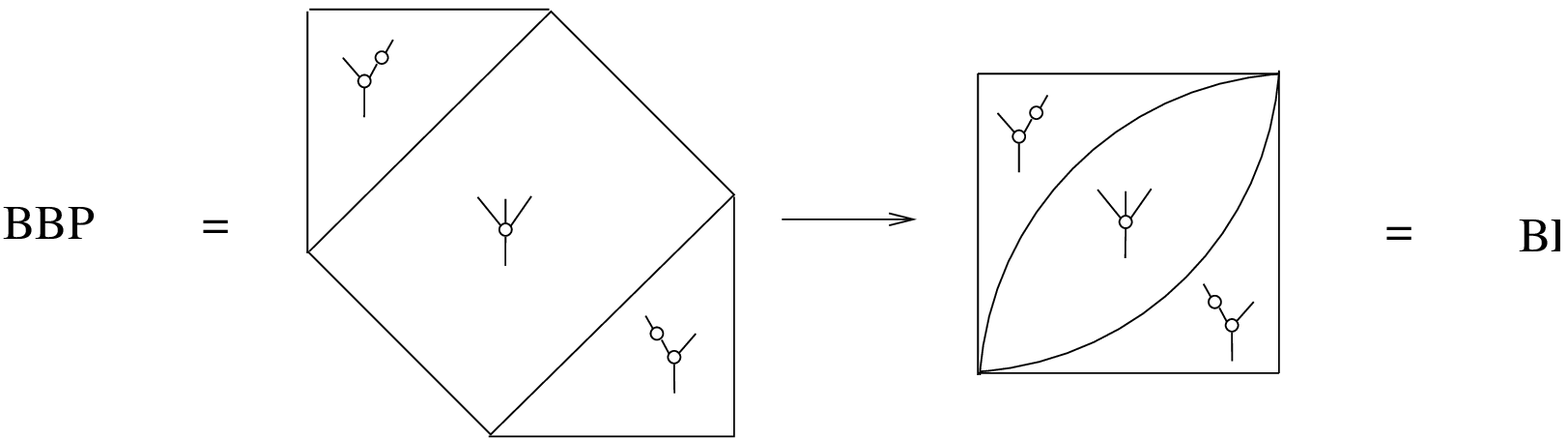}
$$
On the figure above a tree $T$ refers to the piece $\overline{B\Pentagon}(T)$ in $\overline{B\Pentagon}(3)$ or the piece $B\Pentagon(T)$ in $B\Pentagon(3)$, respectively.

The quotient map $\overline{B\Pentagon}(4)\to B\Pentagon(4)$ is harder to visualize. Both the source and the target are unions of 10 pieces labeled by the trees from $\Psi_4$ --- the category of faces of $\Pentagon(4)$. As a little demonstration of  what is going  the figure below shows the cell decomposition of the boundary of $B\Pentagon(4)$.
$$
\includegraphics[width=5cm]{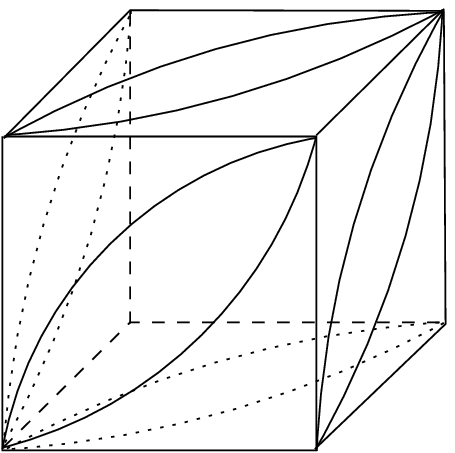}
$$
Notice that the back face consists of only one 2-cell. Under the homeomorphism $B\Pentagon(4)\to\square(4)$ this face gets mapped to the face labeled by the tree \raisebox{-7pt}{\includegraphics[width=.9cm]{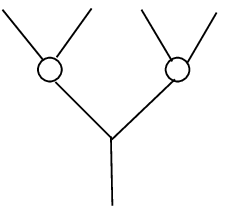}}. All the other faces of the cube above have a cell decomposition homeomorphic to that of $B\Pentagon(3)$. The  trees labeling the corresponding faces of $\square(4)$ are \raisebox{-7pt}{\includegraphics[width=.8cm]{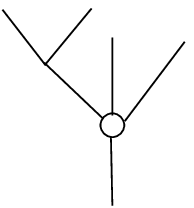}}, \raisebox{-7pt}{\includegraphics[width=.8cm]{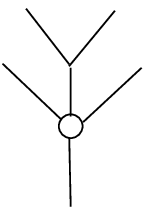}}, \raisebox{-7pt}{\includegraphics[width=.9cm]{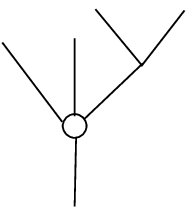}}, \raisebox{-7pt}{\includegraphics[width=.9cm]{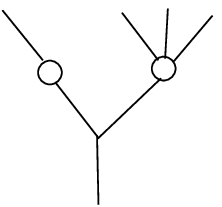}}, \raisebox{-7pt}{\includegraphics[width=.9cm]{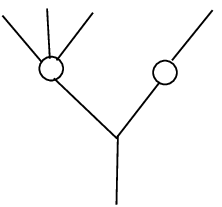}}.

\subsection{$B$-construction and Proof of Theorem~\ref{t:hom_equiv_towers2}}\label{ss:b_construction}
One can notice that for any topological operad $\calP$ with $\calP(0)=\emptyset$  one can define an $\Assoc_{>0}$ bimodule $B\calP$ similarly as one produces $B\Pentagon$ from $\Pentagon$. First one defines a functor $\lambda_\calP\colon\Psi_n\to\Top$ that assigns $T\mapsto\prod_{b\in B(T)}\calP(|b|)$ and uses the operadic composition structure for morphisms. The $n$-th component is defined as
$$
B\calP(n):=\left(\lambda_\calP\otimes_{\Psi_n}\chi_\blacktriangle\right)/\sim,
$$
where the equivalence relation is~\eqref{eq:relation_BPent}.

In particular for $\calP=\Pentagon_N$ the bimodule $B(\Pentagon_N)$ is the subbimodule of $B\Pentagon$ consisting of cells labeled by the trees with encircled beads whose all beads have at most $N$ outgoing vertices. By contrast  the cells of the bimodule $B\Pentagon_N$ are labeled by trees with encircled beads which are product of trees with the property that the sum of outgoing vertices from all beads minus the number of beads is $\leq N$. One obviously has $B\Pentagon_N\subset B(\Pentagon_N)$.

\begin{proof}[Proof of Theorem~\ref{t:hom_equiv_towers2}]
The proof is completely analogous to that of Theorem~\ref{t:first_deloop_stages}. Again we induct over $N$. The map $\zeta_1$ sends $\Omega\, T_1^{\pentagon}(\calO)$
which is a point to $\eta(a_1)\in\calO(1)=T_1^\square(\calO)$. Since $\calO(1)$ is contractible $\zeta_1$ is a homotopy equivalence. The map $\zeta_2$   sends $\Omega\, T_2^{\pentagon}(\calO)=\Omega\calO(2)$  homeomorphically to the preimage of the fibration $T_2^\square(\calO)\to T_1^\square(\calO)$ over $\eta(a_1)$.
Since $\calO(1)=T_1^\square(\calO)$ is contractible, the map $\zeta_2$ is a homotopy equivalence.

Assuming that $\zeta_{N-1}$ is a homotopy equivalence, we will prove that $\zeta_N$ is so. Since both maps $T_N^\square(\calO)\to T_{N-1}^\square(\calO)$ and $\Omega\, T_N^{\pentagon}(\calO)\to \Omega\, T_{N-1}^{\pentagon}(\calO)$ are fibrations, it is enough to check that for any $h\in\Omega\, T_{N-1}^{\pentagon}(\calO)$ the induced map of fibers is a homotopy equivalence. The fiber for the first map is the space of morphisms $\underset{\Assoc_{>0}}\Bimod{}_{\zeta_{N-1}^g}\left((\square_N,\square_{N-1}),\calO\right)=\underset{\Assoc_{>0}}
\Bimod{}_{\zeta_{N-1}^g}\left((B\Pentagon_N,B\Pentagon_{N-1}),\calO\right)$. Using Lemma~\ref{l:extension_space} one can show that the fibers of the second map can also be expressed as  spaces of maps of bimodules $\underset{\Assoc_{>0}}\Bimod{}_{\zeta_{N-1}^g}\left((B\Pentagon_N,B\Pentagon_{N-1/2}),\calO\right)$, where $B\Pentagon_{N-1/2}$ is the bimodule obtained as intersection of $B\Pentagon_{N}$ and the $B$-construction of $\Pentagon_{N-1}$:
$$
B\Pentagon_{N-1/2}=B\Pentagon_{N}\cap B(\Pentagon_{N-1}).
$$
But notice that $B\Pentagon_{N}$ is obtained from $B\Pentagon_{N-1}$ by a free attachment of an $N$-disc in degree $N$ (whose interior is the interior of $B\Pentagon(N)$). While $B\Pentagon_{N-1/2}$ is obtained from $B\Pentagon_{N-1}$ by a free attachment of a punctured disc $B\Pentagon(N)\setminus Int(B\Pentagon(C_N))$.  Lemma~\ref{l:fiber_equiv} finishes the proof.

\end{proof}

\part{Delooping taking into account degeneracies}\label{part2}

In Part~\ref{part1} the partial semicosimplicial totalizations $\hosTot_N\calO(\bullet)$, $N\geq 0$, of a cosimplicial space $\calO(\bullet)$ obtained from a multiplicative operad are described as dooble loop spaces of the spaces of derived morphisms of truncated operads:
$$
\hosTot_N\calO(\bullet)\simeq\Omega^2\widetilde{\Operad}_N\left(\Assoc_{>0}|_N,\calO|_N\right).
$$
In part~\ref{part2} we take into account degeneracies and deloop the partial homotopy totalizations  $\hoTot_N\calO(\bullet)$:
$$
\hoTot_N\calO(\bullet)\simeq\Omega^2\widetilde{\Operad}_N\left(\Assoc|_N,\calO|_N\right).
$$
The constructions in Part~\ref{part2} generalize those from Part~\ref{part1}. Below we outline the plan of Part~\ref{part2} together with the major differences compared to Part~\ref{part1}. In Sections~\ref{s:deg_wb}, \ref{s:cof_tilde_sq}, \ref{s:second_deloop_tilde} respectively, we construct  cofibrant models $\hoTriangle$, $\hoSquare$, and $\hoPentagon$  of $\Assoc$ in the category of weak bimodules over $\Assoc$, bimodules over $\Assoc$, and operads respectively. These objects are naturally filtered and one has that $\hoTriangle_N|_N$, $\hoSquare_N|_N$, $\hoPentagon_N|_N$ are cofibrant models of $\Assoc|_N$ in the corresponding category of $N$-truncated objects. One of the essential differences with Part~\ref{part1} is that $\hoTriangle$, $\hoSquare$, and $\hoPentagon$ are not any more polytopes in each degree, but infinite-dimensional $CW$-complexes. However we will still have that the filtration terms $\hoTriangle_N$, $\hoSquare_N$, and $\hoPentagon_N$ are cellular cofibrant and are freely generated by finitely many cells. To pass from $\hoTriangle_{N-1}$, $\hoSquare_{N-1}$,  $\hoPentagon_{N-1}$ to $\hoTriangle_N$, $\hoSquare_N$, $\hoPentagon_N$ respectively, one has to attach not one but $2^N$ cells. To be precise one has to attach first one cell in degree $N$ of dimension $N$, $N-1$, $N-2$ respectively to the three cases, then one attaches $N\choose 1$ cells in degree $N-1$ and of dimension $N+1$, $N$, $N-1$ respectively, $\ldots$, on the $i-th$ step one attaches $N\choose i$ cells in degree $N-i$ of dimension $N+i$, $N+i-1$, $N+i-2$ respectively. One can see that similarly with the case considered in the first part of the paper the dimension of all generating cells decreases by one each time we pass from $\hoTriangle$ to $\hoSquare$, or from $\hoSquare$ to $\hoPentagon$. This is a crucial fact used to show that the delooping maps $\tilde\xi_n$ and $\tilde\zeta_n$ considered below  are homotopy equivalences.

In Section~\ref{s:first_deloop_tilde} we construct a weak bimodule $Wb\hoSquare$ that has the same properties as $\hoTriangle$: it is a cofibrant model of $\Assoc$, it is naturally filtered so that $Wb\hoSquare_N|_N$ is a cofibrant model of $\Assoc|_N$, and one also needs $2^N$ new cells to get $Wb\hoSquare_N$ from $Wb\hoSquare_{N-1}$ (that are attached in the same order and in the same degrees and are of the same dimensions as in the case of $\hoTriangle$). The essential difference with Part~\ref{part1} is that $Wb\hoSquare$ is not homeomorphic to $\hoTriangle$, but  is only homotopy equivalent to it (as a filtered weak bimodule). By the construction we have a natural map
$$
\tilde\xi_N\colon\Omega\underset{\Assoc}{\Bimod}{}_N(\square_N|_N,\calO|_N)\to \underset{\Assoc}{\WBimod}{}_N(Wb\square_N|_N,\calO_N),
$$
that we show to be a homotopy equivalence in case $\calO(0)\simeq*$. The latter map  describes the first delooping of $\hoTot_N\calO(\bullet)$ as the space of derived morphisms of bimodules:
$$
\Omega\underset{\Assoc}{\widetilde{\Bimod}}{}_N(\Assoc|_N,\calO|_N)\simeq \hoTot_N\calO(\bullet).
$$

In Section~\ref{s:second_deloop_tilde} we construct an $\Assoc$ bimodule $B\hoPentagon$ with the same properties as $\hoSquare$: it is a cofibrant model of $\Assoc$, it is naturally filtered with each $B\hoPentagon_N$ cofibrant and $B\hoPentagon_N|_N$ being a cofibrant model of $\Assoc|_N$  as an $N$-truncated $\Assoc$ bimodule. Similarly $B\hoPentagon_N$ is obtained from $B\hoPentagon_{N-1}$ by a free attachment of $2^N$ cells. Again we will have that $B\hoPentagon$ is not homeomorphic to $\hoSquare$, but only homotopy equivalent to it (as a filtered bimodule). By construction one has a natural map
$$
\tilde\zeta_N\colon\Omega\Operad_N(\hoPentagon_N|_N,\calO_N)\to\underset{\Assoc}{\Bimod}{}_N(B\hoPentagon_N|_N,\calO|_N),
$$
that we show to be a homotopy equivalence in case $\calO(1)\simeq *$. Combining the two results we obtain an equivalence
$$
\Omega^2\widetilde{\Operad}_N(\Assoc|_N,\calO|_N)\simeq\hoTot_N\calO(\bullet),\quad N\geq 0,
$$
assuming $\calO(0)\simeq\calO(1)\simeq *$. Taking the limit $N\to\infty$ one reproduces the Dwyer-Hess' theorem~\cite{DwHe}:
$$
\Omega^2\widetilde{\Operad}(\Assoc,\calO)\simeq\hoTot\calO(\bullet).
$$

\section{Cofibrant model $\widetilde{\triangle}$ of $\Assoc$ as a weak bimodule over itself\\ and associated tower}\label{s:deg_wb}

\subsection{Cofibrant model $\widetilde{\triangle}$}\label{ss:deg_wb}
To recall $\{\triangle(n),\, n\geq 0\}$, is a weak bimodule over $\Assoc_{>0}$ given by a sequence of simplices, where each simplex can be viewed as a configuration space of $n$ points on the interval $[0,1]$:
$$
0\leq t_1\leq t_2\leq \ldots\leq t_n\leq 1.
$$
The right action of $\Assoc_{>0}$ consists in doubling points, and the weak left action consists in adding new points at the ends of $[0,1]$. This action can be naturally extended to an action of $\Assoc$. The zero component $\Assoc(0)=*$ can act only to the right. We define this action
$$
\circ_i\colon \triangle(n)\times \Assoc(0)\to \triangle(n-1)
\eqno(\numb)\label{eq:right_assoc_0_action}
$$
as forgetting the $i$-th point in the configuration.

Obviously $\triangle$ is not cofibrant as a weak bimodule over $\Assoc$, since the degeneracies~\eqref{eq:right_assoc_0_action} do not act freely. Below we define its cofibrant replacement $\widetilde{\triangle}$.

A point in $\widetilde{\triangle}(n)$, $n\geq 0$, is a configuration space of $n$ points on $[0,1]$ (labeled by $1,2,\ldots,n$ in increasing order) together with some number of hairs growing from the interior $(0,1)$ of the segment and that have length $\leq 1$. Such configurations will be also called {\it hairy configurations}.

\vspace{.3cm}

\begin{figure}[h]
\includegraphics[width=5cm]{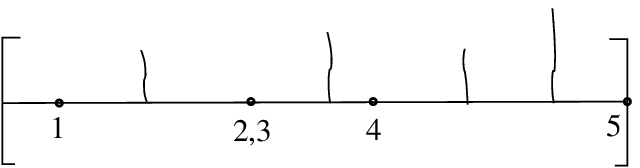}
\caption{An example of a hairy configuration with 7 geometrically distinct points. This configuration is a point in $\widetilde{\triangle}(5)$.}\label{fig6}
\end{figure}

The points labeled by $1,2,\ldots,n$ are allowed to collide with each other and with the endpoints of the interval $[0,1]$, but not with the hairs. Moreover hairs are not allowed to grow from the same point. If a hair shortens to length~0, it disappears. If two hairs collide, only the one of the longer length survives. If a hair collides with one of the labeled vertices or with one of the endpoints of $[0,1]$, then it also disappears.

The left action of $\Assoc$ on $\widetilde{\triangle}=\{\widetilde{\triangle}(n),\, n\geq 0\}$ is defined similarly to that on $\triangle$ --- by adding (labeled) points on the boundary of $[0,1]$. The right action of $\Assoc_{>0}$ is defined by multiplying the labeled points. The right action of $\Assoc(0)$ is defined by replacing the corresponding labeled point by a hair of length~1, in case the labeled point was inside the interval and did not coincide with other labeled points, and by forgetting this point otherwise. It is easy to see that with this definition $\widetilde{\triangle}$ becomes a weak bimodule over $\Assoc$.

We will consider a filtration in $\widetilde{\triangle}$:
$$
\widetilde{\triangle}_0\subset\widetilde{\triangle}_1\subset\widetilde{\triangle}_2\subset\ldots \subset\widetilde{\triangle},
\eqno(\numb)\label{eq:filtr_tilde_triangle}
$$
where $\widetilde{\triangle}_i(n)$ for each $n\geq 0$ consists of configurations having $\leq i$ geometrically distinct points strictly inside the interval $(0,1)$. For example, for a hairy configuration from Figure~\ref{fig6} the number of such points is~7, since the points~2 and~3 are counted for one, and the point~5 does not contribute being outside~$(0,1)$.
Obviously, the action of $\Assoc$ preserves this filtration.

\begin{figure}[h]
\psfrag{hD0}[0][0][1][0]{$\widetilde{\triangle}_1(0)$}
\psfrag{hD1}[0][0][1][0]{$\widetilde{\triangle}_1(1)$}
\psfrag{hD2}[0][0][1][0]{$\widetilde{\triangle}_1(2)$}
\includegraphics[width=8cm]{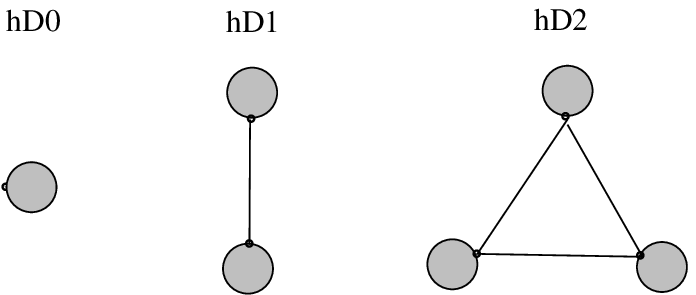}
\caption{$\widetilde{\triangle}_1$ in small degrees.}\label{fig7}
\end{figure}

\begin{proposition}\label{p:cof_tilde_tr}
(1) $\widetilde\triangle$ is a cofibrant model of $\Assoc$ as a weak bimodule over $\Assoc$.

(2) $\widetilde{\triangle}_n|_n$ is a cofibrant model of $\Assoc|_n$ as a truncated weak bimodule over $\Assoc$.
\end{proposition}

\begin{proof}
We start by showing that $\widetilde{\triangle}(i)$ (and $\widetilde{\triangle}_n(i)$) are contractible for any $i$ (any $i\leq n$, respectively). By contracting all the hairs the spaces $\widetilde{\triangle}(i)$ (respectively $\widetilde{\triangle}_n(i)$) are retracted to $\triangle(i)$ (respectively $\triangle_n(i)$). Then we notice that $\triangle(i)=\triangle_n(i)$ are contractible being simplices.

Secondary we have to show that $\widetilde{\triangle}$ (and $\widetilde{\triangle}_n|_n$) are cofibrant, or in other words are obtained by a sequence of free attachments of cells, see Section~\ref{s:col_operads}. We give an argument for $\widetilde{\triangle}$. The case of $\widetilde{\triangle}_N|_N$ follows from Lemma~\ref{l:truncation}.

One starts by putting a 0-cell in degree zero (empty configuration) which generates $\widetilde{\triangle}_0$ whose points are configurations of labeled points all seating at the endpoints of $[0,1]$. Obviously, $\widetilde{\triangle}_0=\triangle_0$ is free. Then one attaches a 1-cell in degree~1 --- this cell corresponds to configurations with only one labeled point and without hairs. This cell generates cells with exactly one geometrically distinct point such that if this point is a hair, then its length is~1. Then we attach a 2-cell in degree~0 corresponding to configurations without label points and with one hair of length $<1$. Attaching this cell finishes creating $\widetilde{\triangle}_1$. In general to pass from $\widetilde{\triangle}_{n-1}$ to $\widetilde{\triangle}_n$ one has to attach 1 cell of dimension $n$ in degree $n$, $n$ cells of dimension $n+1$ in degree $n-1$, ${n\choose 2}$ cells of dimension $n+2$ in degree $n-2$,$\ldots$, $n\choose i$ cells of dimension $n+i$ in degree $n-i$ (these cells contain configurations with $(n-i)$ labeled points and $i$ hairs; we will denote by $\widetilde{\triangle}_{n-1,i}$ the weak bimodule obtained at this point),$\ldots$, ${n\choose n}=1$ cell of dimension $2n$ in degree~0. Thus we made it explicit that the construction is cofibrant. The truncated case follows from~\ref{l:truncation}.
\end{proof}

We also notice that the space $\widetilde{\triangle}(n)$ can be obtained as a quotient
$$
\widetilde{\triangle}(n)=\left.\left(\coprod_{m\geq 0}\,\,\,\coprod_{\alpha\colon\underline{m}\hookrightarrow\underline{n+m}}\triangle(n+m)\times[0,1]^{\alpha(\underline{m})}\right)\right/\sim,
\eqno(\numb)\label{eq:tTriangle_as_union}
$$
where the second disjoint union is taken over the order preserving inclusions. The quotient relation can be easily obtained from the definition given in the beginning of this section.\footnote{We write down this quotient relation in Sections~\ref{ss:tilde_square},~\ref{ss:cof_tilde_wb} in a more general context.} Notice again that the interior of each $\triangle(n+m)\times[0,1]^{\alpha(\underline{m})}$ in~\eqref{eq:tTriangle_as_union} is a generating cell of $\widetilde{\triangle}$ (such cell is attached when passing from $\widetilde{\triangle}_{n+m-1,m-1}$ to $\widetilde{\triangle}_{n+m-1,m}$).

\subsection{Tower $T_\bullet^{\widetilde{\triangle}}(\calO)$ associated to $\widetilde{\triangle}$}\label{ss:tower_tilde_wb}
Given a weak bimodule $\calO$ over $\Assoc$, one can consider a tower of fibrations
$$
T_0^{\widetilde{\triangle}}(\calO)\leftarrow T_1^{\widetilde{\triangle}}(\calO)\leftarrow T_2^{\widetilde{\triangle}}(\calO)\leftarrow T_3^{\widetilde{\triangle}}(\calO)\leftarrow\ldots,
\eqno(\numb)\label{eq:tower_tilde_triangle}
$$
whose stages are spaces of morphisms
$$
T_n^{\widetilde{\triangle}}(\calO)=\underset{\Assoc}{\WBimod}(\widetilde{\triangle}_n,\calO)=\underset{\Assoc}{\WBimod_n}(\widetilde{\triangle}_n|_n,\calO|_n)\simeq
\underset{\Assoc}{\widetilde{\WBimod}}{}_n(\Assoc|_n,\calO|_n).
$$
The limit of this tower is
$$
T_{\infty}^{\widetilde{\triangle}}(\calO)=\underset{\Assoc}{\WBimod}(\widetilde{\triangle},\calO)=\underset{\Assoc}{\widetilde{\WBimod}}(\Assoc,\calO).
$$
Equivalently the above tower is the tower of partial homotopy totalizations
$$
\hoTot_0\calO(\bullet)\leftarrow\hoTot_1\calO(\bullet)\leftarrow\hoTot_2\calO(\bullet)\leftarrow\ldots
$$
converging to the homotopy totalization $\hoTot\calO(\bullet)$.

The advantage of considering the tower~\eqref{eq:tower_tilde_triangle} compared to the tower~\eqref{eq:tower_triangle} is that~\eqref{eq:tower_tilde_triangle} has nicer convergency properties.

The map $T_n^{\widetilde{\triangle}}(\calO)\to T_{n-1}^{\widetilde{\triangle}}(\calO)$ can be decomposed as
$$
T_{n-1}^{\widetilde{\triangle}}(\calO)\leftarrow T_{n-1,0}^{\widetilde{\triangle}}(\calO)\leftarrow T_{n-1,1}^{\widetilde{\triangle}}(\calO)\leftarrow\ldots
\leftarrow T_{n-1,n}^{\widetilde{\triangle}}(\calO)=T_n^{\widetilde{\triangle}}(\calO),
$$
where we denote by $T_{n-1,i}^{\widetilde{\triangle}}(\calO)$ the space of morphisms $\underset{\Assoc}{\WBimod}(\widetilde{\triangle}_{n-1,i},\calO)$. A variation of these intermediate stages will be used in the proof of Theorem~\ref{t:deloop1_stages_tilde}.

To describe the fiber of the map $T_n^{\widetilde{\triangle}}(\calO)\to T_{n-1}^{\widetilde{\triangle}}(\calO)$ we will need the notion of  the total homotopy fiber of a cubical diagram. The description of the fibers is given for completeness of exposition and is not used in the proof of the main results.

\subsection{Cubical diagrams and fiber of $T_n^{\widetilde{\triangle}}(\calO)\to T_{n-1}^{\widetilde{\triangle}}(\calO)$}\label{ss:cubic}
Let $\calP_n$ be the poset of subsets of $\underline{n}=\{1,2,\ldots,n\}$. An $n$-{\it cubical diagram} (or an $n$-{\it cube}) is a functor from $\calP_n$.

Denote by $[0,1]^\bullet$ the $n$-cubical diagram that assigns $[0,1]^S$ to each set $S\subset\underline{n}$. For each inclusion $S\subset T$ the corresponding map
$$
[0,1]^S\to [0,1]^T
$$
sends a function $f\colon S\to [0,1]$ to the function $T\to [0,1]$ defined by
$$
t\mapsto
\begin{cases}
f(t),&t\in S;\\
1,& t\in T\setminus S.
\end{cases}
$$

Denote by $\left([0,1]^S\right)_0$ the subset of $[0,1]^S$ consisting of functions $f\colon S\to [0,1]$, such that $0\in f(S)$. We will also denote by $\left([0,1]^\bullet\right)_0$ the corresponding $n$-cube and by $\left([0,1]^\bullet,\left([0,1]^\bullet\right)_0\right)$ the $n$-cube in the category of pairs of spaces.

\begin{definition}\label{d:tot_fib}
Let $\sfX$ be an $n$-cube of based spaces. We define $tfiber(\sfX)$ the {\it total homotopy fiber} (or simply {\it total
fiber}) of $\sfX$ as the space
$$
tfiber(\sfX)={\Top^{\calP_n}}_*\Bigl(\left([0,1]^\bullet,\left([0,1]^\bullet\right)_0\right),\sfX\bigr)
\eqno(\numb)\label{eq:tot_fib}
$$
of natural transformations $[0,1]^\bullet\to\sfX$, that send each space $\left([0,1]^S\right)_0$ to the base point of $\sfX(S)$.
\end{definition}

The reader can consult \cite[Section~1]{Good-Calc2} for a better introduction to this notion of total fiber.

Notice that in case $\sfX$ is an $n$-cube of non-based spaces, one can still define its total fiber by  choosing a base point in $\sfX(\emptyset)$ (which as an image defines  base points in each $\sfX(S)$). The total fiber thus defined depends only on the connected component of $\sfX(\emptyset)$ where the base point is chosen.

\begin{notation}\label{n:nat_cubes}
More generally given an $n$-cube $(\sfX_1,\sfX_0)$ in pairs of spaces and an $n$-cube $\sfY$ in spaces, and a natural transformation $\sfX_0\stackrel{\beta}{\longrightarrow}\sfY$, we consider the space
$$
{\Top^{\calP_n}}_\alpha\bigl((\sfX_1,\sfX_0),\sfY\bigr)
$$
of natural transformations $\beta\colon\sfX_1\to\sfY$ such that $\beta|_{\sfX_0}=\alpha$.
\end{notation}

We will describe the fiber of the map $T_n^{\widetilde{\triangle}}(\calO)\to T_{n-1}^{\widetilde{\triangle}}(\calO)$ using the above construction.

As before we assume that $\calO$ is a weak bimodule over $\Assoc$. For any finite set $R$ we define $\calO(R)$ as $\calO(|R|)$. We denote by $\calO(\underline{n}\setminus \bullet)$ the $n$-cube that assigns $\calO(\underline{n}\setminus S)$ to any set $S\subset\underline{n}$. On morphisms this cubical diagram is defined using degeneracies (right action of $\Assoc(0)$) in an obvious way.

Let $\triangle(n)\times [0,1]^\bullet$ denote the $n$-cube obtained from $[0,1]^\bullet$ by taking the product of each object with the space $\triangle(n)$. We will also consider $\partial_0\left(\triangle(n)\times [0,1]^\bullet\right)$ -- a cubical subobject of  $\triangle(n)\times [0,1]^\bullet$ defined by
$$
\partial_0\left(\triangle(n)\times [0,1]^S\right):=
\Bigl(\partial\bigl(\triangle(n)\bigr)\times [0,1]^S\Bigr)\cup \Bigl(\triangle(n)\times \left([0,1]^S\right)_0\Bigr).
$$

Let $\beta_S$ denote the map
$$
\triangle(n)\times [0,1]^S\to \widetilde{\triangle}(\underline{n}\setminus S),
$$
which is the composition of the inclusion
$$
i_S\colon \triangle(n)\times [0,1]^S\hookrightarrow \coprod_{m\geq 0}\,\,\,\coprod_{\alpha\colon\underline{m}\hookrightarrow\underline{n-|S|+m}} \triangle(n-|S|+m)\times [0,1]^{\alpha(\underline{m})}
$$
(that sends $\triangle(n)\times [0,1]^S$ identically to the component labeled by the inclusion $\underline{|S|}\simeq S\subset\underline{n}$), and the quotient map
$$
\coprod_{m\geq 0}\,\,\,\coprod_{\alpha\colon\underline{m}\hookrightarrow\underline{n-|S|+m}} \triangle(n-|S|+m)\times [0,1]^{\alpha(\underline{m})}\to\left(\coprod_{m\geq 0}\,\,\,\coprod_{\alpha\colon\underline{m}\hookrightarrow\underline{n-|S|+m}} \triangle(n-|S|+m)\times [0,1]^{\alpha(\underline{m})}\right)\Bigl/\sim,
$$
where the right-hand side is $\widetilde{\triangle}(n-|S|)$, see~\eqref{eq:tTriangle_as_union}. The maps $\beta_S$, $S\subset\underline{n}$ define a natural transformation $\beta\colon \triangle(n)\times [0,1]^\bullet\to \widetilde{\triangle}(\underline{n}\setminus \bullet)$. Notice that preimage of the subcube $\widetilde{\triangle}_{n-1}(\underline{n}\setminus\bullet)\subset \widetilde{\triangle}(\underline{n}\setminus\bullet)$ under $\beta$ is exactly $\partial_0\left(\triangle(n)\times [0,1]^\bullet\right)$.

Let us fix a point $\gamma\in T_{n-1}^{\widetilde{\triangle}}(\calO)$, which is a morphism of weak $\Assoc$ bimodules $\gamma\colon\widetilde{\triangle}_{n-1}\to\calO$. By abuse of notation denote by $\gamma$ the induced natural transformation of cubical diagrams
$$
\gamma\colon \widetilde{\triangle}_{n-1}(\underline{n}\setminus\bullet)\to \calO (\underline{n}\setminus\bullet).
$$

\begin{proposition}\label{p:fiber_tilde_trianle}
The preimage of $\gamma\in T_{n-1}^{\widetilde{\triangle}}(\calO)$, $n\geq 0$,\footnote{We are assuming $T_{-1}^{\widetilde{\triangle}}(\calO)=*$.} under the fibration $T_n^{\widetilde{\triangle}}(\calO)\to T_{n-1}^{\widetilde{\triangle}}(\calO)$ is homeomorphic to
$$
{\Top^{\calP_n}}_{\gamma\circ\beta}\Bigl(\bigl(\triangle(n)\times[0,1]^\bullet,\partial_0\left(\triangle(n)\times [0,1]^\bullet\right)\bigr),\calO(\underline{n}\setminus\bullet)\Bigr).
$$
\end{proposition}

\begin{proof} Direct inspection.
\end{proof}

The following statement can be easily proved from Proposition~\ref{p:fiber_tilde_trianle}.

\begin{proposition}\label{p:fiber_tilde_trianle2}
The preimage of the fibration $T_n^{\widetilde{\triangle}}(\calO)\to T_{n-1}^{\widetilde{\triangle}}(\calO)$, $n\geq 0$, is either empty or homotopy equivalent to $\Omega^n tfiber\left(\calO(\underline{n}\setminus\bullet)\right)$, where (in case $n\geq 1$) the base point of $\calO(n)$ can be taken to be $a_{n+1}\circ_{n+1}p$ with $p\in\calO(0)$ being the image of the map $\widetilde{\triangle}_{0}(0)\to\calO(0)$.\footnote{To recall $\widetilde{\triangle}_{0}(0)$ is a point.} In case the space $tfiber\left(\calO(\underline{n}\setminus\bullet)\right)$ is $(n-1)$-connected, the fiber is always homotopy equivalent to $\Omega^n tfiber\left(\calO(\underline{n}\setminus\bullet)\right)$.
\end{proposition}

\section{Cofibrant model $\widetilde{\square}$ and associated tower $T_\bullet^{\hoSquare}(\calO)$}\label{s:cof_tilde_sq}

\subsection{$\square$ as a bimodule over $\Assoc$}\label{ss:bimod0_sq}
Recall Section~\ref{ss:cof_bim2} where we defined an $\Assoc-\Assoc_{>0}$ bimodule $\square$. In this section we will define a structure of an $\Assoc$ bimodule on it.


Each space $\square(n)$, $n\geq 0$, can be viewed as a configuration space of $n$ points on $\R^1$ labeled by $1,2,\ldots,n$ with the possible distance between consecutive points to be $\leq 1$ (including 0 and 1). This configuration space is quotiented out by translations so that a configuration is determined by a sequence of distances $(d_1,d_2,\ldots,d_{n-1})$ between the points. These distances $d_i$, $i=1\ldots n-1$ are coordinates of the cube $\square(n)$.
The left action of $a_2\in\Assoc(2)$ concatenates configurations putting the distance~1 between them:
\begin{gather*}
\Assoc(2)\times\square(k)\times\square(m)\to\square(k+m);\\
a_2\times (d_1\ldots d_{k-1})\times (d_1'\ldots d_{m-1}')\mapsto (d_1\ldots d_{k-1},1,d_1'\ldots d_{m-1}').
\end{gather*}
The right action of $a_2$ doubles points in the configuration:
\begin{gather*}
\circ_i\colon\square(n)\times\Assoc(2)\to \square(n+1)\\
(d_1\ldots d_{n-1})\times a_2\mapsto (d_1\ldots d_{i-1},0,d_i\ldots d_{n-1}).
\end{gather*}

Finally we define the right action of $\Assoc(0)$ by forgetting the corresponding point in configurations. In case the point was internal, the new distance between the point to the left and the one to the right will be the maximum of two:
$$
(d_1\ldots d_{n-1})\circ_i a_0=
\begin{cases}
(d_2\ldots d_{n-1}),& i=1;\\
(d_1\ldots d_{i-2},\mathrm{max}(d_{i-1},d_i),d_{i+1}\ldots d_{n-1}),& i=2\ldots n-1;\\
(d_1\ldots d_{n-2}),& i=n.
\end{cases}
\eqno(\numb)\label{eq:square_degen}
$$
One can easily see that this definition is consistent and turns $\square$ into a bimodule over $\Assoc$.

\subsection{Cofibrant replacement $\widetilde{\square}$}\label{ss:tilde_square}
In this section we define a cofibrant model of $\Assoc$ in the category of $\Assoc$ bimodules. The $n$-th component $\widetilde{\square}$ of this bimodule is defined as a quotient
$$
\widetilde{\square}(n)=\left.\left(\coprod_{m\geq 0}\,\,\,\coprod_{\alpha\colon\underline{m}\hookrightarrow\underline{n+m}}\square(n+m)\times[0,1]^{\alpha(\underline{m})}\right)\right/\sim.
\eqno(\numb)\label{eq:tSquare_as_union}
$$
We will specify the equivalence relation a little bit later. A point in $\widetilde{\square}(n)$ is a configuration of $n$ points on $\R^1$  labeled by $1\ldots n$ in increasing order and a number of hairs growing from $\R^1$. Again the configurations are considered modulo translations of $\R^1$. One has the following restrictions on configurations:
\begin{itemize}
\item the distance between neighbor points/hairs is $\leq 1$;
\item the length of hairs is $\leq 1$;
\item a hair can not have distance 0 with another hair or with a labeled point;
\end{itemize}
Two labeled points are allowed to have distance 0. If a hair shortens to length 0, it disappears, the distance between its neighbor to the left and the one to the right becomes the maximum of two distances as in formula~\ref{eq:square_degen}. If a hair collides with a labeled point it disappears. If two hairs collide together, only the hair of longer length survives.

The construction of~$\hoSquare$ follows a general framework. Given any bimodule ${\calY}$ over $\Assoc$, we will define another $\Assoc$ bimodule by
$$
\widetilde{\calY}(n) :=\left.\left(\coprod_{m\geq 0}\,\,\,\coprod_{\alpha\colon\underline{m}\hookrightarrow\underline{m+n}}{\calY}(m+n)\times[0,1]^{\alpha(\underline{m})}\right)\right/\sim.
$$
Any element $z\in\widetilde{\calY}(n)$ can be written as a tuple $(y;\tau_{\alpha(1)}\ldots \tau_{\alpha(m)})_\alpha$ labeled by an inclusion $\alpha\colon \underline{m}\hookrightarrow\underline{n+m}$, where $y\in{\calY}(n+m)$, and $0\leq \tau_{\alpha(i)}\leq 1$. Graphically such $z$ can be represented as the element $y$ corked in its entries $\alpha(1),\ldots,\alpha(m)$ by the segments of length $\tau_{\alpha(1)},\ldots,\tau_{\alpha(m)}$ respectively.

One of the equivalence relations is as follows:
$$
(y;\tau_{\alpha(1)},\ldots,\tau_{\alpha(i)}=0,\ldots,\tau_{\alpha(m)})_\alpha\sim (y\circ_{\alpha(i)}a_0,\tau_{\alpha(1)},\ldots,\widehat{\tau_{\alpha(i)}},\ldots,\tau_{\alpha(m)})_{\alpha_{\widehat{i}}},
\eqno(\numb)\label{eq:Yrel1}
$$
where $\alpha_{\widehat{i}}$ denotes the composition of order preserving maps
$$
\alpha_{\widehat{i}}\colon\underline{m-1}\simeq \underline{m}\setminus\{i\}\stackrel{\alpha}{\hookrightarrow} \underline{m+n}\setminus\{\alpha(i)\}\simeq \underline{m+n-1}.
\eqno(\numb)\label{eq:alpha_i}
$$
In the other relations~\eqref{eq:Yrel2},~\eqref{eq:Yrel3} we will be assuming  that $\alpha(i)<j<\alpha(i+1)$.

The other relations are:
\begin{multline}
\left(y;\tau_{\alpha(1)},\ldots,\tau_{\alpha(m)}\right)_\alpha=\left(y\circ_{j}a_2;\tau_{\alpha(1)},\ldots,\tau_{\alpha(i)},\tau,\tau_{\alpha(i+1)},\ldots,
\tau_{\alpha(m)}\right)_{\alpha_{j^0}}=\\
\left(y\circ_ja_2;\tau_{\alpha(1)},\ldots,\tau_{\alpha(i)},\tau,\tau_{\alpha(i+1)},\ldots,
\tau_{\alpha(m)}\right)_{\alpha_{j^+}},\label{eq:Yrel2}
\end{multline}
where $\alpha_{j^0}\colon\underline{m+1}\hookrightarrow\underline{m+n+1}$ (respectively $\alpha_{j^+}\colon\underline{m+1}\hookrightarrow\underline{m+n+1}$) is an order preserving map whose image is $\alpha(\underline{i})\cup\{j\}\cup\left(\alpha(\underline{m}\setminus\underline{i})+1\right)$ (respectively
$\alpha(\underline{i})\cup\{j+1\}\cup\left(\alpha(\underline{m}\setminus\underline{i})+1\right)$.
%
%
%
%
%
%

Finally the last relation is
\begin{multline}
\left(y\circ_ja_2;\tau_{\alpha(1)},\ldots,\tau_{\alpha(i)},\tau,\tau'\tau_{\alpha(i+1)},\ldots,t_{\alpha(m)}\right)_{\alpha_{j^{0+}}}=\\
\left(y;\tau_{\alpha(1)},\ldots,\tau_{\alpha(i)},max(\tau,\tau'),\tau_{\alpha(i+1)},\ldots,\tau_{\alpha(m)}\right)_{\alpha_j},
\label{eq:Yrel3}
\end{multline}
where $\alpha_{j^{0+}}\colon \underline{m+2}\hookrightarrow \underline{m+n+1}$ has image $\alpha_{j^0}(\underline{m+1})\cup\{j+1\}$, and $\alpha_{j}\colon\underline{m+1}\hookrightarrow\underline{m+n+1}$ has image $\alpha(\underline{m})\cup\{j\}$.

The $\Assoc$ bimodule structure is defined in the natural way.

Since this construction is functorial the filtration
$$
\square_0\subset\square_1\subset\square_2\subset\ldots\subset\square
$$
induces filtration of $\Assoc$ bimodules
$$
\widetilde{\square}_0\subset\widetilde{\square}_1\subset\widetilde{\square}_2\subset\ldots\subset\widetilde{\square}.
\eqno(\numb)\label{eq:filtr_tilde_square}
$$

\begin{proposition}\label{p:cof_bimod_tilde}
\textup{(1)} One has that $\widetilde{\square}$ is a cofibrant model of $\Assoc$ in the category of $\Assoc$ bimodules.

\textup{(2)} For any $n\geq 0$, $\widetilde{\square}_n|_n$ is a cofibrant model of $\Assoc|_n$ in the category of $n$-truncated $\Assoc$ bimodules.
\end{proposition}

\begin{proof}
The proof is similar to that of Proposition~\ref{p:cof_tilde_tr}. We will just describe  the sequence in which we attach cells to get $\widetilde{\square}$. The free $\Assoc$ bimodule generated by empty sets in each degree is the bimodule that has only one point $\unit$ in degree zero, which is the result of left action of $\Assoc(0)$. This element corresponds to an empty configuration. The obtained bimodule is the degree zero term $\widetilde{\square}_0$ of the filtration~\eqref{eq:filtr_tilde_square}. Then we attach a 0-cell in degree 1, that corresponds to the configuration of only one labeled point. Acted on by $\Assoc$, this cell produces configurations of hairs and labeled points, in which all the hairs have length~1 and the distance between labeled points/hairs is either~0 or~1. Then we attach a 1-cell in degree zero, that corresponds to the configurations consisting of a single hair of positive length $<1$. The obtained bimodule is $\widetilde{\square}_1$, which consists of configurations of hairs (of any length $\leq 1$) and labeled points with the only condition that the distance between labeled points/hairs is either~0 or~1. To get $\widetilde{\square}_2$ we attach a 1-cell in degree~2 (corresponding to configurations of 2 labeled points without hairs), then two 2-cells in degree~1 (corresponding to configurations consisting of 1 hair and one labeled point), then one 3-cell in degree~0 (corresponding to configurations of 2 hairs without labeled points). In general to pass from $\widetilde{\square}_{n-1}$ to $\widetilde{\square}_n$ we need to attach one $(n-1)$-cell in degree $n$,
then $n\choose 1$ $n$-cells in degree $n-1$,
then $n\choose 2$ $(n+1)$-cells in degree $n-2$,
$\ldots$, then $n\choose i$ $(n+i-1)$-cells in degree $n-i$ (these cells consist of configurations with $i$ hairs and $(n-i)$ labeled points; the bimodule obtained at this point will be denoted $\widetilde{\square}_{n-1,i}$), $\ldots$, at the end we attach one $(2n-1)$-cell in degree~0 (the obtained bimodule is $\widetilde{\square}_{n-1,n}=\widetilde{\square}_{n}$).
\end{proof}

We stress the fact that the interior of each $\square(n+m)\times[0,1]^{\alpha(\underline{m})}$, $n+m>0$, in~\eqref{eq:tSquare_as_union} is a generating cell of the $\Assoc$ bimodule structure.

\subsection{Tower $T_\bullet^{\widetilde{\square}}(\calO)$}\label{ss:tower_tilde_bimod}
Let $\calO$ be a bimodule over $\Assoc$. Define the space $T_n^{\widetilde{\square}}(\calO)$ by
$$
T_n^{\widetilde{\square}}(\calO):=
\underset{\Assoc}{\Bimod}(\widetilde{\square}_n,\calO)=\underset{\Assoc}{\Bimod}{}_n(\widetilde{\square}_n|_n,\calO|_n)\simeq
\widetilde{\underset{\Assoc}{\Bimod}}{}_n
(\Assoc|_n,\calO|_n).
\eqno(\numb)\label{eq:tilde_bimod_stages}
$$
The second equation is a consequence of Lemma~\ref{l:truncation}. These spaces fit together into a tower of fibrations
$$
T_0^{\widetilde{\square}}(\calO)\leftarrow T_1^{\widetilde{\square}}(\calO)\leftarrow T_2^{\widetilde{\square}}(\calO)\leftarrow\ldots,
\eqno(\numb)\label{eq:tower_tilde_Square}
$$
whose (homotopy) limit is $T_\infty^{\widetilde{\square}}(\calO)=\underset{\Assoc}{\Bimod}(\widetilde{\square},\calO)\simeq\widetilde{\underset{\Assoc}{\Bimod}}(\Assoc,\calO)$.

The map $T_n^{\widetilde{\square}}(\calO)\to T_{n-1}^{\widetilde{\square}}(\calO)$ can be decomposed as
$$
T_{n-1}^{\widetilde{\square}}(\calO)\leftarrow T_{n-1,0}^{\widetilde{\square}}(\calO)\leftarrow T_{n-1,1}^{\widetilde{\square}}(\calO)\leftarrow\ldots
\leftarrow T_{n-1,n}^{\widetilde{\square}}(\calO)=T_n^{\widetilde{\square}}(\calO),
$$
where we denote by $T_{n-1,i}^{\widetilde{\square}}(\calO)$ the space of morphisms $\underset{\Assoc}{\Bimod}(\widetilde{\square}_{n-1,i},\calO)$. These intermediate stages will be used in the proof of Theorems~\ref{t:deloop1_stages_tilde} and~\ref{t:deloop2_stages_tilde}.

Now we will describe the fiber of the fibration $T_n^{\widetilde{\square}}(\calO)\to T_{n-1}^{\widetilde{\square}}(\calO)$. The Propositions~\ref{p:fiber_tilde_square},~\ref{p:fiber_tilde_square2} below are given for completeness of exposition and are not used in the proofs of other results.

Since any $\Assoc$ bimodule $\calO$ can be viewed as a right module over $\Assoc$, by a construction similar to Section~\ref{ss:cubic} we will define an $n$-cubical diagram $\calO(\underline{n}\setminus\bullet)$. Let $\gamma\in T_{n-1}^{\widetilde{\square}}(\calO)$. By abuse of notation we will denote by $\gamma$ the induced morphism of cubes:
$$
\gamma\colon\widetilde{\square}_{n-1}(\underline{n}\setminus\bullet)\to\calO(\underline{n}\setminus\bullet).
$$

Let $\square(n)\times [0,1]^\bullet$ denote the $n$-cube obtained from $[0,1]^\bullet$ by taking the product of each object with the space $\square(n)$. We will also consider $\partial_0\left(\square(n)\times [0,1]^\bullet\right)$ -- a cubical subobject of  $\square(n)\times [0,1]^\bullet$ defined by
$$
\partial_0\left(\square(n)\times [0,1]^S\right):=
\Bigl(\partial\bigl(\square(n)\bigr)\times [0,1]^S\Bigr)\cup \Bigl(\square(n)\times \left([0,1]^S\right)_0\Bigr).
$$

Let $\beta_S$ denote the map
$$
\square(n)\times [0,1]^S\to \widetilde{\square}(\underline{n}\setminus S),
$$
which is the composition of the inclusion
$$
i_S\colon \square(n)\times [0,1]^S\hookrightarrow \coprod_{m\geq 0}\,\,\,\coprod_{\alpha\colon\underline{m}\hookrightarrow\underline{n-|S|+m}} \square(n-|S|+m)\times [0,1]^{\alpha(\underline{m})}
$$
(that sends $\square(n)\times [0,1]^S$ identically to the component labeled by the inclusion $\underline{|S|}\simeq S\subset\underline{n}$), and the quotient map
$$
\coprod_{m\geq 0}\,\,\,\coprod_{\alpha\colon\underline{m}\hookrightarrow\underline{n-|S|+m}} \square(n-|S|+m)\times [0,1]^m\to\Bigl(\coprod_{m\geq 0}\,\,\,\coprod_{\alpha\colon\underline{m}\hookrightarrow\underline{n-|S|+m}} \square(n-|S|+m)\times [0,1]^{\alpha(\underline{m})}\Bigr)\Bigl/\sim,
$$
where the target is $\widetilde{\square}(n-|S|)$. The maps $\beta_S$, $S\subset\underline{n}$ define a natural transformation $\beta\colon \square(n)\times [0,1]^\bullet\to \widetilde{\square}(\underline{n}\setminus \bullet)$. Notice that the preimage of the subcube $\widetilde{\square}_{n-1}(\underline{n}\setminus\bullet)\subset \widetilde{\square}(\underline{n}\setminus\bullet)$ under $\beta$ is exactly $\partial_0\left(\square(n)\times [0,1]^\bullet\right)$.

\begin{proposition}\label{p:fiber_tilde_square}
The preimage of $\gamma\in T_{n-1}^{\widetilde{\square}}(\calO)$ under the fibration $T_n^{\widetilde{\square}}(\calO)\to T_{n-1}^{\widetilde{\square}}(\calO)$, $n\geq 1$, is homeomorphic to
$$
{\Top^{\calP_n}}_{\gamma\circ\beta}\Bigl(\bigl(\square(n)\times[0,1]^\bullet,\partial_0\left(\square(n)\times [0,1]^\bullet\right)\bigr),\calO(\underline{n}\setminus\bullet)\Bigr).
$$
\end{proposition}

\begin{proof} Direct inspection.
\end{proof}

The following statement can be easily proved from Proposition~\ref{p:fiber_tilde_square}.

\begin{proposition}\label{p:fiber_tilde_square2}
The preimage of the fibration $T_n^{\widetilde{\square}}(\calO)\to T_{n-1}^{\widetilde{\square}}(\calO)$, $n\geq 2$,\footnote{In case $n=1$ the preimage of $T_1^{\widetilde{\square}}(\calO)\to T_{0}^{\widetilde{\square}}(\calO)$ is the homotopy fiber of the degeneracy map $\calO(1)\to\calO(0)$, where $\calO(0)$ is based in~$\unit$.} is either empty or homotopy equivalent to $\Omega^{n-1} tfiber(\calO(\underline{n}\setminus\bullet))$, where for the base point of $\calO(n)$ one can choose $p\circ_1a_{n}$ with $p\in\calO(1)$ being the image of the map $\widetilde{\square}_{0,0}(1)\to\calO(1)$.\footnote{To recall $\widetilde{\square}_{0,0}(1)$ is a point.}
 In case $tfiber(\calO(\underline{n}\setminus\bullet))$ is $(n-2)$-connected, the preimage is always homotopy equivalent to $\Omega^{n-1} tfiber(\calO(\underline{n}\setminus\bullet))$.
\end{proposition}

\section{First delooping. Cofibrant models $\widetilde{Wb\square}=Wb\widetilde{\square}$}\label{s:first_deloop_tilde}

\subsection{$Wb\square$ as a weak bimodule over $\Assoc$}\label{ss:w_square_ass_bimod}
To recall $Wb\square$ is a weak $\Assoc_{>0}$ bimodule constructed from $\square$, see Section~\ref{s:first_deloop}. Since $\square$ is naturally acted on from the right by $\Assoc(0)$, the weak bimodule $Wb\square$ also inherits this action. If $\bar x=(x_1,\ldots,x_k;t_1,\ldots,t_k)\in Wb\square(n_1+\ldots +n_k)$ (where $x_i\in\square(n_i)$, $1\leq i\leq k$, $0\leq t_1\leq \ldots \leq t_k\leq 1$), then
$$
\bar x\circ_\ell a_0:=(x_1,\ldots,x_{i'}\circ_{j'}a_0,x_{i'+1},\ldots,x_k;t_1,\ldots,t_k),
$$
where $n_1+\ldots+n_{i'-1}<\ell\leq n_1+\ldots+n_{i'}$, $j'=\ell-n_1-\ldots-n_{i'-1}$, and $n_{i'}\geq 2$. If $n_{i'}=1$, one defines
$$
\bar x\circ_\ell a_0:= (x_1,\ldots,x_{i'-1},x_{i'+1},\ldots,x_k;t_1,\ldots,t_{i'-1},t_{i'+1},\ldots,t_k).
\eqno(\numb)\label{eq_tilde_wsq_unit1}
$$
One can see that the above operations together with the weak bi-action of $\Assoc_{>0}$ turn $Wb\square$ into a weak bimodule over $\Assoc$.

\begin{figure}[h]
\psfrag{s1}[0][0][1][0]{$\scriptstyle s_1$}
\psfrag{s2}[0][0][1][0]{$\scriptstyle s_2$}
\psfrag{D2}[0][0][1][0]{$\scriptstyle \triangle(2)$}
\psfrag{WS2}[0][0][1][0]{$\scriptstyle Wb\square(2)$}
\includegraphics[width=8cm]{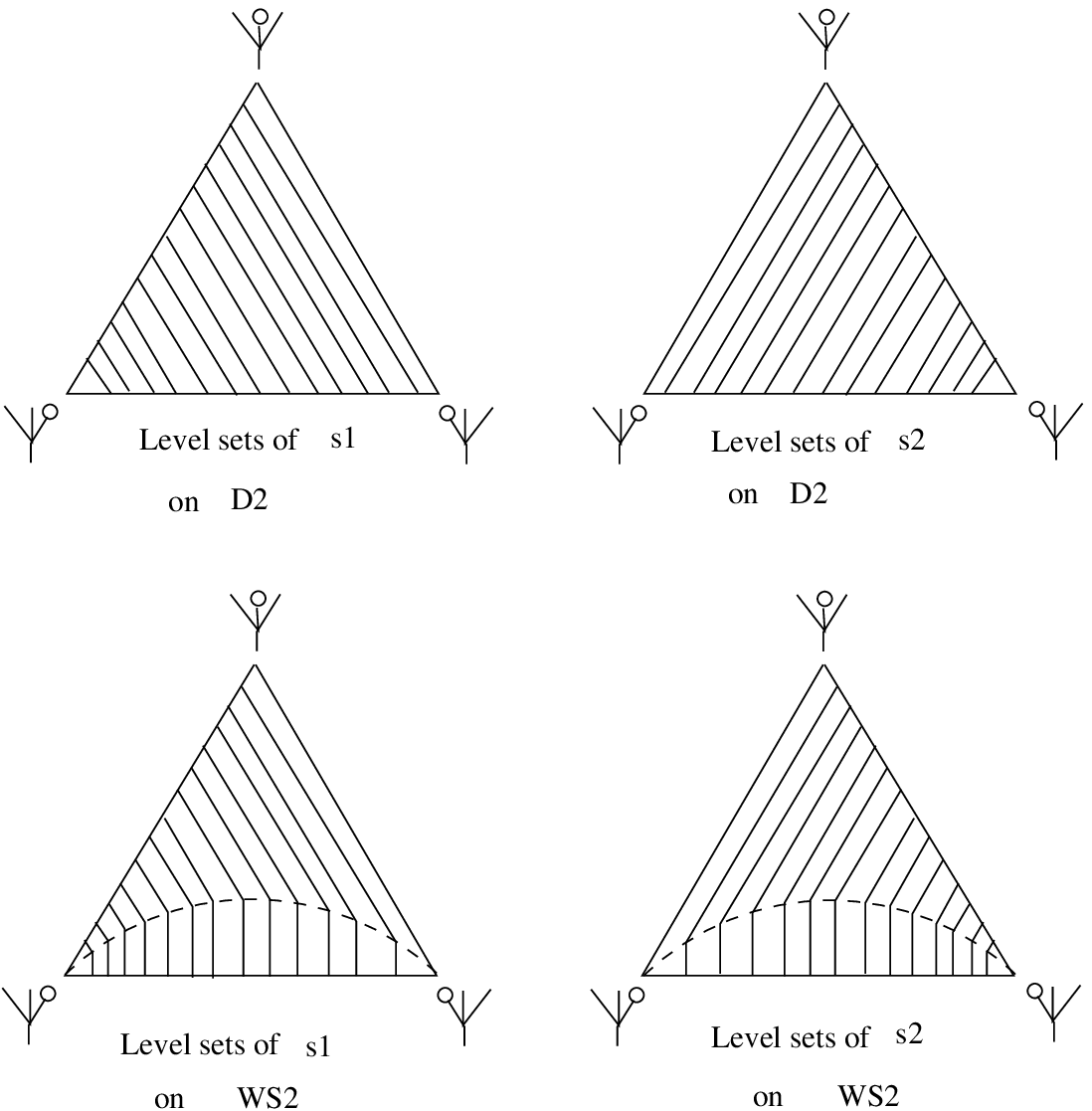}
\caption{}\label{fig8}
\end{figure}

Even though $\triangle$ and $Wb\square$ are homeomorphic as weak $\Assoc_{>0}$ bimodules (Proposition~\ref{p:w_square}), they are no more homeomorphic as weak $\Assoc$ bimodules. For example, one can compare the degeneracy maps $s_1(-):=(-)\circ_1 a_0$ and $s_2(-):=(-)\circ_2 a_0$ on $\triangle(2)$ and $Wb\square(2)$, see Figure~\ref{fig8}.

In the case of $\triangle(2)$ a level set of $s_1$ can have at most one point of intersection with a level set of $s_2$. In the case of $Wb\square(2)$ such level sets can have a segment in common. This comes from the fact that the maps
\begin{gather*}
\circ_1\colon \square(2)\times \Assoc(0)\to \square(1),\\
\circ_2\colon \square(2)\times \Assoc(0)\to \square(1),
\end{gather*}
are the same maps since $\square(1)$ is a point.

\subsection{Cofibrant weak $\Assoc$ bimodule $\widetilde{Wb\square}$ and associated tower $T_\bullet^{\widetilde{Wb\square}}(\calO)$}\label{ss:cof_tilde_wb}
The construction of $\widetilde{\triangle}$ and $\widetilde{\square}$ works in the same way to make a cofibrant weak $\Assoc$ bimodule $\widetilde{Wb\square}$ from $Wb\square$:
$$
\widetilde{Wb\square}(n):=\left.\left(\coprod_{m\geq 0}\,\,\,\coprod_{\alpha\colon\underline{m}\hookrightarrow\underline{n+m}}Wb\square(n+m)\times [0,1]^{\alpha(\underline{m})}\right)\right/\sim.
\eqno(\numb)\label{eq_tilde_wsq_disj}
$$
The quotient relations are the same as in~\eqref{eq:Yrel1},~\eqref{eq:Yrel2},~\eqref{eq:Yrel3}. In addition one has
\begin{gather}
(a_2\circ_2y;\tau,\tau_{\alpha(1)},\ldots,\tau_{\alpha(m)})_{\alpha_{1^0}}\sim(y;\tau_{\alpha(1)},\ldots,\tau_{\alpha(m)})_\alpha;\label{eq:Yrel5}\\
(a_2\circ_1y;\tau_{\alpha(1)},\ldots,\tau_{\alpha(m)},\tau)_{\alpha_{(m+n+1)^0}}\sim(y;\tau_{\alpha(1)},\ldots,\tau_{\alpha(m)})_\alpha\label{eq:Yrel6},
\end{gather}
where $\alpha_{1^0}\colon \underline{m+1}\hookrightarrow \underline{m+n+1}$ has image $\{1\}\cup (\alpha(\underline{m})+1)$, and $\alpha_{(m+n+1)^0}\colon
\underline{m}\hookrightarrow \underline{m+n+1}$ has image $\alpha(\underline{m})\cup\{m+n+1\}$.

%
%
%

 Since this construction is fonctorial,  Filtration~\eqref{eq:filt_w_sq} in $Wb\square$ induces filtration
$$
\widetilde{Wb\square}_0\subset \widetilde{Wb\square}_1\subset \widetilde{Wb\square}_2\subset \widetilde{Wb\square}_3\subset \ldots.
\eqno(\numb)\label{eq_tilde_wsquare_filtr}
$$

The following proposition is analogous to Proposition~\ref{p:cof_tilde_tr}.

\begin{proposition}\label{p:cof_tilde_w_sq}
\textup{(i)} $\widetilde{Wb\square}$ is a cofibrant model of $\Assoc$ as a weak bimodule over itself.

\textup{(ii)} for any $n\geq 0$, $\widetilde{Wb\square}_n|_n$ is a cofibrant model of $\Assoc|_n$ as an $n$-truncated weak $\Assoc$ bimodule.
\end{proposition}

\begin{remark}\label{r:fiber_tower_w_sq}
The weak bimodule $\widetilde{Wb\square}$ has a coarsening which is obtained by the same combinatorial sequence of cell attachments as $\widetilde{\triangle}$.
To recall $Wb\square(n)$ is an $n$-disc. Similarly to $\widetilde{\triangle}$, see Section~\ref{s:deg_wb}, one will have that the interior cell of each $Wb\square(n+m)\times[0,1]^{\alpha(\underline{m})}$ in~\eqref{eq_tilde_wsq_disj} is a generating cell of $\widetilde{Wb\square}$.
To pass from $\widetilde{Wb\square}_{n-1}$ to $\widetilde{Wb\square}_{n}$ one needs to attach one cell of dimension $n$ in degree $n$, then $n$ cells of dimension $n+1$ in  degree $n-1$,$\ldots$, then ${n\choose i}$ cells of dimension $n+i$ in degree $n-i$ (the corresponding cells are the interiors of $Wb\square(n)\times[0,1]^{\alpha(\underline{i})}$; we will denote by $\widetilde{Wb\square}_{n,i}$ the intermediate weak bimodule obtained at this point and which is similar to $\widetilde{\triangle}_{n,i}$),$\ldots$, then $1$ cell of dimension $2n$ in degree~0.
\end{remark}

As before for any $\Assoc$ weak bimodule $\calO$ we obtain the tower whose stages are
$$
T_n^{\widetilde{Wb\square}}(\calO):=\underset{\Assoc}{\WBimod}(\widetilde{Wb\square}_n,\calO)=\underset{\Assoc}{\WBimod}{}_n(\widetilde{Wb\square}_n|_n,\calO).
\eqno(\numb)\label{eq:stages_tilde_wsquare}
$$

\begin{proposition}\label{p:equiv_tild_wb}
One has that $\hoTriangle$ is homotopy equivalent to $\widetilde{Wb\square}$ as a filtered weak bimodule over $\Assoc$. As a consequence  the towers $T_\bullet^{\widetilde{Wb\square}}(\calO)$ and $T_\bullet^{\widetilde{\triangle}}(\calO)$ are also homotopy equivalent.
\end{proposition}

\begin{proof}
This follows immediately from Propositions~\ref{p:cof_tilde_tr},~\ref{p:cof_tilde_w_sq}, and Lemma~\ref{l:trunc_maps}.
\end{proof}

We will also need intermediate mapping spaces
$$
T_{n,i}^{\widetilde{Wb\square}}(\calO):=\underset{\Assoc}{\WBimod}(\widetilde{Wb\square}_{n,i},\calO)
\eqno(\numb)\label{eq:substages_tilde_wsq}
$$
that will be used in the proof of Theorem~\ref{t:deloop1_stages_tilde}.

\subsection{Weak bimodule $Wb\widetilde{\square}$}\label{ss:f_deloop_tilde}
In this section we will describe an alternative construction for $\widetilde{Wb\square}$. In this construction it will be denoted by $Wb\widetilde{\square}$. Using this alternative description of $\widetilde{Wb\square}=Wb\widetilde{\square}$ we construct the homotopy equivalence of towers $\Omega\, T_\bullet^{\widetilde{\square}}(\calO)\stackrel{\tilde\xi_\bullet}{\longrightarrow}T_\bullet^{\widetilde{Wb\square}}(\calO)$.

Recall the category $\Xi_n$. Let $\widetilde{\Xi}_n$ denote a similar category of trees with the same morphisms, but now we allow trees with univalent and bivalent beads.

\begin{figure}[h]
\includegraphics[width=8cm]{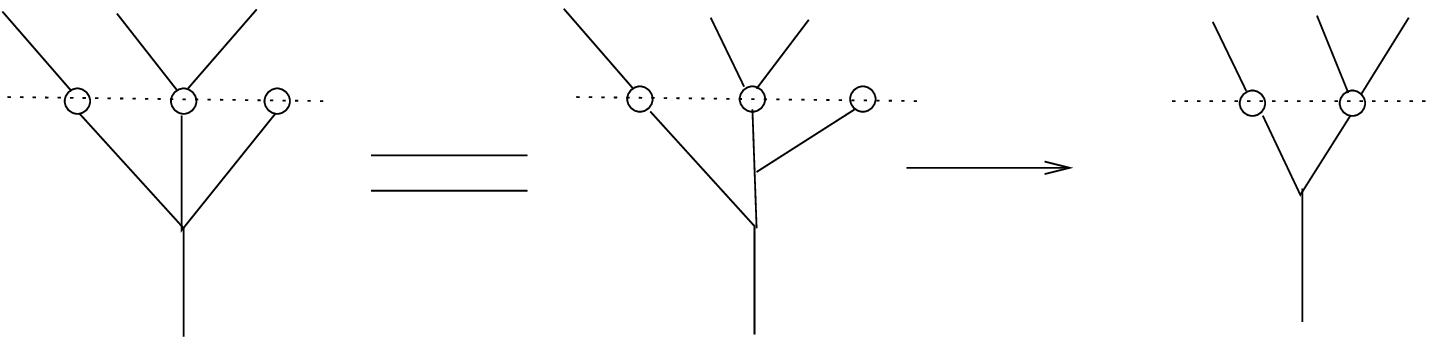}
\caption{An example of a morphism in $\widetilde{\Xi}_3$.}\label{fig9}
\end{figure}

The category $\widetilde{\Xi}_n$ has infinitely many elements. Using $\widetilde{\square}$ as a bimodule over $\Assoc$, we will define a covariant functor
$$
\lambda_{\widetilde{\square}}\colon\widetilde{\Xi}_n\to\Top
$$
by sending a tree $T$ to $\prod_{b\in B(T)}\widetilde{\square}(|b|)$, where $|b|$ denotes the number of outgoing vertices of a bead $b$. On the morphisms this functor is defined using the left action of $\Assoc_{>0}$ on $\widetilde{\square}$.

We also define a contravariant functor
$$
\chi_\blacktriangle\colon\widetilde{\Xi}_n\to\Top,
$$
which assigns to any tree $T$ in $\widetilde\Xi_n$ the simplex $\Delta^{B(T)}$, which is the space of order preserving maps from the set $B(T)$ of beads to $[0,1]$.

Then we define $Wb\widetilde{\square}(n)$ as
$$
\left.\left(\lambda_{\widetilde{\square}}\otimes_{\widetilde{\Xi}_n}\chi_\blacktriangle\right)\right/\sim,
$$
where the equivalence relations are the same as in Section~\ref{s:first_deloop} (defined after Remark~\ref{r:ev_w_sq}) plus in addition we impose
$$
(x_1,\ldots,x_i=\unit,\ldots,x_k;t_1,\ldots,t_k)\sim (x_1,\ldots,\widehat{x_{i}},\ldots,x_k;t_1,\ldots,\widehat{t_{i}},\ldots,t_k).
\eqno(\numb)\label{eq_tilde_wsq_unit2}
$$

Similarly as we defined a structure of a weak $\Assoc_{>0}$ bimodule on $Wb\square$, one defines a structure of a weak $\Assoc$ bimodule on $Wb\widetilde{\square}=\{Wb\widetilde{\square}(n),\, n\geq 0\}$.

\begin{proposition}\label{p:equiv_wtilde=tildew}
The weak $\Assoc$ bimodules $Wb\widetilde{\square}$ and $\widetilde{Wb\square}$ are naturally homeomorphic.
\end{proposition}

\begin{proof}
The homeomorphism is clear from the construction. In both cases a typical element of degree $n$ can be written as a tuple
 $((x_1,\ldots,x_k);(t_1,\ldots,t_k),(\tau_{\alpha(1)},\ldots,\tau_{\alpha(m)}))_\alpha$,
where $x_i\in\square(n_i)$, $n_i\geq 1$; $(t_1,\ldots,t_k)\in\Delta^k$; $\alpha\colon\underline{m}\hookrightarrow\underline{m+n}$,
 $(\tau_{\alpha(1)},\ldots,\tau_{\alpha(m)})\in [0,1]^{\alpha(\underline{m})}$, and one also has $n_1+\ldots +n_k=n+m$. It is easy to check that the equivalence relations for
 $\widetilde{Wb\square}$ and for $Wb\widetilde{\square}$ are  the same. As example the relation~\eqref{eq:Yrel1} together with~\eqref{eq_tilde_wsq_unit1} imply~\eqref{eq_tilde_wsq_unit2}.
\end{proof}

Starting from this point we will not distinguish between $\widetilde{Wb\square}$ and $Wb\widetilde{\square}$. In particular one will always mean
$$
Wb\widetilde{\square}_n=\widetilde{Wb\square}_n,\quad
Wb\widetilde{\square}_{n,i}=\widetilde{Wb\square}_{n,i},\quad
T_n^{Wb\widetilde{\square}}(\calO)=T_n^{\widetilde{Wb\square}}(\calO),\quad
T_{n,i}^{Wb\widetilde{\square}}(\calO)=T_{n,i}^{\widetilde{Wb\square}}(\calO).
$$

\subsection{$Wb$-construction and first delooping}\label{ss:tilde_w_constr}

Now let $\calO$ be an $\Assoc$ bimodule endowed with a morphism of $\Assoc$ bimodules
$$
p\colon\Assoc\to\calO.
$$
Due to this map, $\calO$ inherits a structure of a weak $\Assoc$ bimodule.

Given any $\Assoc$ bimodule $\calQ$ one can produce a weak $\Assoc$-bimodule $Wb(\calQ)$ in the same way as in the previous section we constructed $Wb\widetilde{\square}$ from $\widetilde{\square}$. The $Wb$ construction has the property that one can define a natural map
$$
 \tilde\xi_\calQ\colon\Omega\,\underset{\Assoc}{\Bimod}(\calQ,\calO)\to \underset{\Assoc}{\WBimod}(Wb(\calQ),\calO),
$$
where the base point in $\underset{\Assoc}{\Bimod}(\calQ,\calO)$ is the composition $\calQ\to\Assoc\stackrel{p}{\to}\calO$.

Applying the $Wb$-construction to $\widetilde{\square}_n$ and $\widetilde{\square}_{n,i}$ one obtains the weak bimodules $Wb(\widetilde{\square}_n)$ and $Wb(\widetilde{\square}_{n,i})$ that are regarded as subobjects of $Wb\widetilde{\square}$. Due to the inclusions $Wb\widetilde{\square}_n\subset Wb(\widetilde{\square}_n)$ and $Wb\widetilde{\square}_{n,i}\subset Wb(\widetilde{\square}_{n,i})$ the maps $\tilde\xi_{\widetilde{\square}_n}$ and $\tilde\xi_{\widetilde{\square}_{n,i}}$ can be composed with the restriction maps to get
\begin{gather*}
\tilde\xi_n\colon \Omega T_n^{\widetilde{\square}}(\calO)=\Omega\, \underset{\Assoc}{\Bimod}(\widetilde{\square}_n,\calO)\to\underset{\Assoc}{\WBimod}(Wb\widetilde{\square}_n,\calO)=T_n^{Wb\widetilde{\square}}(\calO),\\
\tilde\xi_{n,i}\colon \Omega T_{n,i}^{\widetilde{\square}}(\calO)=\Omega\, \underset{\Assoc}{\Bimod}(\widetilde{\square}_{n,i},\calO)\to\underset{\Assoc}{\WBimod}(Wb\widetilde{\square}_{n,i},\calO)=T_{n,i}^{Wb\widetilde{\square}}(\calO).
\end{gather*}


\begin{theorem}\label{t:deloop1_stages_tilde}
Each map $\tilde\xi_n\colon \Omega\, T_n^{\hoSquare}(\calO)\to T_n^{Wb\widetilde{\square}}(\calO)$ is a homotopy equivalence for any $n$ and any pointed $\Assoc$ bimodule $\calO$ with $\calO(0)\simeq *$.
\end{theorem}

Above by a \textit{pointed} $\Assoc$ bimodule we mean a bimodule endowed with a morphism from $Assoc$. In other words it is a bimodule in the category of pointed spaces.

\begin{proof}
We will prove that in the settings of the theorem all the maps $\tilde\xi_{n,i}$ are homotopy equivalences. The result will follow since $\tilde\xi_n=\tilde\xi_{n-1,n}$.

For a convenience we will use the notation  $\tilde\xi_{n-1,-1}:=\tilde\xi_{n-1}$, $T_{n-1,-1}^{\widetilde{\square}}(\calO):=T_{n-1}^{\widetilde{\square}}(\calO)$, $\widetilde{\square}_{n-1,-1}:=\widetilde{\square}_{n-1}$, and so on.

One has a commutative diagram
$$
\xymatrix{
\Omega T_{n-1,i}^{\widetilde{\square}}(\calO)\ar[r]^{\tilde\xi_{n-1,i}}\ar[d]&T_{n-1,i}^{Wb\widetilde{\square}}(\calO)\ar[d]\\
\Omega T_{n-1,i-1}^{\widetilde{\square}}(\calO)\ar[r]^{\tilde\xi_{n-1,i-1}}&T_{n-1,i-1}^{Wb\widetilde{\square}}(\calO)
}.
$$

Assuming that $\tilde\xi_{n-1,i-1}$ is a homotopy equivalence we will prove that $\tilde\xi_{n-1,i}$ is also a homotopy equivalence. Since both
  $\Omega T_{n-1,i}^{\widetilde{\square}}(\calO)\to \Omega T_{n-1,i-1}^{\widetilde{\square}}(\calO)$ and $T_{n-1,i}^{Wb\widetilde{\square}}(\calO)\to T_{n-1,i-1}^{Wb\widetilde{\square}}(\calO)$ are fibrations it is enough to show that $\tilde\xi_{n-1,i}$ induces a homotopy equivalence on fibers.

  To recall $Wb\widetilde{\square}_{n-1,i}$ is obtained from $Wb\widetilde{\square}_{n-1,i-1}$ by a free attachment of $n\choose i$ cells in degree $n-i$. These cells are exactly the interior of $\coprod_{\alpha\colon\underline{i}\hookrightarrow\underline{n}}Wb\square(n)\times[0,1]^{\alpha(\underline{i})}$ in~\eqref{eq_tilde_wsq_disj}. Consider
$$
Wb\widetilde{\square}_{n-1,i-1/2}:=Wb(\widetilde{\square}_{n-1,i-1})\cap Wb\widetilde{\square}_{n-1,i}.
$$
One can see that $Wb\widetilde{\square}_{n-1,i-1/2}$ is obtained from $Wb\widetilde{\square}_{n-1,i-1}$ by a free attachment of $n\choose i$ punctured discs
$\coprod_{\alpha\colon\underline{i}\hookrightarrow\underline{n}}\left(Wb\square(n)\setminus Int(Wb\square(C_n))\right)\times[0,1]^{\alpha(\underline{i})}$.

On the other hand it is easy to see that given $\tilde g\in \Omega T_{n-1,i-1}^{\hoSquare}(\calO)$ the fiber over $\tilde g$   in $\Omega T_{n-1,i}^{\hoSquare}(\calO)$ is
homeomorphic to the space of maps $\underset{\Assoc}{\WBimod}\, _{\tilde\xi_{n-1}^{\tilde g}}\left((Wb\widetilde{\square}_{n-1,i},Wb\widetilde{\square}_{n-1,i-1/2}),\calO\right)$. While the fiber over its image $\tilde\xi_{n-1,i-1}^{\tilde g}\in
T_{n-1,i-1}^{Wb\hoSquare}(\calO)$ in $T_{n-1,i-1}^{Wb\hoSquare}(\calO)$ is the space $\underset{\Assoc}{\WBimod}\, _{\tilde\xi_{n-1}^{\tilde g}}\left((Wb\widetilde{\square}_{n-1,i},Wb\widetilde{\square}_{n-1,i-1}),\calO\right)$. It follows from Lemma~\ref{l:fiber_equiv}, that the inclusion
$$
\underset{\Assoc}{\WBimod}\, _{\tilde\xi_{n-1}^{\tilde g}}\left((Wb\widetilde{\square}_{n-1,i},Wb\widetilde{\square}_{n-1,i-1/2}),\calO\right)\hookrightarrow
\underset{\Assoc}{\WBimod}\, _{\tilde\xi_{n-1}^{\tilde g}}\left((Wb\widetilde{\square}_{n-1,i},Wb\widetilde{\square}_{n-1,i-1}),\calO\right)
$$
is a homotopy equivalence.
\end{proof}

\begin{theorem}\label{t:first_deloop_tilde}
In the settings of Theorem~\ref{t:deloop1_stages_tilde} the induced map of limits
$$
\Omega\,\underset{\Assoc}{\Bimod}(\hoSquare,\calO)\stackrel{\tilde\xi_\infty}{\longrightarrow}\underset{\Assoc}{\WBimod}(Wb\hoSquare,\calO).
\eqno(\numb)\label{eq:first_deloop_xi_tilde}
$$
is a homotopy equivalence.
\end{theorem}
\begin{proof}
The proof is similar to that of Theorem~\ref{t:first_deloop}.
\end{proof}

\section{Second delooping and the operad {\sc $\widetilde{\Pentagon}$}}\label{s:second_deloop_tilde}
\subsection{Operad {\sc $\widetilde{\Pentagon}$}}\label{ss:tilde_pentagon}

It is possible to add a point in the zero component of the operad $\Pentagon$ so that the new sequence of spaces $\Pentagon=\{\Pentagon(n),\, n\geq 0\}$ will still have a structure of an operad. But it turns out that if we will then  try to use the method of Subsection~\ref{ss:tilde_square} to make a cofibrant operad from $\Pentagon$, then the obtained sequence of spaces
$$
\left.\left(\coprod_{m\geq 0}\,\,\,\coprod_{\alpha\colon\underline{m}\hookrightarrow\underline{n+m}}\Pentagon(n+m)\times[0,1]^{\alpha(\underline{m})}\right)\right/\sim
$$
(where the equivalence relation is defined as in Subsection~\ref{ss:tilde_square}) will not have a natural structure of an operad.

On the other hand we can use the Boardmann-Vogt $W$-construction~\cite{BoaVog} to produce a cofibrant replacement of the operad $\Assoc$. The $n$-th component $\widetilde{\Pentagon}(n)$ of this cofibrant operad is the space of planar metric trees with $n$ leaves. A point in this space is a planar rooted tree whose inner vertices can have any valence and edges have lengths between 0 and 1.  The root in such trees is supposed to have valence~1. Neither the leaf edges, nor the root one are assigned a length. There are two relations:  If the length of an edge becomes~0, this edge is contracted. Another relation is that if a tree has an inner vertex of valence~2, then this tree is equivalent to the tree obtained by removing this vertex and making a single edge from the two ones adjacent to this vertex. If both edges were internal ones, the length of the new edge will be the maximum of two lengths, otherwise the edge is not assigned any length since it must a leaf edge or a root one.
The latter relation corresponds to the fact that the vertex of valence 2 corresponds to the identity element of the operad $\Assoc$.

\begin{figure}[h]
\includegraphics[width=3cm]{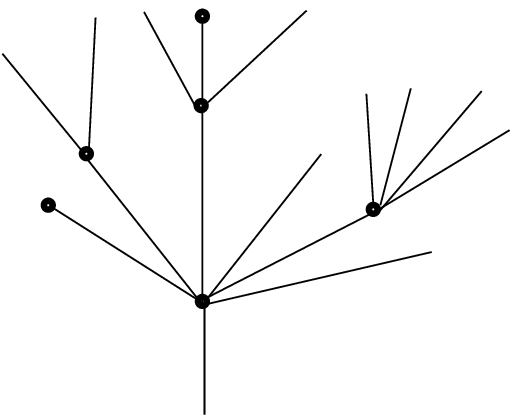}
\caption{An example of a tree from $\widetilde{\Pentagon}(10)$.}\label{fig10}
\end{figure}

We will depict the internal vertices of metric trees by black dots to distinguish these trees from the ones previously used in the paper.

The operadic composition is defined by grafting trees and assigning length~1 to the edges connecting the trees in the composition. The identity element with respect to this operadic structure is the tree having only one leaf and no internal vertices.

The subspaces of $\widetilde{\Pentagon}(n)$, $n\geq 0$, consisting of trees whose all internal vertices have valence $\geq 3$ form a suboperad which is naturally isomorphic to $\Pentagon$ (since it is also the Boardmann-Vogt resolution of $\Assoc_{>0}$).

We say that a tree is {\it prime} if it does not have edges of length~1 and also the tree is not the identity tree. If the tree is not prime and is not the identity it is said {\it composite}. In other words a composite tree is a tree having at least one edge of length~1. Any composite tree can be uniquely decomposed into an operadic composition of prime trees called {\it prime components}. We say that a prime tree is in the $i$-th filtration term if the number of its leaves plus the number of its univalent internal vertices is $\leq i$.
A composite tree is said to be in the $i$-th filtration term if all its prime components are in the $i$-th filtration term. This defines a filtration of operads in $\widetilde{\Pentagon}$:
$$
\widetilde{\Pentagon}_0\subset\widetilde{\Pentagon}_1\subset\widetilde{\Pentagon}_2\subset\widetilde{\Pentagon}_3\subset\ldots\subset\widetilde{\Pentagon}
$$
All the components of the term $\widetilde{\Pentagon}_0$ are empty except the  one in degree one $\widetilde{\Pentagon}_0(1)=\{id\}$. It is a free operad generated by empty sets. 
The first term $\widetilde{\Pentagon}_1$ has in addition one degree zero operation $\widetilde{\Pentagon}_1(0)=\{\,
\raisebox{-4pt}{\includegraphics[width=.13cm]{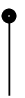}}\,
\}$. The other components are the same. In other words $\widetilde{\Pentagon}_1$ is obtained from $\widetilde{\Pentagon}_0$ by attaching one zero-cell in degree zero. The term
$\widetilde{\Pentagon}_2$ is obtained from $\widetilde{\Pentagon}_1$ by attaching one 0-cell in degree~2, corresponding to the tree~
\raisebox{-4pt}{\includegraphics[width=.4cm]{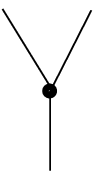}}, then attaching two 1-cells in degree 1 that consist of metric trees~
\raisebox{-4pt}{\includegraphics[width=.4cm]{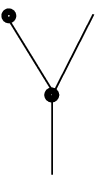}} and~\raisebox{-4pt}{\includegraphics[width=.4cm]{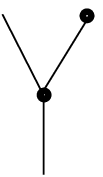}}, then attaching one 2-cell in degree~0 that consist of trees~\raisebox{-4pt}{\includegraphics[width=.4cm]{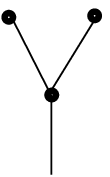}}. In general $\widetilde{\Pentagon}_{n}$ is obtained from $\widetilde{\Pentagon}_{n-1}$, $n\geq 2$, by a sequence of free cell attachments. We attach first one $(n-2)$-cell in degree~$n$ -- the interior of the corresponding cell consists of prime trees with $n$ leaves and without univalent internal vertices. This cell is exactly the interior of~$\Pentagon(n)$. Then one attaches  $n$ $(n-1)$-cells in degree $(n-1)$. The interior of these cells consist of prime trees with $(n-1)$ leaves and one internal univalent vertex. Each such cell is naturally described by $\Pentagon(n)\times [0,1]$. On the $i$-th step of this procedure we freely attach $n\choose i$ $(n-2+i)$-cells in degree $(n-i)$ (the interiors of these cells consist of prime trees with $(n-i)$ leaves and $i$ univalent beads;  the operad obtained at this point will be denoted by $\widetilde{\Pentagon}_{n-1,i}$).   One has $\widetilde{\Pentagon}_{n-1,n}=\widetilde{\Pentagon}_n$.

\begin{proposition}\label{p:operad_cof_tilde}
\textup{(i)} The operad $\widetilde{\Pentagon}$ is a cofibrant model of $\Assoc$.

\textup{(ii)} For $n\geq 1$, the $n$-truncated operad $\widetilde{\Pentagon}_n|_n$ is a cofibrant model of $\Assoc|_n$.

\end{proposition}

\begin{proof}
An explicit sequence of cell attachment is described above. The credit for the result~(i) should be given to Boardmann-Vogt~\cite{BoaVog}, since what we described above  is a partial case of their construction. For the part~(ii) according to Lemma~\ref{l:truncation} we are only left to show that $\widetilde{\Pentagon}_n(k)$ is contractible for any $k\leq n$. We say that a metric tree $T\in\widetilde{\Pentagon}_n(k)$ is {\it quasi-prime} if each of its edges of length~1 is adjacent to an inner vertex of valence~1 and also $T\neq \raisebox{-4pt}{\includegraphics[width=.13cm]{b0.eps}}$. Every metric tree (except \raisebox{-4pt}{\includegraphics[width=.13cm]{b0.eps}}) has a unique operadic decomposition into quasi-prime components. Define a filtration in  $\widetilde{\Pentagon}_n(k)$ whose $i$-th term $(\widetilde{\Pentagon}_n(k))_i$ consists of metric trees whose all quasi-prime components have $\leq i$ inner vertices of valence~1. By contracting all the edges adjacent to inner vertices of valence~1 one can show that $(\widetilde{\Pentagon}_n(k))_i$ can be retracted to $(\widetilde{\Pentagon}_n(k))_{i-1}$. Since the 0-th filtration term
$$
(\widetilde{\Pentagon}_n(k))_0=
\begin{cases}
\{\,\raisebox{-4pt}{\includegraphics[width=.13cm]{b0.eps}}\,\}, & k=0;\\
\Pentagon_n(k)=\Pentagon(k), & k>1;
\end{cases}
$$
is contractible, we obtain by induction that all the filtration terms are also contractible. On the other hand each inclusion of filtration terms is an NDR pair (being an inclusion of $CW$-complexes), as a consequence the colimit (union) of this filtration sequence coincides up to a homotopy equivalence with the homotopy colimit. The latter one is contractible since all the terms in the filtration are contractible.
\end{proof}

To end the section we notice that for each $n\geq 0$ we have a natural surjective map

$$
\sigma_{n,\widetilde{\pentagon}}\colon
\left(\coprod_{m\geq 0}\,\,\,\coprod_{\alpha\colon\underline{m}\hookrightarrow\underline{n+m}}\Pentagon(n+m)\times[0,1]^{\alpha(\underline{m})}\right)\to\widetilde{\Pentagon}(n),
\eqno(\numb)\label{eq:tilde_pent}
$$
that sends $(T;\tau_{\alpha(1)},\ldots,\tau_{\alpha(m)})_\alpha$ to a  tree obtained from a metric tree $T\in\Pentagon(n+m)$ by corking its entries $\alpha(1),\ldots,\alpha(m)$ by the segments of the lengths $\tau_{\alpha(1)},\ldots,\tau_{\alpha(m)}$ respectively (this produces $m$ new univalent internal vertices in $T$). 

This means that we can also view $\widetilde{\Pentagon}(n)$ as
$$
\widetilde{\Pentagon}(n)=\left.\left(\coprod_{m\geq 0}\,\,\,\coprod_{\alpha\colon\underline{m}\hookrightarrow\underline{n+m}}\Pentagon(n+m)\times[0,1]^{\alpha(\underline{m})}\right)\right/\sim.
\eqno(\numb)\label{eq:tPentagon_as_union}
$$
But the equivalence relations here are different from those in Section~\ref{ss:tilde_square}. Consider for example the tree below
$$
\psfrag{t1}[0][0][1][0]{$\tau_1$}
\psfrag{t2}[0][0][1][0]{$\tau_2$}
\psfrag{t1=0}[0][0][1][0]{$\tau_1\to 0$}
\psfrag{t2=0}[0][0][1][0]{$\tau_2\to 0$}
\includegraphics[width=6cm]{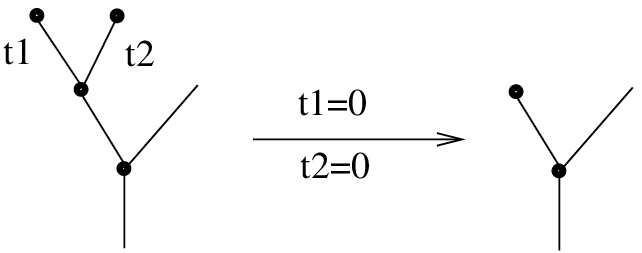}
$$
which is an interior point of $\Pentagon(3)\times[0,1]^{\alpha(\underline{2})}$, where $\alpha\colon\underline{2}\to\underline{3}$ is defined by $\alpha(i)=i$, $i=1,\,2$.
When $\tau_1$ and $\tau_2$ become~0 the limiting tree will still have a univalent internal vertex. Which means that even when all $\tau_{\alpha(i)}$ are zero the attaching map does not go to some point in $\Pentagon(n)$, but to some other cell of the form $\Pentagon(n+m')\times[0,1]^{\alpha'(\underline{m'})}$ with $m'<m$.
This actually explains why the proof of Proposition~\ref{p:operad_cof_tilde}~(2) is more complicated compared to the proof of the analogous statement~\ref{p:cof_tilde_tr} in which $\hoTriangle_n$ is immediately retracted to $\triangle_n$.

Again we have that the interior of each $\Pentagon(n+m)\times [0,1]^{\alpha(\underline{m})}$, $n+m>1$, in~\eqref{eq:tPentagon_as_union} is a generating cell of $\hoPentagon$ of dimension $n+2m-2$.

\subsection{Tower associated to {\sc $\widetilde{\Pentagon}$}}\label{ss:tower_tilde_pentag}
Given a topological operad $\calO$ one can define a tower
$$
 T_1^{\widetilde{\pentagon}}(\calO)\leftarrow T_2^{\widetilde{\pentagon}}(\calO)\leftarrow T_3^{\widetilde{\pentagon}}(\calO)\leftarrow\ldots,
$$
where each stage is defined by
$$
T_n^{\widetilde{\pentagon}}(\calO):=\Operad(\widetilde{\Pentagon}_n,\calO)=\Operad_n(\widetilde{\Pentagon}_n|_n,\calO|_n)\simeq
\widetilde{\Operad}_n(\Assoc|_n,\calO|_n).
$$
The (homotopy) limit of the above tower is $T_\infty^{\widetilde{\pentagon}}(\calO)=\Operad(\widetilde{\Pentagon},\calO)\simeq\widetilde{\Operad}(\Assoc,\calO)$.

We will describe the fibers of the maps between stages in terms of the spaces of morphisms of cubical diagrams. Propositions~\ref{p:fiber_tilde_pentag},~\ref{p:fiber_tilde_pentag2} below are given for completeness of exposition and will not be used in the sequel.

Let us fix an element $p$ in $\calO(0)$. Using this element one can define a cubical diagram $\calO(\underline{n}\setminus\bullet)$, where each map in the diagram is a composition with $p$. In practice $p$ will be the image of~\raisebox{-4pt}{\includegraphics[width=.13cm]{b0.eps}} under a map $\gamma\colon\widetilde{\Pentagon}_1\to\calO$.

We will also consider the cubical diagram $\widetilde{\Pentagon}(\underline{n}\setminus\bullet)$ where~\raisebox{-4pt}{\includegraphics[width=.13cm]{b0.eps}} is the element with which we consider the compositions in order to define the maps in the cubical diagram.

 Let $\gamma\in T_{n-1}^{\widetilde{\pentagon}}(\calO)$. By abuse of notation we will denote by $\gamma$ the induced morphism of cubes:
$$
\gamma\colon\widetilde{\Pentagon}_{n-1}(\underline{n}\setminus\bullet)\to\calO(\underline{n}\setminus\bullet).
$$

Let $\Pentagon(n)\times [0,1]^\bullet$ denote the $n$-cube obtained from $[0,1]^\bullet$ by taking the product of each object with the space $\Pentagon(n)$. We will also consider $\partial_0\left(\Pentagon(n)\times [0,1]^\bullet\right)$ -- a cubical subobject of  $\Pentagon(n)\times [0,1]^\bullet$ defined by
$$
\partial_0\left(\Pentagon(n)\times [0,1]^\bullet\right)=
\Bigl(\partial\bigl(\Pentagon(n)\bigr)\times [0,1]^S\Bigr)\cup \Bigl(\Pentagon(n)\times \left([0,1]^S\right)_0\Bigr).
$$

Let $\beta_S$ denote the map
$$
\Pentagon(n)\times [0,1]^S\to \widetilde{\Pentagon}(\underline{n}\setminus S),
$$
which is the composition of the inclusion
$$
i_S\colon \Pentagon(n)\times [0,1]^S\hookrightarrow \coprod_{m\geq 0}\,\,\,\coprod_{\alpha\colon\underline{m}\hookrightarrow\underline{n-|S|+m}} \Pentagon(n-|S|+m)\times [0,1]^{\alpha(\underline{m})}
$$
(that sends $\Pentagon(n)\times [0,1]^S$ identically to the component labeled by the inclusion $\underline{|S|}\simeq S\subset\underline{n}$), and the quotient map
$$
\coprod_{m\geq 0}\,\,\,\coprod_{\alpha\colon\underline{m}\hookrightarrow\underline{n-|S|+m}} \Pentagon(n-|S|+m)\times [0,1]^{\alpha(\underline{m})}\to\Bigl(\coprod_{m\geq 0}\,\,\,\coprod_{\alpha\colon\underline{m}\hookrightarrow\underline{n-|S|+m}} \Pentagon(n-|S|+m)\times [0,1]^{\alpha(\underline{m})}\Bigr)\Bigl/\sim,
$$
where the target space is $\widetilde{\Pentagon}(n-|S|)$. The maps $\beta_S$, $S\subset\underline{n}$ define a natural transformation
$$
\beta\colon \Pentagon(n)\times [0,1]^\bullet\to \widetilde{\Pentagon}(\underline{n}\setminus \bullet).
$$
Notice that the preimage of the subcube $\widetilde{\Pentagon}_{n-1}(\underline{n}\setminus\bullet)\subset \widetilde{\Pentagon}(\underline{n}\setminus\bullet)$ under $\beta$ is exactly $\partial_0\left(\Pentagon(n)\times [0,1]^\bullet\right)$.

\begin{proposition}\label{p:fiber_tilde_pentag}
The preimage of $\gamma\in T_{n-1}^{\widetilde{\pentagon}}(\calO)$ under the fibration $T_n^{\widetilde{\pentagon}}(\calO)\to T_{n-1}^{\widetilde{\pentagon}}(\calO)$,
$n\geq 2$,\footnote{One also has $T_{1}^{\widetilde{\pentagon}}(\calO)=\calO(0)$.} is homeomorphic to
$$
Nat_{\gamma\circ\beta}\Bigl(\bigl(\Pentagon(n)\times[0,1]^\bullet,\partial_0\left(\Pentagon(n)\times [0,1]^\bullet\right)\bigr),\calO(\underline{n}\setminus\bullet)\Bigr).
$$
\end{proposition}

\begin{proof} Direct inspection.
\end{proof}

The following statement can be easily proved from Proposition~\ref{p:fiber_tilde_pentag}.

\begin{proposition}\label{p:fiber_tilde_pentag2}
The preimage of the fibration $T_n^{\widetilde{\pentagon}}(\calO)\to T_{n-1}^{\widetilde{\pentagon}}(\calO)$, $n\geq 3$,  is either empty or homotopy equivalent to $\Omega^{n-2} tfiber(\calO(\underline{n}\setminus\bullet))$, where for the degeneracy maps in the cubical diagram one chooses compositions with $p\in\calO(0)$ which is the image of~\raisebox{-4pt}{\includegraphics[width=.13cm]{b0.eps}}, for the base point in $\calO(n)$ one can  choose $\underbrace{q\circ_1\ldots\circ_1 q}_{n-1}$ with $q$ being the image of~\raisebox{-4pt}{\includegraphics[width=.4cm]{b2.eps}} in $\calO(2)$.\footnote{For $n=2$ the preimage can be described as the space $tfiber(\calO(\underline{2}\setminus\bullet))$ where each $\calO(1)$ is pointed in $id$, and $\calO(0)$ is poined in $p$.} In case $tfiber(\calO(\underline{n}\setminus\bullet))$ is $(n-3)$-connected, the fiber is always $\Omega^{n-2} tfiber(\calO(\underline{n}\setminus\bullet))$.
\end{proposition}

\subsection{Bimodule {\sc$B\widetilde{\Pentagon}$}}\label{ss:bimod_B_tilde_Pentagon}
In this section we construct an $\Assoc$ bimodule $B\widetilde{\Pentagon}$ which is constructed from $\widetilde{\Pentagon}$ in the same way as $B\Pentagon$ is constructed from $\Pentagon$. It will have the  property that it is a cofibrant replacement of $\Assoc$ as a bimodule over itself, and also one will have a natural homotopy equivalence $\tilde\zeta\colon\Omega\,\Operad(\widetilde{\Pentagon},\calO)\to\underset{\Assoc}{\Bimod}(B\widetilde{\Pentagon},\calO)$ which will be constructed in the next section.

Let $\tilde\Psi_n$ be the category of planar trees with $n$ leaves and internal vertices of any valence. The root is supposed to have valence~1. Morphisms in $\tilde\Psi_n$ are given by contraction of edges in the trees.

Given any non-$\Sigma$ operad $\calO$ one can define a functor
$$
\lambda_{\calO}\colon\tilde\Psi_n\to\Top
$$
that assigns $\prod_{b\in B(T)}\calO(|b|)$ to any tree $T\in\Psi_n$, where $B(T)$ is the set of inner vertices of $T$, and $|b|$ denotes the number of outgoing edges.
We will be actually interested at the moment only in the functor $\lambda_{\widetilde{\pentagon}}$

We also define a contravariant functor
$$
\chi_\blacktriangle\colon\tilde\Psi_n\to\Top
$$
in exactly the same way as we do it in Subsection~\ref{ss:PB}. To be precise $\chi_\blacktriangle$ assigns to any planar tree the space of order preserving maps from the set of its beads to~$[0,1]$.

Then we define
$$
B\widetilde{\Pentagon}(n):=\left.\left(\lambda_{\widetilde{\pentagon}}\otimes_{\tilde\Psi_n}\chi_\blacktriangle\right)\right/\sim,
$$
where the equivalence relations are the same as in Subsection~\ref{ss:PB} plus one additional relation that takes into account the identity element $id\in\widetilde{\Pentagon}(1)$. Namely, if it happens that the identity element $id$ is assigned to a   vertex $b$ of a tree $T\in\tilde\Psi_n$ (which in particular means that $|b|=1$), then this element in $\lambda_{\widetilde{\pentagon}}\otimes_{\tilde\Psi_n}\chi_\blacktriangle$ is equivalent to the element whose underlying tree is obtained from $T$ by replacing $b$ and two adjacent edges by a single edge. All the other vertices and labels stay unchanged:
$$
\psfrag{id}[0][0][1][0]{$(id,t)$}
\psfrag{sim}[0][0][1][0]{$\sim$}
\includegraphics[width=2cm]{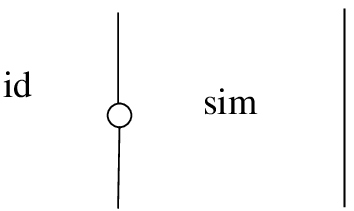}
$$

For a tree $\widetilde{T}\in\tilde\Psi_n$ the image of $\lambda_{\widetilde{\pentagon}}(\widetilde{T})\times \chi_\blacktriangle(\widetilde{T})$ in $B\widetilde{\Pentagon}(n)$ will be denoted by $B\widetilde{\Pentagon}(\widetilde{T})$.

The sequence of spaces $B\widetilde{\Pentagon}$ has a natural structure of a bimodule over $\Assoc$. This structure is absolutely analogous to the construction given in Subsection~\ref{ss:PB}.

We will define a filtration in $B\widetilde{\Pentagon}$ of $\Assoc$ bimodules:

$$
B\widetilde{\Pentagon}_0\subset B\widetilde{\Pentagon}_1\subset B\widetilde{\Pentagon}_2\subset B\widetilde{\Pentagon}_3\subset\ldots\subset B\widetilde{\Pentagon}
\eqno(\numb)\label{eq:filt_tilde_B_Pent}
$$
and a subfiltration
$$
B\widetilde{\Pentagon}_{n-1}=B\widetilde{\Pentagon}_{n-1,-1}\subset B\widetilde{\Pentagon}_{n-1,0}\subset B\widetilde{\Pentagon}_{n-1,1}\subset B\widetilde{\Pentagon}_{n-1,2}\subset\ldots\subset B\widetilde{\Pentagon}_{n-1,n}=B\widetilde{\Pentagon}_n
\eqno(\numb)\label{eq:subfilt_tilde_B_Pent}
$$
in each inclusion $B\widetilde{\Pentagon}_{n-1}\subset B\widetilde{\Pentagon}_n$. First notice that the cell decomposition of $\widetilde{\Pentagon}$ induces a natural cell structure in $B\widetilde{\Pentagon}$. Its description is completely analogous to the description of the induced cell structure in $B\Pentagon$, see Subsection~\ref{ss:PB}. This induced cell structure in $B\widetilde{\Pentagon}$ has a natural coarsening which we will describe below. Combinatorially the latter structure is obtained by attaching cells of the same dimension and in the same degrees as in the case of $\widetilde{\square}$, see Section~\ref{ss:tilde_square}. Which means in particular that $B\widetilde{\Pentagon}_{n-1,i}$ is obtained from $B\widetilde{\Pentagon}_{n-1,i-1}$ by a free attachment of $n\choose i$ cells of dimension $n-1+i$ in degree $n-i$.

To construct such cell structure in $B\widetilde{\Pentagon}$ for each $n\geq 0$ define a map
$$
\sigma_{n,B\widetilde{\pentagon}}\colon
\Bigl(\coprod_{m\geq 0}\,\,\,\coprod_{\alpha\colon\underline{m}\hookrightarrow\underline{n+m}} B\Pentagon(n+m)\times [0,1]^{\alpha(\underline{m})}\Bigr)
\to
B\widetilde{\Pentagon}(n).
\eqno(\numb)\label{eq:tilde_wsq_disj}
$$
The filtration term $B\widetilde{\Pentagon}_{N-1,i}$ will be defined as a subbimodule of $B\widetilde{\Pentagon}$ generated by the images (under the above map) of the interiors  of $B\Pentagon(n+m)\times [0,1]^{\alpha(\underline{m})}$ with $n+m\leq N-1$, or $n+m=N$ and $m\leq i$. The subbimodule $B\widetilde{\Pentagon}_N$ is defined as $B\widetilde{\Pentagon}_{N-1,N}=B\widetilde{\Pentagon}_{N,-1}$.

Let ${\mathbf X}=(X,\bar\tau)_{\alpha}\in B\Pentagon(n+m)\times [0,1]^{\alpha(\underline{m})}$ be an element for which we need to define
$\sigma_{n,B\widetilde{\pentagon}}({\mathbf X})$, where $X=(x_b,t_b)_{b\in B(T)}\in B\Pentagon(T)$ with $T$ being a tree in $\Psi_{n+m}$, and $x_b\in\Pentagon(|b|)$,
$t_b\in[0,1]$ (one views $t_\bullet$ as an order preserving function from the set of beads $B(T)$ to $[0,1]$), one also has
$\bar\tau=(\tau_i)_{i\in\alpha(\underline{m})}$. The element $\sigma_{n,B\widetilde{\pentagon}}({\mathbf X})$ will lie in $B\widetilde{\Pentagon}(\widetilde{T})\subset B\widetilde{\Pentagon}(n)$, where
$\widetilde{T}$ is obtained from $T$ by a procedure specified below. Briefly speaking to obtain $\widetilde{T}$ each leaf of $T$ labeled by $\alpha(i)$ will be either cut off (in case $\tau_{\alpha(i)}$ is sufficiently small) or replaced by a univalent bead (if $\tau_{\alpha(i)}$ is sufficiently big). First we define a function $U_T\colon \underline{m+n}\to B(T)$ that assigns to a leaf of $T$ the bead to which this
leaf is attached. Then  we define a set
$$
S_{\mathbf X}:=\{s\in\alpha(\underline{m})\, |\, \tau_s< t_{U_T(s)}\}.
$$
 The tree $\widetilde{T}\in\tilde\Psi_n$ is obtained from $T$ by cutting off all the leaves that are from the set $S_{\mathbf X}$ and by replacing each leaf from $\alpha(\underline{m})\setminus S_{\mathbf X}$ by a bead. Thus the set of beads of $\widetilde{T}$ is
$$
B(\widetilde{T})=B(T)\cup (\alpha(\underline{m})\setminus S_{\mathbf X}).
$$
In coordinates $\sigma_{n,B\widetilde{\pentagon}}({\mathbf X})$ is written as
$$
\sigma_{n,B\widetilde{\pentagon}}({\mathbf X})=(x_{\widetilde{b}},t_{\widetilde{b}})_{\widetilde{b}\in B(\widetilde{T})}=
\left(\sigma_{|b|,\widetilde{\pentagon}}\left(\left(x_b,(\tau_i/t_b)_{i\in U_T^{-1}(b)\cap S_{\mathbf X}}\right)_{\rho_b}\right),t_b\right)_{b\in B(T)}\times (\,\raisebox{-4pt}{\includegraphics[width=.13cm]{b0.eps}}\,,\tau_i)_{i\in
\alpha(\underline{m})\setminus S_{\mathbf X}}.
$$
In the above $\rho_b$ stands for the inclusion $\rho_b\colon \underline{|U_T^{-1}(b)\cap S_{\mathbf X}|}\hookrightarrow \underline{|b|}$, where the left-hand side $\underline{|U_T^{-1}(b)\cap S_{\mathbf X}|}$ stands for the (ordered) set of leaves in $S_{\mathbf X}$ that are adjacent to the vertex $b$, and the right-hand side $\underline{|b|}$ stands for the (ordered) set of outgoing edges from $b$. To recall the map $\sigma_{n,\widetilde{\pentagon}}$ is from~\eqref{eq:tilde_pent}.

The proof of the following result is analogous to that of Proposition~\ref{p:operad_cof_tilde}:

\begin{proposition}\label{p:cof_B_tilde_pentag}
(1) The $\Assoc$ bimodule $B\widetilde{\Pentagon}$ is a cofibrant model for $\Assoc$.

(2) The $n$-truncated $\Assoc$ bimodule $B\widetilde{\Pentagon}_n|_n$ is a cofibrant model for $\Assoc|_n$ .
\end{proposition}

%

%
%

\subsection{Tower associated to {\sc $B\widetilde{\Pentagon}$}, $B$-construction, and second delooping}\label{ss:tower_B_tilde_pent}

Given an $\Assoc$ bimodule $\calO$, filtration~\eqref{eq:filt_tilde_B_Pent} induces a tower $T_\bullet^{B\widetilde{\pentagon}}(\calO)$ of fibrations, with $n$-th stage
$$
T_n^{B\widetilde{\pentagon}}(\calO):=\underset{\Assoc}{\Bimod}(B\widetilde{\Pentagon}_n,\calO)=
\underset{\Assoc}{\Bimod}{}_n(B\widetilde{\Pentagon}_n|_n,\calO).
$$

\begin{proposition}\label{p:equiv_tild_b}
One has that $\hoSquare$ is homotopy equivalent to $B\hoPentagon$ as a filtered bimodule over $\Assoc$. As a consequence  the towers $T_\bullet^{\hoSquare}(\calO)$ and $T_\bullet^{B\widetilde{\pentagon}}(\calO)$ are also homotopy equivalent.
\end{proposition}

\begin{proof}
This follows immediately from Propositions~\ref{p:cof_bimod_tilde},~\ref{p:cof_B_tilde_pentag}, and Lemma~\ref{l:trunc_maps}.
\end{proof}

It follows from Propositions~\ref{p:cof_bimod_tilde},~\ref{p:cof_B_tilde_pentag}, and Lemma~\ref{l:cof_maps} that $\hoSquare$ is homotopy equivalent to $B\hoPentagon$ as a filtered bimodule over $\Assoc$. As a consequence the towers  $T_\bullet^{\widetilde{\square}}(\calO)$ and $T_\bullet^{B\widetilde{\pentagon}}(\calO)$ are homotopy equivalent. In the sequel we will need to use the spaces
$$
T_{n,i}^{B\widetilde{\pentagon}}(\calO):=\underset{\Assoc}{\Bimod}(B\widetilde{\Pentagon}_{n,i},\calO).
$$

It is easy to see that given any operad $\calP$, one can define an $\Assoc$ bimodule $B\calP$ in the same way as one constructed $B\widetilde{\Pentagon}$ from
$\hoPentagon$. Moreover if $\calO$ is a pointed operad, one can define an evaluation map
$$
\tilde{\zeta}_\calP\colon \Omega\, \Operad (\calP,\calO)\to\underset{\Assoc}{\Bimod}(B\calP,\calO).
$$
Due to the inclusion of $\Assoc$-bimodules $B\widetilde{\Pentagon}_n\subset B(\widetilde{\Pentagon}_n)$, $B\widetilde{\Pentagon}_{n,i}\subset B(\widetilde{\Pentagon}_{n,i})$, the maps $\tilde\zeta_{\widetilde{\pentagon}_n}$, $\tilde\zeta_{\widetilde{\pentagon}_{n,i}}$ can be composed with the restriction maps to get maps
$$
\tilde\zeta_n\colon \Omega T_n^{\widetilde{\pentagon}}(\calO)=\Omega\, \Operad(\widetilde{\Pentagon}_n,\calO)\to \underset{\Assoc}{\Bimod}(B\widetilde{\Pentagon}_n,\calO)=
T_n^{B\widetilde{\pentagon}}(\calO);
$$
$$
\tilde\zeta_{n,i}\colon \Omega T_{n,i}^{\widetilde{\pentagon}}(\calO)=\Omega\,\Operad(\widetilde{\Pentagon}_{n,i},\calO)\to \underset{\Assoc}{\Bimod}(B\widetilde{\Pentagon}_{n,i},\calO)=
T_{n,i}^{B\widetilde{\pentagon}}(\calO).
$$

\begin{theorem}\label{t:deloop2_stages_tilde}
Each map $\tilde\zeta_n\colon \Omega\, T_n^{\widetilde{\pentagon}}(\calO)\to T_n^{B\widetilde{\pentagon}}(\calO)$ is a homotopy equivalence for any $n\geq 1$ and any pointed operad $\calO$ with $\calO(1)\simeq *$.
\end{theorem}

By a \textit{pointed} operad we mean an operad in based spaces.

\begin{proof}
 One has that $\tilde\zeta_1$ sends  $\Omega\,T_1^{\widetilde{\pentagon}}(\calO)=\Omega\calO(0)$ homeomorphically to the preimage of $\id\in\calO(1)$ under the natural fibration map $T_1^{B\widetilde{\pentagon}}(\calO)\to\calO(1)$. Thus $\tilde\zeta_1$ is also a homotopy equivalence.
For the other stages we will use the induction.

We will prove that in the settings of the theorem all the maps $\tilde\zeta_{n,i}$ are homotopy equivalences. The result will follow since $\tilde\zeta_n=\tilde\zeta_{n-1,n}$.

 The proof is very much similar to that of Theorem~\ref{t:deloop1_stages_tilde}. We will use again the notation  $\tilde\zeta_{n-1,-1}:=\tilde\zeta_{n-1}$, $T_{n-1,-1}^{\widetilde{\pentagon}}(\calO):=T_{n-1}^{\widetilde{\pentagon}}(\calO)$, $\widetilde{\Pentagon}_{n-1,-1}:=\widetilde{\Pentagon}_{n-1}$, and so on.

Assuming that $\tilde\zeta_{n-1,i-1}$ is a homotopy equivalence we will prove that $\tilde\zeta_{n-1,i}$ has the same property. Again by the induction it is enough to show that $\tilde\zeta_{n-1,i}$ sends the fibers of $\Omega T_{n-1,i}^{\widetilde{\pentagon}}(\calO)\to\Omega T_{n-1,i-1}^{\widetilde{\pentagon}}(\calO)$ to the fibers of $T_{n-1,i}^{B\widetilde{\pentagon}}(\calO)\to T_{n-1,i-1}^{B\widetilde{\pentagon}}(\calO)$ by a homotopy equivalence.

  To recall $B\widetilde{\Pentagon}_{n-1,i}$ is obtained from $B\widetilde{\Pentagon}_{n-1,i-1}$ by a free attachment of $n\choose i$ cells in degree $n-i$. These cells are exactly the interior of $\coprod_{\alpha\colon\underline{i}\hookrightarrow\underline{n}}B\Pentagon(n)\times[0,1]^{\alpha(\underline{i})}$ in~\eqref{eq_tilde_wsq_disj}. Consider
$$
B\widetilde{\Pentagon}_{n-1,i-1/2}:=B(\widetilde{\Pentagon}_{n-1,i-1})\cap B\widetilde{\Pentagon}_{n-1,i}.
$$
One can see that $B\widetilde{\Pentagon}_{n-1,i-1/2}$ is obtained from $B\widetilde{\Pentagon}_{n-1,i-1}$ by a free attachment of $n\choose i$ punctured discs
$\coprod_{\alpha\colon\underline{i}\hookrightarrow\underline{n}}\left(B\Pentagon(n)\times [0,1]^{\alpha(\underline{i})}\right)\setminus Int(D_{\alpha})$.
 Where $D_{\alpha}$ is a closed subdisc of the same dimension in (the subdisc) $B\Pentagon(C_n)\times [0,1]^{\alpha(\underline{i})}$ that consists of points
 $\left((x,t),(\tau_i)_{i\in\alpha(\underline{m})}\right)$,  with $x\in\Pentagon(n)$, $t\in [0,1]$, $(\tau_j)_{j\in\alpha(\underline{i})}\in
 [0,t]^{\alpha(\underline{i})}$.

On the other hand using Lemma~\ref{l:extension_space} it is easy to see that given $\tilde h\in \Omega T_{n-1,i-1}^{\widetilde{\pentagon}}(\calO)$ the fiber over $\tilde h$   in $\Omega T_{n-1,i}^{\widetilde{\pentagon}}(\calO)$ is
homeomorphic to the space of maps $\underset{\Assoc}{\Bimod}\, _{\tilde\zeta_{n-1}^{\tilde h}}\left((B\widetilde{\Pentagon}_{n-1,i},B\widetilde{\Pentagon}_{n-1,i-1/2}),\calO\right)$. While the fiber over its image $\tilde\zeta_{n-1,i-1}^{\tilde h}\in
T_{n-1,i-1}^{B\widetilde{\pentagon}}(\calO)$ is the space $\underset{\Assoc}{\Bimod}\, _{\tilde\zeta_{n-1}^{\tilde h}}\left((B\widetilde{\Pentagon}_{n-1,i},B\widetilde{\Pentagon}_{n-1,i-1}),\calO\right)$. It follows from Lemma~\ref{l:fiber_equiv}, that the inclusion
$$
\underset{\Assoc}{\Bimod}\, _{\tilde\zeta_{n-1}^{\tilde h}}\left((B\widetilde{\Pentagon}_{n-1,i},B\widetilde{\Pentagon}_{n-1,i-1/2}),\calO\right)\to
\underset{\Assoc}{\Bimod}\, _{\tilde\zeta_{n-1}^{\tilde h}}\left((B\widetilde{\Pentagon}_{n-1,i},B\widetilde{\Pentagon}_{n-1,i-1}),\calO\right)
$$
is a homotopy equivalence.
\end{proof}

\begin{theorem}\label{t:second_deloop_tilde}
In the settings of Theorem~\ref{t:deloop2_stages_tilde} the induced map of limits
$$
\Omega\,\Operad(\Pentagon,\calO)\stackrel{\tilde\zeta_\infty}{\longrightarrow}\underset{\Assoc}{\Bimod}(B\hoPentagon,\calO).
\eqno(\numb)\label{eq:second_deloop_xi_tilde}
$$
is a homotopy equivalence.
\end{theorem}
\begin{proof}
The proof is similar to that of Theorem~\ref{t:first_deloop}.
\end{proof}

\subsection*{Acknowledgement} The author is grateful to M. Kontsevich, P.~Lambrechts, P.~Salvatore, and D.~Stanley discussions with whom deeply influenced this work. He also thanks Max Planck Institute for Mathematics, Universit\`a  di Roma Tor Vergata, Universit\'e Catholique de Louvain, and University of Virginia for hospitality.


\end{document}